\documentclass[a4paper,oneside,english,reqno]{amsart}
\usepackage[utf8]{inputenx}
\usepackage[russian,english]{babel}

\usepackage{amsmath,amsfonts,amssymb,mathrsfs,amscd,comment,latexsym,bbold}
\usepackage{amsthm, mathtools, stmaryrd}
\usepackage{textcomp, url, enumerate, latexsym, graphicx, varioref, hyperref, multirow, layout, siunitx, booktabs}
\usepackage{pdflscape}
\usepackage{verbatim}

\usepackage[matrix,arrow,curve]{xy}
\usepackage{tikz-cd}

\usepackage{enumitem}

\usepackage{csquotes}

\usepackage[top = 3cm, bottom = 4cm, left = 3cm, right = 3cm]{geometry}

\usepackage[pagewise,displaymath,mathlines]{lineno}
\usepackage{microtype}

\RequirePackage[backend    = biber,
                sortcites  = true,
                giveninits = true,
                doi        = false,
                isbn       = false,
                url        = false,
                maxnames   = 6,
                style      = numeric]{biblatex}
\DeclareFieldFormat*{title}{#1}
\usepackage{cleveref}
\usepackage{amsthm}
\usepackage{biblatex}
\renewbibmacro{in:}{}

\numberwithin{equation}{section}
\theoremstyle{plain}
\newtheorem{theorem}[equation]{Theorem}
\newtheorem*{theorem*}{Theorem}

\newtheorem{corollary}[equation]{Corollary}
\newtheorem*{corollary*}{Corollary}
\newtheorem{proposition}[equation]{Proposition}
\newtheorem*{nonumproposition}{Proposition}

\newtheorem{lemma}[equation]{Lemma}

\theoremstyle{definition}

\newtheorem{definition}[equation]{Definition}

\newtheorem{example}[equation]{Example}

\theoremstyle{remark}
\newtheorem{remark}[equation]{Remark}

\theoremstyle{plain} % just in case the style had changed
\newcommand{\thistheoremname}{}
\newtheorem*{genericthm}{\thistheoremname}

\newtheoremstyle{named}{}{}{\itshape}{}{\bfseries}{.}{.5em}{\thmnote{#3's }#1}
\DeclareMathOperator{\coker}{coker}
\DeclareMathOperator{\Coker}{Coker}

\DeclareMathOperator{\Cofib}{Cofib}
\DeclareMathOperator{\Ker}{Ker}
\DeclareMathOperator{\codim}{codim}
\DeclareMathOperator{\chark}{char}

\DeclareMathOperator{\coCone}{coCone}

\newcommand{\Fr}{\mathrm{Fr}}

\newcommand{\ZF}{{\mathbb Z\mathrm{F}}}
\newcommand{\ovZF}{\overline{\mathbb Z\mathrm{F}}}

\newcommand{\A}{\mathbb A}

\newcommand{\fr}{\mathrm{fr}}
\newcommand{\nis}{\mathrm{Nis}}

\newcommand{\Ab}{\mathrm{Ab}}

\newcommand{\PP}{\mathbb P}

\newcommand{\GL}{\mathrm{GL}}

\newcommand{\Sm}{\mathrm{Sm}}

\newcommand{\Aff}{\mathrm{Aff}}

\newcommand{\SmAff}{\mathrm{SmAff}}

\newcommand{\Et}{\mathrm{Et}}
\newcommand{\Sch}{\mathrm{Sch}}

\newcommand{\zar}{\mathrm{zar}}

\newcommand{\Nis}{\mathrm{Nis}}

\newcommand{\Fin}{\mathrm{Fin}}

\newcommand{\ZShS}{\mathbb Z\mathrm{Sh}(S)}
\newcommand{\ZShfrS}{\mathbb Z\mathrm{Sh}(\ZF_*(S))}
\newcommand{\Ext}{\mathrm{Ext}}
\newcommand{\Pre}{\mathrm{Pre}}
\newcommand{\nrFr}{\mathrm{Fr}^\mathrm{1th}}
\newcommand{\nrFrpair}{\mathrm{Fr}^{\mathrm{1th},\pair}}
\newcommand{\ZPre}{\mathbb Z\mathrm{Pre}}

\newcommand{\SH}{\mathbf{SH}}
\newcommand{\DM}{\mathbf{DM}}

\newcommand{\codimcalVWprimer}{W^\prime\in\calV^{[r]}}
\newcommand{\codimcalVWprimermo}{W^\prime\in\calV^{[r-1]}}

\newcommand{\Xx}{{X_{x}}}
\newcommand{\Ueta}{{U_{\underline{\eta}}}}
\newcommand{\calV}{\Ueta}
\newcommand{\calO}{\mathcal O}
\newcommand{\calL}{\mathcal L}
\newcommand{\calI}{\mathcal I}
\newcommand{\ovZ}{\overline{Z}}
\newcommand{\ovcalZ}{\overline{\mathcal Z}}
\newcommand{\ovS}{\overline{S}}
\newcommand{\ovX}{\overline{X}}

\newcommand{\calZ}{\mathcal{Z}}

\newcommand{\ovC}{\overline{C}}
\newcommand{\val}{\mathrm{val}}
\newcommand{\id}{\mathrm{id}}
\newcommand{\can}{\mathrm{can}}
\newcommand{\pr}{\mathrm{pr}}

\newcommand{\Xhinf}{X_{\overline{\infty}}}
\newcommand{\Xhpinf}{X^\prime_{\overline{\infty}}}
\newcommand{\pri}{\prime}

\newcommand{\pair}{\mathrm{pair}}
\newcommand{\Set}{\mathrm{Set}}

\newcommand{\unH}{\underline{H}}

\newcommand{\Spec}{\mathrm{Spec}\;}

\newcommand{\Sh}{\mathrm{Sh}}

\newcommand{\XTheta}{{X^{(0)}}}%{\Theta}

\newcommand{\Corkdash}{$\mathrm{Cor}(k)$}%{$\mathrm{Cor}$}%{$\mathrm{Cor}(k)$-}

\author{Andrei Druzhinin}
\address{Andrei Druzhinin, \\
Chebyshev Laboratory, St. Petersburg State University,\\
a.druzhinin@spbu.ru,
and \\
St. Petersburg Department of Steklov Mathematical Institute 
of Russian Academy of Sciences, Russia\\
andrei.druzh@gmail.com
}

\begin{document}
% article 

\selectlanguage{english}
% \title[SHI theorem via compactified homotopies and local smooth schemes over DVR]{
% Strict %$\A^1$-
% homotopy invariance theorem via compactified homotopies and correspondences, 
% and generic fibre of local essentially smooth schemes over DVR.
% }
% \title[SHI theorem via compactified homotopies and local smooth schemes over DVR]{
% Strict homotopy invariance via compactified homotopies, 
% and fibres of essentially smooth schemes over
% one-dimensional base schemes%
% .
% % DVR
% }
% \title[SHI theorem via compactified homotopies and smooth schemes fibres]{
\title[SHI theorem via compactified homotopies]{
Strict homotopy invariance via compactified homotopies and correspondences, 
and fibres of essentially smooth schemes over
one-dimensional base schemes.
% DVR.
}
% \title[SHI theorem via compactified homotopies and generic fibres over DVR]{
% Strict $\A^1$-homotopy invariance theorem via compactified homotopies and correspondences, and generic fibres of local essentially smooth schemes over DVR.
% }
% \title[Strict $\A^1$-homotopy invariance theorem via compactified homotopies]{Strict $\A^1$-homotopy invariance theorem via compactified homotopies and correspondences, and generic fibres of local essentially smooth .}

\begin{abstract}

We develop the technique of compactified correspondences and homotopies over 
one-dimensional base schemes,
and
illuminate 
the perfectness 
and 
the inverting of characteristic 
assumptions 
from 
the celebrating Voevodsky's strict homotopy invariance theorem and 
its framed correspondences generalisation over an arbitrary base field.
The assumption in this crucial theorem for Voevodsky's motives theory %over a field 
was kept from the origins of the study, 
and came later into more modern theory of framed motives by Garkusha-Panin.
Applying 
the technique,
we obtain
also
analogs of 
Gersten and Nisnevich conjectures for 
Cousin complexes of generalised motivic cohomotopies over a field,
and acyclicity of Cousin complexes on generic fibres of essentially smooth local schemes over 
one-dimensional base schemes.
% DVR. %one-dimensional local base schemes.
\end{abstract}

\maketitle
\tableofcontents
\section{Introduction. }

The philosophy of motivic $\A^1$-homotopy theory 
\cite{Morel-Voevodsky}
suggests to consider 
the affine line $\A^1_k\in\Sm_k$ over a base field $k$ by analogy 
to 
the topological interval $[0,1]\in \mathrm{Top}$
in algebraic topology.
The $p$-torsion, where $p$ % = \operatorname{char} k
is the characteristic of the base field $k$, 
is a mysterious part of the study. %in the $\A^1$-homotopy theory;
%todo:changed.
Because of finite \'etale coverings of degrees $p^l$,
in distinct to $[0,1]\in \mathrm{Top}$, the affine line $\A^1_k$ is not one-connected.
%Е% In distinct to $[0,1]\in \mathrm{Top}$, the affine line $\A^1_k$ is not one-connected. 
%
% In distinct to 
% $[0,1]\in \mathrm{Top}$,
% % the topological interval $[0,1]$, 
% the affine line $\A^1_k$ 
% % $\A^1_k\in\Sm_k$
% % is not one-connected. 
% is not one-connected, because of finite \'etale coverings of degrees $p^l$.
% because of finite \'etale coverings of degrees $p^l$,
% and at the same time, 
% somehow 
% for looking formally independent reasons
At the same time, 
%todo %edited
for reasons that look formally independent,
various technical obstructions 
% for a set of reasons which looks formally independent,
% often 
appear inside proofs of 
% remarkable 
complicated structural results on the motivic categories.
% New techniques or even ideology
% are needed to study the borderline of theories 
% that hold 
% with $\mathbb Z$-coefficients over any field.
% A well know example
% is
% Voevodsky's motives theory \cite{Voe-hty-inv,Voe-motives,Voe-cancel,cycles-book} on $\DM(k)$.
% A
%todo%edited
One well-know example
is
the theory \cite{Voe-motives,cycles-book} %Voe-hty-inv, Voe-cancel,
% on 
regarding 
Voevodsky's motives category $\DM(k)$.
% that 
% requires either perfectness assumption on $k$ or inverting of $p$ in the coefficients of the category $\DM(k)$.
%
% The proof of the milestone technical result, 
% namely, the strict homotopy invariance theorem for presheaves with transfers \cite[Theorem 5.6]{Voe-hty-inv},
% % accumulated techniques of 
% accumulated ideas of earlier 
% works on injectivity and purity properties and Gersten conjecture
% % , see
% \cite{BlochOgus,zbMATH05778592}
% % ,
% % combining them with new developed transfers techniques.
%
% The proof of the milestone technical result, 
% namely, the strict homotopy invariance theorem 
% % for presheaves with transfers 
% \cite[Theorem 5.6]{Voe-hty-inv},
The proof of the strict homotopy invariance theorem 
\cite[Theorem 5.6]{Voe-hty-inv},
which is a milestone technical result, 
was based on 
an upgrade of powerful algebro geometric ideas 
from
% oriented to 
% regarding
% injectivity, purity and 
% Gersten
% conjecture 
proofs of
Gersten
conjecture 
% from 
in
\cite{BlochOgus,zbMATH05778592}
combined with a new %useful 
instrument of finite correspondences. 
% made an upgrade of 
% % combined 
% powerful
% algebro geometric
% ideas 
% oriented to
% % regarding
% injectivity, purity and Gersten conjecture
% % injectivity and purity properties and Gersten conjecture
% from 
% % earlier 
% % works 
% \cite{BlochOgus,zbMATH05778592}
% combining
% them
% %additionally
% with 
% a new 
% useful
% % the 
% instrument of finite correspondences.
%
% combining them with 
% the %
% % new 
% instrument of finite correspondences.
% regarding itself,
% and provided new transfers techniques regarding itself,
% and provided new 
% % powerful 
% techniques itself,
Nevertheless, % at the same time,
% it
% the %is 
the %is 
approach % argument 
required 
either perfectness assumption on $k$, 
or inverting of $p$ %the characteristic 
in the coefficient ring, 
see \Cref{sect:Prehistory}.
% see %historical 
% a %brief 
% review %discussion %etails 
% in \Cref{sect:Prehistory}.
% \Cref{sect:Prehistory}.
% for a review.
Our article develops a modern technique % instrument 
of 
% compactified %higher-dimensional 
% $d$-dimensional correspondences and 
compactified homotopies
and 
compactified %higher-dimensional 
$d$-dimensional correspondences
over 
one-dimensional base schemes.
% DVR
This allows to illuminate the assumption, 
% on the base field $k$ in the above result,
because 
% treati
the obstruction was hidden, 
from some viewpoint, 
at 
% relative 
% $X^h_x$-infinity %$\infty$ 
infinity %$\infty$ 
 % over $X^h_x$
of $\A^1_{X^h_x}$, $X\in\Sm_k$, $x\in X$.
% See \Cref{sect:KeyConstrunctionandApplications,sect:ProofStrategy}
% for details and other applications.
% of $\A^1_{X^h_x}$
%  over $X^h_x$, $X\in\Sm_k$, $x\in X$.
% of $\A^1_{X^h_x}$.
% and prove the above result over an arbitrary base field. 
% from this result.
% one milestone technical result, 
% namely, the strict homotopy invariance theorem for presheaves with transfers \cite[Theorem 5.6]{Voe-hty-inv}. 
% Covering %moreover 
% % the generalisation of the latter theorem 
% the result for framed transfers \cite[Theorem 1.1]{hty-inv} as well, 
% % our article 
% % Covering framed transfers as well, see \cite[Theorem 1.1]{hty-inv}, 
% we automatically extend the generality of all results of Garkusha-Panin's framed motives theory \cite{Framed}, 
% % that applies to the stable motivic homotopy category $\SH(k)$ \cite{Voe98,Jardine-spt,morel-trieste,mot-functors}.
% which %extends generality of %applies 
% allows to apply 
% the approaches of Voevodsky's motives theory on $\DM(k)$ \cite{Voe-motives}
% to the stable motivic homotopy category $\SH(k)$ \cite{Voe98,Jardine-spt,morel-trieste,mot-functors}.
Covering framed correspondences context, 
% our result generalises 
we have generalised \cite[Theorem 1.1]{hty-inv} 
and immediately extended %s 
% the generality of 
all results of 
Garkusha-Panin's framed motives theory \cite{Framed}, 
which 
applies
% lifts 
the approach of Voevodsky's motives theory on $\DM(k)$ %\cite{Voe-motives}
to the stable motivic homotopy category $\SH(k)$ \cite{Morel-Voevodsky,mot-functors,Jardine-spt,Voe98,morel-trieste}.
See \Cref{sect:KeyConstrunctionandApplications,sect:ProofStrategy}
for details and other applications.

\begin{theorem*}[\Cref{th:strhominv}]%cor:strhominv%, \Cref{cor:HniscongHzar}, \Cref{cor:ZarStrictHomotopyInvariance}
Let $F$ be an $\A^1$-invariant quasi-stable framed linear presheaf over a field $k$, and $F_\nis$ denote the associated Nisnevich sheaf, and similarly $F_\zar$ for Zariski topology.
Then for any $X\in \Sm_k$, the canonical projection induces isomorphisms 
on %of 
Nisnevich cohomologies
\begin{equation}\label{eq:strA1hominvNis}H^n_\nis(\A^1_k\times X,F_\nis)\cong H^n_\nis(X,F_\nis),\quad n\in \mathbb Z\end{equation}
and similarly for Zariski cohomologies,
% \begin{equation}\label{eq:HzarA1timesXFzarcongHzarXFzar}H^n_\zar(\A^1_k\times X,F_\zar)\cong H^n_\zar(X,F_\zar),\end{equation}
and %even more, 
for each $n\in\mathbb Z$, there are canonical isomorphisms
% $H^n_\zar(X,F_\zar)\cong H^n_\nis(X,F_\nis)$ for each $n\in\mathbb Z$.
\begin{equation}\label{eq:HniscongHzar}H^n_\zar(X,F_\zar)\cong H^n_\nis(X,F_\nis),\quad n\in \mathbb Z.\end{equation}
% In particular, the claim holds for presheaves with classical Voevodsky's transfers with respect to the finite correspondences $\mathrm{Cor}(k)$.
\end{theorem*}

\subsection{Strict homotopy invariance}\label{sect:Prehistory}

In the topological setting,
homotopy classes of 
% some things 
%objects
something 
% always 
tautologically form a discrete set. %, 
% which represents a homotopy invariant presheaf, 
% on $\Top$
% sounds as a tautology in the topological setting.
%
% At the same time,
% sounds as a tautology in the .
% while
In the motivic homotopy theory,
this property
is called %stro
strict homotopy invariance
and relates to 
the fundamental complexity 
formed by the combination of $\A^1$-homotopy and Nisnevich localisations
in the construction of the $\A^1$-motivic localisation. 
% %via the transfers structure 
% in the constructions of $\HH(k)$, $\SH(k)$, and $\DM(k)$.
%
% in the this setting.
% \cite{Voe-motives,Voe-cancel,MVW,Morel-Voevodsky,Jardine-spt,morel-trieste,mot-functors}
% Various strict homotopy invariance theorems 
% aim to 
% solve this complexity 
% in the following sense.
% While the general construction of $\A^1$-motivic localisation on the category of a simplicial Nisnevich sheaves in \cite{Morel-Voevodsky} requires infinitely many iterative 
% % procedure
% % of 
% $\A^1$-invariantisations and Nisnevich localisations, 
% %Suslin's 
% % $F\mapsto F(\Delta^\bullet_k\times_k -)$
% % $\A^1$-homotopy and Nisnevich localisations, 
% the strict homotopy invariance results 
% % allow to 
% reduce it 
% to single $\A^1$-invariantisation
% under curtain assumptions on the sheaf. 
Under curtain assumptions on 
% the class of 
% the sheaves or 
the motivic spaces
various strict homotopy invariance theorems 
reduce infinite iterative 
composition %procedure 
of 
$\A^1$-invariantisations and Nisnevich 
localisations 
% in the construction of $\A^1$-motivic localisation
to the single iteration.
% $\A^1$-invariantisation 
% % procedure. 
% % This procedure % which 
% % is 
% given by
% % the totoalisation of 
% the bi-simplicial presheaf $F(\Delta^\bullet_k\times_k -)$.
%
% under 
% % some 
% curtain assumptions on $F$.
%
% The strict homotopy invariance for representable theories provided by
% There is for tecn\hnical computing
The theorem discussed in the article
regards %a sheaf 
% $\A^1$-invariant 
structure of 
% sheaves
% % $F$ 
% with 
% % that have %certain additional 
% % structure %basically 
% % called 
transfers, %or framed transfers,
% see \cite{Voe-hty-inv,Voe-notes,hty-inv,Framed},
which 
% Such sheaves 
% have shown themself 
% as 
gives 
a 
universal
computational framework
% form 
% a computational tool 
for 
% precise 
% fibrant replacements with respect to injective
% stable motivic model structures 
motivic localisation functors
and hom-groups,
% mapping groups %sets 
% in the homotopy categories
% \cite{Voe-motives,Framed}
% in a universal sense,
see
\cite{Voe-hty-inv,Voe-notes,hty-inv,Framed}.
\cite[Theorem 5.6]{hty-inv} is the result for $F$ with \Corkdash-transfers over a perfect base field $k$. Suslin proved in 
\cite{Sus-nonperftalk,Sus17nonperfect} the result with $\mathbb Z[1/p]$-coefficients for non-perfect base fields. 
The argument is a reduction to the prefect base field case known early for algebraic K-theory as was
explained in \cite{Sus17nonperfect}.
The framed motives theory by Garkusha and Panin \cite{Framed}
provides a computation for the stable motivic homotopy types in %the category 
$\SH(k)$ % being 
based on the unpublished notes \cite{Voe-motives} by Voevodsky.
Notes \cite{Voe-motives} introduced the framed correspondences as an instrument that would allow to study $\SH(k)$ in parallel to $\DM(k)$ \cite{Voe-motives}.
Large part of %proofs in 
framed motives theory 
\cite{Framed,hty-inv,framed-cancel,ConeTheGNP,surj-etale-exc,Nesh-FrKMW,DrKyllfinFrpi00,ElmantoKhannonperfect}
upgrades the respective results
% ones %in the case of framed correspondences of the arguments originally 
% used 
for \Corkdash-correspondences,
and upgrading the original Voevodsky's reasoning
% argument 
%from \cite{Voe-hty-inv} 
%to the case of framed correspondences
Garkusha and Panin in \cite{hty-inv} covered 
initially the case of an infinite perfect base field of odd characteristic, 
and 
additional arguments in 
% the additional works 
\cite{surj-etale-exc,DrKyllfinFrpi00,five-authors} 
covered the generality of arbitrary perfect fields.
% achieve the same generality as was known for \Corkdash-correspondences.
% and works \cite{surj-etale-exc,DrKyllfinFrpi00,five-authors,ElmantoKhannonperfect,notes,SHetfrk} 
% allowed to achieve the same generality as was known for \Corkdash-correspondences.
Note that the resolution of singularities over $k$ required for the cancellation theorem in \cite{Voe-motives} was illuminated by Voevodsky in \cite{Voe-cancel}, and due to the upgrade of the latter argument for framed correspondences in \cite{framed-cancel},
the strict homotopy invariance theorem \cite[Theorem 1.1]{hty-inv} is the only reason of assumptions on the base field 
for results of \cite{Framed,BigFrmotives,FramedGamma}, and \cite{five-authors}.
% %
% % Summarising early known cases and facts regarding the main theorem above we say that:
%
Here is the summary of early known %cases and 
facts: 
% regarding \Cref{th:strhominv}
(1) $n=0$ by \cite{Voe-hty-inv,hty-inv};
(2) for $F$ being an inverse image of presheaves over perfect fields;
(3) 
for presheaves of $\mathbb Z[1/p]$-modules 
% over a field for characteristic $p$ 
with \Corkdash-transfers \cite{Sus17nonperfect};
% for \Corkdash-correspondences \cite{Sus17nonperfect};
(4) it was known 
the equivalence of the natural isomorphism \eqref{eq:strA1hominvNis} 
with the one for Zariski topology \cite{Sus-nonperftalk,five-authors,SHzar},
% and \eqref{eq:HzarA1timesXFzarcongHzarXFzar}, 
% and the from the result implication to \eqref{eq:HniscongHzar}, 
though 
% the isomorphism \eqref{eq:HniscongHzar} was not.
isomorphisms $H^n_\zar(X,F_\zar)\cong H^n_\nis(X,F_\nis)$ were not.

The general case for anyone Zariski or Nisnevich topology was an entire mystery, and our own opinion on whether the claim holds 
was several times changed during years of intensive work
and periodical looking for a disproof
before the argument was found.
Our proof is elementary. 
The only external 
ingredients are 
\'etale excision and sheafification theorems 
\cite[Theorem 3.10]{surj-etale-exc}, \cite[Proposition 16.2]{hty-inv},

\subsection{Key construction and applications}\label{sect:KeyConstrunctionandApplications} %further %proposition
% From some side, the obstruction was hidden a $\infty$ 
% of $\A^1_{X^h_x}$.
% of relative affine 
% Let us note,
While 
the article subject
% the subject of the article 
belongs to the classical area of $\A^1$-homotopy theory over base fields,
the developed instrument 
is applicable
% applies 
simultaneously %to
in the following two modern 
directions: 
% in motivic homotopy theory:
(1) 
studies 
over 
positive-dimensional
base schemes 
oriented 
on 
a computational results for
% lower negative
the lowest non-trivial
%negative
stable motivic homotopy groups 
and 
cohomologies of Cousin or Gersten complexes,
see \cite{DKO:SHISpecZ,ColumnsCousinlocesssmXoverB}, 
and 
% (1) 
% studies 
% over 
% positive-dimensional
% base schemes 
% oriented 
% on 
% a computation of 
% % lower
% % the 
% stable motivic homotopy groups 
% and 
% % the 
% cohomologies of Cousin or Gersten complexes,
% % over base schemes, 
% see \cite{DKO:SHISpecZ,ColumnsCousinlocesssmXoverB}, 
% and 
(2) 
studies of 
$\square$-homotopy theories and categories 
\cite{logDMk,binda2024logarithmicmotivichomotopytheory,zbMATH07341096,zbMATH07341097,MotiveswithmodulusIII}.
% \cite{MotiveswithmodulusIII}, %MotiveModulePairs %S. Saito, Purity of reciprocity sheaves, Adv. in Math. 366 (2020), 107067,70 pp.
% \cite{binda2024logarithmicmotivichomotopytheory},
% \cite{zbMATH07341096},
% \cite{logDMk}.
% ,
% \cite{zbMATH07654583}.
Note that 
(1)
the perfectness assumption 
excludes 
such natural class of base schemes as 
$\A^1\times B$ for any non-zero characteristic scheme $B$,
% of positive characteristic,
% and
(2)
$\square$-homotopy analog of 
$\DM_{\acute{e}t}(k,\mathbb Z/p)$ is non-trivial,
which additionally increases the interest with respect to the $p$-torsion
in Nisnevich motivic categories. 
To achieve 
\Cref{th:strhominv}
% the above result 
we consider 
$\A^1_k$ 
as a base scheme 
with 
the compacitifcation 
$B=\PP^1_k$ 
and 
the point $z=\infty$
and prove
the following 
% crucial 
injectivity theorem,
which is called Proposition because of the technical form, while it is the main result of the article in the appropriate sense.
%
% To achieve the above theorem 
% prove the following injectivity theorem
% applying
% to the base scheme $B=\PP^1_k$ with the point $z=\infty$,
% and .
% we prove 
% the following injectivity theorem
% and
% apply it
% % in the proof of the above theorem 
% to the base scheme $B=\PP^1_k$ and $z=\infty$.
\begin{nonumproposition}[\protect{\Cref{cor:contractinggenpointhomoveretaSmX}}]
Let $B$ be a one-dimensional %irreducible 
scheme, 
$z\in B^{(1)}$, $\eta\in B^{(0)}$.
% with a closed point $z\in B$ and a generic point $\eta\in B$.
% and $B_{(z)}$ denote the local scheme at $z$.
Let $X\in \Sm_B$, $x\in X\times_B z$.
Denote by $U=\Xx$ %=\Spec \mathcal O_{X,x}
the local scheme at $x$,
and define $\calV=U\times_B \eta$.
Then for any closed subscheme $Z$ in $X_{\underline{\eta}}=X\times_B \eta$, there are 
framed correspondences 
$c\in \ZF_N(\calV\times\A^1, X_{\underline{\eta}})$, $c^\prime\in \ZF_N(\calV\times\A^1, X_{\underline{\eta}}-Z)$
such that \begin{itemize}
\item[(1)] $c \circ i_0=\sigma^N_{\calV}\can_{\underline{\eta}}$, $c\circ i_0 = j\circ c^\prime$,
where $i_0,i_1\colon \calV\to \calV\times\A^1$ are the zero and unit sections, $\can_{\underline{\eta}}\colon \calV\to X_{\underline{\eta}}$ is the canonical morphism, $j\colon X_{\underline{\eta}}-Z\to X_{\underline{\eta}}$ is the open immersion,
\item[(2)] $c^{-1}(Z)$ is finite over $\calV$, see \eqref{sect:corrpreimZ} for $c^{-1}(Z)$. 
\end{itemize}
\end{nonumproposition}
% We called the above result ``injectivity theorem'' because it immediately implies 
% injectivity
% \[E^l(\calV)\to E^l(\calV^{(0)}),\]
% where $\calV^{(0)}$ denote the union of generic points of $\calV$,
% for an $\SH(B)$-representable cohomology theory $E^*$,
% and b
% Then the following result follows immediately
% Part (1)
% of Proposition \ref{cor:contractinggenpointhomoveretaSmX}
% % the proposition above 
% immediately
% implies the following result
% by the argument for the Gersten conjecture in \cite{Voe-hty-inv}.
% See
% \cite{BlochOgus} for the notion of Cousin complexes.
%
% To demonstrate applications of Part (1) of Proposition \ref{cor:contractinggenpointhomoveretaSmX}
% we formulate the %following 
% result,
% which follows 
% % directly
% % immediately %implies 
% by the argument for the Gersten conjecture 
% % like 
% in \cite{Voe-hty-inv}.
% An application of Part (1) of Proposition \ref{cor:contractinggenpointhomoveretaSmX} is the result,
% which follows directly by the argument for the Gersten conjecture in \cite{Voe-hty-inv}.
A direct application of Part (1) of Proposition \ref{cor:contractinggenpointhomoveretaSmX} is the result,
which follows by the argument for the Gersten conjecture in \cite{Voe-hty-inv}.
\begin{theorem*}[\protect{\Cref{th:CousinsubsetsLocEssSm}}]
Under the assumptions of
\Cref{cor:contractinggenpointhomoveretaSmX},
for any $\SH(B)$-representable cohomology theory $E^*$,
the Cousin complex
\cite{BlochOgus}
% , see \cite{BlochOgus}, 
\begin{equation}\label{eq:Cous(calV,E)} %intro:
0\to E^0(\calV)\to E^0(\calV^{(0)})\to %\bigoplus_{y\in \calV^{(1)}}E_y^{1}(\calV)\to 
\dots \to\bigoplus_{y\in \calV^{(r)}}E_y^{r}(\calV)\to\dots 
% \dots \bigoplus_{x\in \calV^{(d)}}E_x^{d}(\calV)\to 0
% ,    
\end{equation}
% see \cite{BlochOgus},
% where $\calV^{(c)}$ is the set of points of codimension $c$,
% and see \cite{BlochOgus} for the definition, 
is acyclic%,
% for the scheme $\calV=U\times_B \eta$ for an essentially smooth local scheme $U$ of dimension $d$, 
.
% see
% \cite{BlochOgus} for the notion of Cousin complexes. 
% in particular, the claim holds for $\SH^{S^1}(\mathrm{Corr}^\fr(B))$-representable or $\SH(B)$- cohomology theories.
\end{theorem*}
As shown in \cite[\S 1.2, 1.3.3]{ColumnsCousinlocesssmXoverB}
the above theorem has further applications: 
(1) 
the acyclicity of the Cousin complex
\begin{equation*}\label{eq:Cous(U,E)}
0\to E^0(U)\to E^0(U^{(0)})\to \bigoplus_{x\in U^{(1)}}E_x^{1}(U)\to 
\dots \to\bigoplus_{x\in U^{(c)}}E_x^{c}(U)\to\dots 
\end{equation*}
above the terms $E_x^{1}(U)$
\footnote{
Over a DVR,
an alternative proof
is provided by
\cite{Panin_MovLoemmasAdz} %by Panin
},
and 
(2) %Section 1.3.3
isomorphisms
\[\pi_{i+j,j}(\Sigma^\infty_{\PP^1} Y)(U)\cong 0, \quad Y\in\Sm_B,\, i<-1,\, j\in\mathbb Z,\]
which is 
a Zariski local improvement of 
the stable connectivity theorem \cite{ConnBase,ConnDodekindDomains,SS}
generalising 
% Morel's results \cite{Mor0,Mor1} with respect to 
results for $\SH(k)$ from \cite{Mor0,Mor1} by Morel with respect to 
\cite{Ayo06} by Ayoub.
Part (2) of Proposition \ref{cor:contractinggenpointhomoveretaSmX}
has no role for \Cref{th:CousinsubsetsLocEssSm}
and has a crucial role for \Cref{th:strhominv}.
Part (2) provides a compactification of the $\A^1$-homotopy from Part (1).
% Part (2) of
% Proposition \ref{cor:contractinggenpointhomoveretaSmX},
% providing a compacitifaction of the constructed $\A^1$-homotopy,
% has no role for \Cref{th:CousinsubsetsLocEssSm}
% and has a crucial role for \Cref{th:strhominv}. 
%
% Part (2) of
% Proposition \ref{cor:contractinggenpointhomoveretaSmX}
% has no role for %deduction of 
% \Cref{th:CousinsubsetsLocEssSm}
% and
% has a crucial role 
% % in our approach
% for \Cref{th:strhominv}
% providing 
% % a compacitifaction of 
% % being related
% % relating
% % to 
% a compacitifaction of 
% the constructed $\A^1$-homotopy. 
In detail,
% it 
the claim
means that 
the compactification of the homotopy
sends % takes
the closed point $(\infty,x) \in\PP^1\times\Xx$
% maps 
to %inside 
the complement of the closure of $Z$ 
in the compactification of $X_{\underline{\eta}}$
over $B$.
The key role in the construction
% in our method 
belongs to 
% We consider
certain %appropriate 
projective compactifications of 
$B$-schemes
over the one-dimensional local irreducible base scheme $B$  
giving rise to compacitfied homotopies of so called focused compactified 
framed correspondences,
% of dimensions ,
see \Cref{def:compactifiedfocusedFr}.
% Saying compacitfied homotopies 
% we mean 
% Saying compacitfied homotopy
% we mean 
% that the data defining an $\A^1$-homotopy of correspondences % compactified correspondences 
% naturally extends along the embedding of $\A^1\times U$ to $\PP^1\times U$, $U\in\EssSm_B$.
A compacitfication of an $\A^1$-homotopy 
above
means 
an extension of the data defining the homotopy %of correspondences % compactified correspondences 
along the embedding of $\A^1\times U$ to $\PP^1\times \overline{U}$
for a given
% closed immersio
compactification
$U\to\overline{U}$.
% In other words
% This is 
% this
We
construct it
as 
a focused compactified
one-dimensional
framed 
correspondence
equipped with a special morphism to 
the pair
$(\A^1\times U,\PP^1\times \overline{U})$.
% $(\A^1,\PP^1)$.
% $\PP^1$
% as the compactification.
% We
% consider 
% compactified
% $\A^1$-homotopies 
% as compactified
% one-dimensional
% framed correspondences
% equiped with a special morphism to $\PP^1$.
%
% , 
% where
% $\overline{U}$ is a compcitifcation of
% $U\in\EssSm_B$.
% an %appropriate
% extension of the data defining an $\A^1$-homotopy % compactified correspondences 
% along the embedding of $\A^1\times B$ to $\PP^1\times B$.
% Compatifications of correspondences 
A compatification of a $d$-dimensional framed correspondence 
% include two parts
includes two parts,
% according to our terminology 
% called 
% according to our terminology 
we call
``compactification'' and ``focusing''.
The role of ``compactification'' is similar to such instruments 
used over base fields
as
Quillen's trick from \cite{Quillen} %construction 
or standard triples in \cite{Voe-hty-inv}.
The ``focusing'' 
is a positive-dimensional base schemes specific part;
% setting; 
% concen%above 
% An example of ``focusing'' %above 
an example
could be given by %the spectrum of DVR 
% a one-dimensional local irreducible base scheme $B$ itself
the base scheme $B$ itself
considered as a compactification of its generic point $\eta$,
or
a smooth projective $B$-scheme %$X$ 
as a compactification of the generic fibre. % $X\times_B\eta$. 
\subsubsection{Directions of 
further applications.}
%$\square$-homotopy 
Being equipped by the construction
with compactifications in the %both of the 
above senses,
both of which are important in the $\square$-homotopy motivic studies,
% in distinct to the previously known ones
our homotopies and correspondences
form an instrument
to obtain 
certain analogs of
\Cref{th:strhominv,th:CousinsubsetsLocEssSm}
for $\square$-homotopy invariant theories. 
% Here 
% we
We 
mean a respective
generalisations of results 
in \cite{zbMATH07183731,zbMATH07654583}
on reciprocity sheaves,
which are 
supposed 
to be 
a natural
%is a 
further development 
of this work.

\subsection{Proof strategy and overview}\label{sect:ProofStrategy}% About the p
% \subsection{Plan of the proof}
%While it could be considered that 
%the direct way to cover the case of an arbitrary field in the motivic homotopy theory 
%is to consider more general definitions and constructions that would be applicable to spectra $\operatorname{Spec} K$ of finitely generated fields that are not essentially smooth, the proof is obtained by technique of compctifified correspondences and homotopies.
% In distinct 
% to the %often used 
% induction on the dimension of $X$ and 
% consequent reduction %of the case of henselian local scheme $X$ 
% to the case of $X = \operatorname{Spec} K$, %that was needed in the previous strategy,
% The %classical 
% arguments in %from 
The %classical 
proof strategy 
of graded isomorphism
% \begin{equation*}\label{eq:copy:strA1hominvNis}H^i_\nis(\A^1_k\times X,F_\nis)\cong H^i_\nis(X,F_\nis),\end{equation*}
% \begin{equation}\label{eq:copy:strA1hominvNislh}H^*_\nis(\A^1_k\times X^h_x,F_\nis)\cong F_\nis(X^h_x),\end{equation}
$H^*_\nis(\A^1_k\times X,F_\nis)\cong F_\nis(X)$
% $H^*_\nis(\A^1_k\times X^h_x,F_\nis)\cong F_\nis(X^h_x)$
in %from 
\cite{Voe-hty-inv,hty-inv}
% ,
includes 
% has two 
the following % elem
parts:
(1)
% an isomorphism $H^i_\nis(\Spec K\times\A^1,F)\simeq 0$ for any field extension $K/k$,
an isomorphism 
\begin{equation}\label{eq:HnisSpecKtAFcongzero}
% H^i_\nis(\Spec K\times\A^1,F)\simeq 0
H^*_\nis(\Spec K\times\A^1,F)\simeq F(\Spec K)
\end{equation}
for any field extension $K/k$,
% the vanishing of cohomologies $H^i(\Spec K\times\A^1,F)$, for a field extension $K/k$,
(2) 
% an induction on $\dim X$, %the dimension of 
% $X\in\Sm_k$. 
an induction on $\dim X$, %the dimension of 
via the codimension filtration on $X$, 
$X\in\Sm_k$.
% uses
% the induction on the dimension of $X\in\Sm_k$,
% and the vanishing of cohomologies $H^i(\Spec K\times\A^1,F)$, for a field extension $K/k$. 
% %and 
% % consequent reduction %of the case of henselian local scheme $X$ 
% % to the case of $X = \operatorname{Spec} K$, %that was needed in the previous strategy,
% The present proof does not use such induction, 
% and
% splits the claim differently on cohomologies 
%
% The present strategy does not use such induction, 
% and
% splits the proof of an isomorphism $H^*_\nis(X^h_x\times\A^1,F)\simeq 0$
% % , where $X^h_x$ is a local henselian scheme of $X\in \Sm_k$ at $x\in X$, to the 
% into the vanishing of
% the cohomologies 
% in certain senses
% % on $X^h_x\times\PP^1_k$
% ``out of infinity'' and ``at infinity''. 
%
%
% So $X\times\A^1$ is split 
% % by
% with respect to codimension filtration on $X$.
%
% into suspensions of $\Spec K\times\A^1$.
% So $X\times\A^1$ is split into suspensions of $\Spec K\times\A^1$.
The present strategy
splits the proof 
% of \eqref{eq:copy:strA1hominvNislh}
oppositely
% an isomorphism 
% $H^*_\nis(X^h_x\times\A^1,F)\simeq 0$
% $H^*_\nis(X^h_x\times\A^1,F)\simeq F(X^h_x)$
into 
the vanishing of the cohomologies
% ``out of infinity'' and ``at infinity'' 
on $X^h_x\times\A^1$
``out of infinity'' and ``at infinity''
% compactified by $X^h_x\times\PP^1$
in certain senses.
% , see \Cre{}
%
% The present
% proof does not use the mentioned induction, and
% splits the claim on cohomologies 
% $H^*(X^h_x\times\A^1,F)$, where $X^h_x$ is a local henselian scheme of $X\in \Sm_k$ at $x\in X$, to the 
% cases of cohomologies out of the infinity and at the infinity in appropriate senses.
%
% The case of cohomologies out of the infinity can be covered by the direct generalisation of the techniques developed earlier
% over base fields to a henselian local base scheme. 

% Cohomologies ``out of infinity'' %on $X^h_x\times\A^1$, 
% are studied 
Triviality of cohomologies ``out of infinity''
% , %on $X^h_x\times\A^1$, i.e. the isomorphism
means the isomorphism
%dispaly% \begin{equation*}\label{eq:copy:cohNisFinunNis}
$
    {H}^*_{\Fin_{X^h_x}}(X^h_x\times\A^1,F)\cong F(X^h_x), %\times\A^1
$
% \end{equation*}
see \Cref{def:FinS} and \Cref{th:FinSuppLocHensbaseafflinetriviality}, 
% and
% is 
proven 
by a 
% straightforward
% direct
generalisation of 
the techniques 
used earlier % from \cite{Voe-hty-inv,hty-inv}
over the base fields 
to the local henselian 
base schemes setting,
% Cohomologies ``out of infinity'' 
% on $X^h_x\times\A^1$
% are studies 
% % covered 
% by a direct
% generalisation of the techniques of the proof of \eqref{eq:HnisSpecKtAFcongzero}
% % $H^i_\nis(\Spec K\times\A^1,F)\simeq 0$ 
% from \cite{Voe-hty-inv,hty-inv}
% from the base fields 
% to the local henselian base schemes.
% N
namely,
% we do this
% straightforwardly
% for
% the proof of \eqref{eq:HnisSpecKtAFcongzero}
% from \cite{GWStrHomInv,ccorrs}.
% we use such
% a straightforward
% generalisation
% of 
% % the argument 
% the proof of \eqref{eq:HnisSpecKtAFcongzero}
% % in \cite{GWStrHomInv,ccorrs}.
% from \cite{GWStrHomInv,ccorrs}.
the proof of \eqref{eq:HnisSpecKtAFcongzero}
% in \cite{GWStrHomInv,ccorrs}.
from \cite{GWStrHomInv,ccorrs}
adimts such a straightforward generalisation.
% there is 
% such
% a straightforward
% generalisation
% % is do
% % of 
% % the argument 
% for
% the proof of \eqref{eq:HnisSpecKtAFcongzero}
% % in \cite{GWStrHomInv,ccorrs}.
% from \cite{GWStrHomInv,ccorrs}.
% % for the proof of 
% to prove \eqref{eq:HnisSpecKtAFcongzero}.
% replacing 
% $\Spec K$ to $X^h_x$.
% such %the 
% a
% lifting 
% the proof of \eqref{eq:HnisSpecKtAFcongzero}
% from \cite{GWStrHomInv,ccorrs}
% from $\Spec K$ to $X^h_x$.
%
% such %the 
% a
% lifting 
% the proof of \eqref{eq:HnisSpecKtAFcongzero}
% from \cite{GWStrHomInv,ccorrs}
% from $\Spec K$ to $X^h_x$.
% Cohomologies ``out of infinity'' 
% on $X^h_x\times\A^1$
% are studies 
% % covered 
% by a direct
% generalisation of the techniques of the proof of \eqref{eq:HnisSpecKtAFcongzero}
% % $H^i_\nis(\Spec K\times\A^1,F)\simeq 0$ 
% from \cite{Voe-hty-inv,hty-inv}
% from the base field $K$ 
% to the base scheme $X^h_x$ setting.
% from the base field to the base scheme setting.
%
% This implies 
% a strict homotopy invariance theorem
% for curtain subtopology of the Nisnevich topology.
% This implies 
% a strict homotopy invariance theorem
% for curtain subtopology of the Nisnevich topology.

The vanishing of the cohomologies ``out of infinity'' means an isomorphism 
\[H^*_\nis((X^h_x\times \PP^1)\times_{\PP^1}(\A^1),F)\simeq F((X^h_x\times \PP^1)\times_{\PP^1}(\A^1))\]
% \[H^*_\nis((X^h_x\times \PP^1)\times_{\PP^1}(\A^1),F)\simeq 0\]
The proof uses two ingredients: 
a moving homotopy %lemma 
provided by
\Cref{cor:contractinggenpointhomoveretaSmX} 
applied over a ``compactified base scheme'' $(\A^1_k,\PP^1_k)$,
% and the $\A^1$-invariance of cohomologies ``out of infinity'' mentioned above
% applied 
% via the criterion 
and 
the partial strict homotopy invariance
provided by the above result on cohomologies ``out of infinity''. 
% % $\A^1$-invariance result of cohomologies ``out of infinity'' mentioned 
% and 
% the above result on cohomologies ``out of infinity'' 
% % $\A^1$-invariance result of cohomologies ``out of infinity'' mentioned 
% applied 
% using \Cref{cor:preimFinFrCintr->Inj}.
% % \Cref{prop:FrA1HomFinCovPreimage,cor:preimFinFrCintr->Inj}.
\Cref{cor:contractinggenpointhomoveretaSmX} is proven by 
% a moving lemma based on 
the technique of compactifications discussed above. The construction of the homotopy requires that the base field contains three different elements used to control values of the compactified homotopy at $0$, $1$, and $\infty$. 
% This proves the claim for a field that has at least 3three elements, and the general claim follows by the reduction from finite base fields to infinite ones from \cite{DrKyllfinFrpi00,five-authors}. 
The claim for any base field follows by the argument
% reduction from finite base fields to infinite ones 
from \cite{DrKyllfinFrpi00,five-authors}.
The argument for strict homotopy invariance theorem %\Cref{th:strhominv} 
in \Cref{sect:StrHomInvTheorem} 
summarises 
% the triviality of cohomoligies ``out of infinity'' and ``at infinity'' 
% provided by 
results of
\Cref{sect:outofinffinitesupSHI,sect:InjThbyFinSupHom,sect:homotpiesatinfinity}.
\Cref{sect:CousinExactness} deduces the result on the Cousin complexes.
Recollection 
of framed correspondences and some used homotopies
is started in \Cref{sect:FrCorPresh} and extended in \Cref{sect-app:FrHomotopies}.
% Definitions of framed correspondences and presheaves
% and some used framed correspondences and homotopies
% are recollected in \Cref{sect:FrCorPresh} and \Cref{sect-app:FrHomotopies}.
% \Cref{sect:Fr} recollects the definitions of framed correspondences, and presheaves,
% and \Cref{sect-app:FrHomotopies} recollects some used framed correspondences and homotopies from other articles.
\Cref{sect-app:RelDim} 
formulates the used result 
on the relative dimension 
of irreducible schemes
over one-dimensional schemes
with non-empty generic fibre.
% contains the formulation of the result on the relative dimensions over one-dimensional schemes used in \Cref{sect:homotpiesatinfinity}.
\Cref{sect-app:SerreTh} collects results on line bundles.
\subsection{Conventions}
%, notation, and notions
% \eqref{sect:corrpreimZ}
% Given a correspondence that defines a span of schemes $X\xleftarrow{p} S\rightarrow Y$,
% for a closed subscheme $Z$ in $Y$,
% $c^{-1}(Z)=\overline{p(S\times_Y X)}$

% Let us list used notation of standard objects, and definitions of some extra or modified notions which we use in the text.
\subsubsection{Categories of schemes, relative dimension}\label{sect:not:Schemes}
% Depending on the context, we usually use notation $S$, $B$, or $V$ for base schemes in \Cref{sect:FrCorPresh}, \Cref{sect:outofinffinitesupSHI}, \Cref{sect:homotpiesatinfinity} %sect:CompactifiedCorrespondences
% respectively.
% We consider noetherian separated base schemes.
$\Sch_S$, $\Aff_S$, $\Sm_S$, $\SmAff_S$, $\Et_S$ are the categories of schemes, affine schemes, smooth schemes, smooth affine schemes, and \'etale schemes over a base scheme $S$.
We say a scheme $X$ has pure dimension $d$ and write $\dim X$, if each irreducible component of $X$ has dimension $d$.
We say that $X\in \Sch_S$ is equidimensional of pure relative dimension $d$ over $S$, 
and write 
% \begin{equation}\label{eq:dimSX=d}
$\dim_S X=d$
% \end{equation}
if for each $z\in S$, $\dim X\times_S z=d$.
 % $X\in \Sch_V$ is equidimensional of relative dimension $d$ over $V$, if for each $z\in V$, $\dim X\times_V z=d$
We say that a closed subscheme $X^\prime$ of $X$ has positive relative codimension, 
and write $\codim_{X/S}(X^\prime)>0$,
if $\codim (X^\prime\times_V z)/(X\times_S z)>0$ for each $z\in S$. 
$T_{X/S}$ denotes the tangent sheaf of a scheme $X$ over a scheme $S$.
Usually 
in \Cref{sect:FrCorPresh}, \Cref{sect:outofinffinitesupSHI}, \Cref{sect:homotpiesatinfinity},
% sect:CompactifiedCorrespondences
the base schemes 
are denoted $S$, $B$ and $V$ respectively.

% , and by $\Sm_S$ the category of smooth schemes, and by $\SmAff_S$ the category of smooth affine ones. %, $\EssSm_B$
% symbol $V$ is used as the base for the relative affine line $\A^1\times V$,
% and if we consider $\A^1\times V$ as a scheme over $\A^1$ or $\PP^1$, then we usually denote the base scheme by $B$.
%
%\subsubsection{Relative dimension}
%

% We denote the site given by a Grothendieck topology $\tau$ on a category $C$ by $C^\tau$.
% In particular, we write $\Et_S^\nis$ and $\Sm_S^\nis$ for the sites of the Nisnevich topology on the categories of \'etale and smooth schemes over a scheme $S$ respectively.

\subsubsection{Functions, sections, subschemes}\label{sect:not:FSSubschemes}

% For a presheaf %or a sheaf 
% $F$ on a site $\mathcal C$, we denote by 
% $F(X)$ or $\Gamma(X,F)$ the set of sections on $X\in \mathcal C$.
$F(X)=\Gamma(X,F)$ is the set of sections of 
a presheaf $F$ on $X$.
%disp
% In particular, for an invertible sheaf or a line bundle $\mathcal L$ on a scheme $X$, we write $\Gamma(U,\mathcal L)$ for the group of sections on an open subscheme $U$ in $X$.
We use notation 
% \[%todo%display.
$\mathcal O(X)=\mathcal O_X(X)=\Gamma(X,\mathcal O_X)$,
% \]
where $\mathcal O_X$ the sheaf of regular functions on the small Zariski site over $X$.
% $\overline{Z}=\mathrm{Cl}_{Z}$ 
% denotes the closure of $Z$ in $X$.
% The closure of $Z$ in $X$ is denoted
% \[\overline{Z}=\mathrm{Cl}_{X}(Z)\]

Given morphism $f\colon Z\to X$, 
\begin{equation}\label{eq:Cl=overline}
    \mathrm{Cl}_{X}(Z)=\mathrm{Cl}_{f,X}(Z)=\mathrm{Cl}_{X}(f(Z))=\overline{f(Z)}
    % X_f(W)=\overline{Z}=\mathrm{Cl}_{X}(Z)
\end{equation}
is the closure of the image;
we skip $f$ above, 
when it is defined by the context
or $Z$ is a constructible subscheme of $X$.
% is %injective 
% an immersion. 
% denotes the 
% \eqref{sect:corrpreimZ}

Given a scheme $X$, we write
% \[i_0,i_1\colon \calV\to \calV\times\A^1\] 
\[i_0\colon X\simeq X\times \{0\}\to X\times\A^1,\quad i_1\colon X\simeq X\times \{1\}\to X\times\A^1\] 
for the zero and unit sections.
% $i_0,i_1\colon \calV\to \calV\times\A^1$ are the zero and unit sections, 

Given a correspondence that defines a span of schemes $X\xleftarrow{p} S\rightarrow Y$
for a closed subscheme $Z$ in $Y$, 
$c^{-1}(Z)$ is the closure of the image of $S\times_Y X$ along $p$, i.e.
\begin{equation}\label{sect:corrpreimZ}
c^{-1}(Z)=\mathrm{Cl}_{X}({p(S\times_Y X)}).
\end{equation}

For a closed subscheme $Z\subset X$, we denote by $\mathcal I_{Z/X}=I_{X}(Z)$ 
the sheaf of ideals in $\mathcal O_X$ formed by functions that vanish on $Z$, and $I_{Z/X}=I_{X}(Z)\subset \mathcal O(X)$ denotes the global sections.
% For a closed subscheme $Z\subset X$, $\mathcal I_{Z/X}$ denotes the ideal sheaf in $\mathcal O_X$ given by functions that vanish on $Z$, and $I_{Z/X}\subset \mathcal O(X)$ denotes the global sections.
The vanishing locus of a sheaf of ideals $I\subset \calO(X)$ is denoted by $Z(I)$.
The vanishing locus of a section $s\in \Gamma(X,\mathcal V)$ of a vector bundle $\mathcal V$ is $Z(s)$.
% ; note that $Z(s)=Z(I)$, where $I$ is locally generated by $s$. % in view of the trivialisation of $V$. 
%
% The common vanishing locus of a set $s$ of sections $s_j$ of vector bundles $\mathcal V$, $j\in J$, is denoted by $Z(s)$ or $Z(s_j)_{J\in J}$.
%
% We write %\[s\big|_{Z}=r^*(s)\in r^*(F)\] 
% $s\big|_{Z}=r^*(s)\in r^*(F)$ %todo%display
% for the 
% inverse image of $s\in \Gamma(X,F)$ 
% for a presheaf $F$ on a small Zariski site of $X$
% along an injective morphism $r\colon Z\to X$.
We write %\[s\big|_{Z}=r^*(s)\in r^*(F)\] 
$s\big|_{Z}=r^*(s)$ %todo%display
for the inverse image of $s\in \Gamma(X,F)$ 
% for a presheaf $F$ 
% on a small Zariski site of $X$
along a given 
constructible injective morphism $r\colon Z\to X$.
% injective morphism $r\colon Z\to X$.

We denote 
by 
$X^{(c)}$
the set of points 
and 
by 
$X^{[c]}$ 
the set of closed subschemes
of codimension $c$ in a scheme $X$.
We use notation 
$\mathcal O_{X,x}$ 
for 
the local ring of $X$ at a point $x\in X$, 
$\mathcal O^h_{X,x}$ 
for
the henselisation,
and 
% We use notation 
\[\begin{array}{ll}
X_x = \Spec\mathcal O_{X,x},&
X^h_x = \Spec\mathcal O^h_{X,x},
\\
\can\colon X_x\to X,& 
\can^h\colon X^h_x\to X,    
\end{array}\]
for schemes and canonical morphisms.
% We denote 
% by 
% $\mathcal O_{X,x}$ the local ring of 
% % a scheme 
% $X$ at a point $x\in X$ 
% and 
% by 
% $\mathcal O^h_{X,x}$ the henselisation,
% and write 
% % We use notation 
% \[\begin{array}{ll}
% X_x = \Spec\mathcal O_{X,x},&
% X^h_x = \Spec\mathcal O^h_{X,x},
% \\
% \can\colon X_x\to X,& 
% \can^h\colon X^h_x\to X    
% \end{array}\]
% for schemes and canonical morphisms.
% for the corresponding affine schemes.
Similarly, %notation 
$\mathcal O_{X,Z}$ and $\mathcal O^h_{X,Z}$
% we use for 
are
the localisation and the henselisation of $\mathcal O_X$ at $\calI_{Z/X}$ for a closed subscheme $Z$ in $X$,
% We write
$X_Z = \Spec\mathcal O_{X,Z}$,
$X^h_Z = \Spec\mathcal O^h_{X,Z}$.

For $X\in\Sch_B$, $\nu\in B$,
we use notation
\[X_{\underline{\nu}}=X\times_B \nu.\]

\subsubsection{Cohomologies with respect to a topology.}\label{subsect:Htaunotations}

Given morphism $\widetilde X\to X$, there are the \v{C}ech complex $\check C_{\widetilde X}(X,F)$ and the extended \v{C}ech complex $\underline{\check C}_{\widetilde X}(X,F)$:
\[\begin{array}{llll}
\phantom{.}[\dots 0\to0&\to  &F(\widetilde X)\to& \dots\to F(\widetilde X^{\times (l-1)})\to F(\widetilde X^{\times l})\to \dots ],
\\ 
\phantom{.}[\dots 0\to F(X)&\to &F(\widetilde X)\to& \dots\to F(\widetilde X^{\times (l-1)})\to F(\widetilde X^{\times l})\to \dots ],
\end{array}\]
where $F(\widetilde X)$ is located in the zeroth deg.
For a Grothendieck topology $\tau$ on a category $\mathcal C$, and 
a presheaf $F$ on $\mathcal C$, we denote by $F_\tau$ the sheafification,
and use notation
% write 
% \begin{equation}
$H_\tau^*(X,F)=H_\tau^*(X,F_\tau)$.
% \end{equation} 
% % \begin{equation}H_\tau^*(X,F_\tau)=H_\tau^*(X,F)\end{equation} 
% for cohomologies of $F_\tau$ on $X\in \mathcal C$. 
% % sometimes without subscript $\tau$ for $F$ to short notation.
% % Denote by 
% % \begin{equation}H_\tau^*(X,F_\tau)=H_\tau^*(X,F)\end{equation} the cohomologies of $F_\tau$ on $X\in \mathcal C$;
% % we write sometimes $F$ without subscript above to short notation.
% %
% % We use also the notion $\unH^*_\tau(X,F)$ 
% % for what we call reduced cohomoligies defined as follows:
\begin{definition}
We call by reduced cohomologies the following groups
\begin{equation}\label{eq:redcohunHtau}
\begin{cases}
\unH^i_\tau(X, F)=0, i<-1,\\
\unH^{-1}_\tau(X, F)=\ker(F(X)\to F_\tau(X)),\\
\unH^0_\tau(X, F)=\coker(F(X)\to F_\tau(X)),\\
\unH^i_\tau(X, F)=H^i_\tau(X,F), i>0.
\end{cases}
\end{equation}
\end{definition}
\begin{remark}
Note that the following conditions are equivalent
% The following conditions are equivalent
% \[\unH^*_\tau(X,F)=0\Longleftrightarrow\begin{cases}H^0_\tau(X, F_\tau)=F(X), \\H^i_\tau(X, F_\tau)=0, i>0.\end{cases}\]
\[\unH^*_\tau(X,F)=0\Longleftrightarrow H^*_\tau(X,F)\cong F(X)\Longleftrightarrow\begin{cases}H^0_\tau(X, F_\tau)\cong F(X), \\H^i_\tau(X, F_\tau)=0, &i>0.\end{cases}\]
% \[\left[\unH^*_\tau(X,F)=0\right]\Longleftrightarrow\left[\begin{cases}H^0_\tau(X, F_\tau)=F(X), \\H^i_\tau(X, F_\tau)=0, i>0.\end{cases}\right]\]
\end{remark}

% \[\begin{array}{lllll}
% \check C_{\widetilde X}(X,F)&=&[\dots 0\to0\to  &F(\widetilde X)\to& \dots\to F(\widetilde X^{\times (l-1)})\to F(\widetilde X^{\times l})\to \dots ],\\ \underline{\check C}_{\widetilde X}(X,F)&=&
% [\dots 0\to F(X)\to &F(\widetilde X)\to& \dots\to F(\widetilde X^{\times (l-1)})\to F(\widetilde X^{\times l})\to \dots ]
% \end{array}\]
% where $F(\widetilde X)$ is located in the zeroth deg.

\subsubsection{Sites, and sheaves.}\label{sect:not:PreShaeves}
Given a Grothendieck topology $\tau$ on a category $C$, 
we denote the site 
by $C^\tau$.
% In particular, we write 
% $\Et_S^\nis$, $\Sm_S^\nis$, $\Sch^\nis_S$ for the 
% %big and small Nisnevich sites over $S$ respectively.
% Nisnevich sites with underlying  categories $\Et_S$, $\Sm_S$, $\Sch_S$ respectively.
$\Pre(S)$, $\Sh_\tau(S)$ and $\mathbb Z\Pre(S)$ and $\mathbb Z\Sh_\tau(S)$ denote 
the categories of presheaves and $\tau$-sheaves 
and additive presheaves and $\tau$-sheaves 
% of abelian groups 
of abelian groups 
on $\Sch_S$.
%
% $\Pre(S)$ and $\mathbb Z\Pre(S)$ denote the categories of presheaves and additive presheaves of abelian groups on $\Sch_S$.
% $\Sh_\tau(S)$ and $\mathbb Z\Sh_\tau(S)$ denote the categories of %Nisnevich 
% $\tau$-sheaves 
% and 
% %Nisnevich 
% $\tau$-sheaves of abelian groups on $\Sch_S$.
Given $F\in\Pre(S)$, and an $X\in\Sch_S$,
denote by $F_X\in\Pre(X)$ 
the inverse image of $F$ along the structure morphism $X\to S$.
% the restriction of $F$ along the functor $\Sch_X\to \Sch_S$.
% the base change of $F$ along the structure morphism $X\to S$.
We let us write $H_{\tau}^*(X,F)$ for $H_{\tau}^*(X,F_X)$ for a topology $\tau$ on $\Sch_X$.

\subsubsection{Notation for the Nisnevich topology}

% We write $\Et_S^\nis$, $\Sm_S^\nis$, $\Sch_S$ for the 
% %big and small Nisnevich sites over $S$ respectively.
% big and small Nisnevich sites 
% on the categories $\Et_S$, $\Sm_S$, $\Sch_S$ respectively.

Also we use notation $\Nis^X$ for the Nisnevich topology on $\Sch_X$ over a scheme $X$. 
%We usually use superscript or subscript '$\nis$'.
We let us write $F_\nis$ instead of $F_{\Nis^X}$ for a presheaf $F$ on $\Sch_X$.
We use the notation $\Nis^X_Z$ for the subtopology of $\Nis^X$ that is trivial over the open subscheme $X-Z$, for a closed subscheme $Z$, 
and also use the notation $\underline\Nis^X$ and $\underline\Nis^X_Z$ for certain subtopologies of $\Nis^X$ and $\Nis^X_Z$,
see \Cref{def:underNisV}.
% defined in 
%sect:(un)NisXZ

\subsubsection{Category of open pairs and linear additivisation}\label{sect:Pair}

% For any category $\mathcal C$ denote by $\mathcal C^\mathrm{arrow}$ that category
% % with objects being morphisms $\alpha\colon c_0\to c_1$ in $\matcal C$
% % and morphisms $\alpha\to \alpha^\prime$ given by pairs of morphisms $(f_0,f_1)$, $f_0$ .
% with objects being morphisms of $\mathcal C$ and morphisms given by commutative squares.
% $\alpha_1$ of arrows in $\mathcal C$, i.e.
% the category of functors $\mathrm{Func}([0,1],\mathcal C)$ from the category defined by the two element ordered set $[0,1]$.  
\begin{definition}\label{def:Cpair}
For any category $\mathcal C$ denote by $\mathcal C^\mathrm{arrow}$ that category
with objects being morphisms of $\mathcal C$ and morphisms given by commutative squares.
Given a %n essentially surjective 
functor of categories 
$\Sm_S\to \mathcal C$,
consider the subcategory of $\mathcal C^\mathrm{arrow}$
spanned by open immersians
and denote by $\mathcal C^\mathrm{pair}$ 
the quotient category 
with morphisms 
\[\mathcal C^\pair( (X_1,U_1), (X_2,U_2) )  = \Cofib ( \mathcal C( X_1 , U_2 ) 
\xrightarrow{c}
% \xrightarrow{\alpha\mapsto (j_2\circ \alpha, \alpha\circ j_1)} 
\mathcal C^\mathrm{arrow}( (X_1,U_1) , (X_2,U_2) ) ),\]
where
$c\colon{\alpha\mapsto (j_2\circ \alpha, \alpha\circ j_1)}$, 
$j_1\colon U_1\to X_1$, $j_2\colon U_2\to X_2$%
.
\end{definition}
\begin{definition}\label{def:ZFstar}
Given a category $\mathcal C$ enriched over pointed sets, 
the \emph{linear additivisation} 
% of $\mathcal C$ 
is 
the universal colimit functor
$\mathcal C\to\mathbb Z\mathcal C^{\oplus}$ to an additive category $\mathbb Z\mathcal C^{\oplus}$ enriched over abelian groups. %$\Ab$.
% quotient category of $\mathbb Z\mathcal C$ such that $\mathcal C(X_1\amalg X_2,Y) = \mathcal C(X_1,Y)\oplus \mathcal C(X_2,Y)$.
% Denote by $\ZF_*(S)$ and $\ZF_*^\pair(S)$ the additivisations of the categories $\Fr_+(S)$ and $\Fr_+^\pair(S)$.
% Let $\mathcal C$ be a category enriched over pointed sets. The \emph{linear additivisation} of $\mathcal C$ is the universal quotient category of $\mathbb Z\mathcal C$ such that $\mathcal C(X_1\amalg X_2,Y) = \mathcal C(X_1,Y)\oplus \mathcal C(X_2,Y)$.
% % Denote by $\ZF_*(S)$ and $\ZF_*^\pair(S)$ the additivisations of the categories $\Fr_+(S)$ and $\Fr_+^\pair(S)$.
\end{definition}

\section{Framed correspondences and presheaves}\label{sect:FrCorPresh}

% \Cref{sect:Fr}
The section
recollects some %of 
definitions %of framed correspondences and framed presheaves introduced in 
from \cite{Voe-notes,Framed,hty-inv}
and % \Cref{sect:Fr1th} 
recalls the notion of first-order framed correspondences from \cite{SmModelSpectrumTP}, 
which is between framed correspondences and normally framed correspondences introduced independently in \cite{GNThomSpectra,five-authors}.
% % \Cref{sect:FrPairs,sect:ZFr}
% Other subsections
% recollect notions of correspondences of pairs, framed presheaves, and framed linear presheaves.
\begin{remark}
    First-order framed correspondences allow us to make to some technical simplification, see %in 
    \Cref{prop:NisZFr}, %\Cref{sect:InjThbyFinSupHom},  
    %of arguments in next sections
    though 
    the original framed correspondences from \cite{Framed}
    % they is not necessary 
    would be  
    % are
    enough
    for our strategy.
    % of main results.
    % The use of first-order framed correspondences leads to some technical simplification of arguments in next sections, though it is not necessary.
\end{remark}
\begin{remark}
    % Note that t
    The definition of correspondences of pairs 
    %introduced 
    here %provided by \Cref{sect:Pair} 
    is slightly different from the one in \cite{hty-inv}.
\end{remark}
% \cite{GNThomSpectra}
% % by Neshitov, see \cite{GNThomSpectra} %Nesh:nfr_GN%todo
% % , 
% and %by Hoyois, Elmanto, Khan, Sosnilo, Yakerson in 
% \cite{five-authors}.

Throughout the section $S$ denotes a base scheme.

% \subsection{
% Framed correspondences 
% and first-order framed correspondences
% of schemes}
% \label{sect:Fr}
\subsection{
Correspondences categories and presheaves with transfers}
\label{sect:Fr}

% We recall the following notions from \cite{Voe-notes,Framed}.
% We recall the following notion from \cite{SmModelSpectrumTP}.
We recall the following notions from \cite{Voe-notes,Framed}
and \cite{SmModelSpectrumTP}.

\begin{definition}\label{def:explFr}

(1)
An \emph{explicit framed correspondence of level $n$} from $X$ to $Y$ in $\Sch_S$ 
is 
% a 
data $(Z,V,\varphi,g)$ %set of 
given by 
an \'etale neighbourhood $V$ of a closed subscheme $Z$ in $\A^n_X$ finite over $X$,
$\varphi=(\varphi_1,\dots,\varphi_n)\in \mathcal O_V(V)^{\oplus n}$, 
such that %the common vanishing locus 
$Z(\varphi)$ in $X$ equals $Z$,
and a regular map $g\colon V\to Y$.

(2)
\emph{A first-order framed correspondence} of level $n$
from $X$ to $Y$ is data
$(Z,\tau,g)$ given by a closed subscheme $Z$ in $\A^n_X$,
a trivialisation of the conormal sheaf of $Z$
$
\tau\colon \mathcal I_{Z/\A^n_X}/\mathcal I^2_{Z/\A^n_X}\simeq \mathcal O_Z^{\times n}
$,
and a regular map $g\colon \mathcal Z(I^2_{Z/\A^n_X})\to Y$.
Note that, the pair $(\tau,g)$ defines the morphism of schemes
\begin{equation}\label{eq:taugTtwoAY}
    (\tau,g)\colon Z(I^2_{Z/\A^n_X})\to \A^n\times Y
\end{equation}
such that 
$Z(I^2_{Z/\A^n_X})\times_{\A^n_X}(0\times Y)\cong Z$.
Denote by $\nrFr_n(X,Y)$ the set of such $(Z,\tau,g)$ painted at the element with $Z=\emptyset$.

\end{definition}

\begin{definition}

(1)
Denote by $\Fr_n(X,Y)$ the set of equivalence classes 
with respect to the equivalence relation
$(Z,V_1,\varphi_1,g_1)\sim(Z,V_2,\varphi_2,g_2)$ 
% from $X$ to $Y$ of level $n$ 
% are equivalent, 
% if there is 
given
by
common shrinkings $V$ of $V_1$ and $V_2$ %of $Z$ in $\A^n_X$ 
such that the inverse images of $(\varphi_1,g_1)$ and $(\varphi_2,g_2)$ are equal.
Define the category $\Fr_+(X, Y )$ with objects being the same as in $\Sm_S$, 
morphisms given by pointed sets $\Fr_+(X,Y) = \bigvee_n\Fr_n(X, Y )$ %\limits
pointed at the elements with $Z=\emptyset$,
see \cite[Definitoin 2.3]{Framed}.

(2)
Define the \emph{pointed category} $\nrFr_+(S)$ with morphisms $\bigvee_n\nrFr_n(X,Y)$,
and composite morphism 
for 
$(Z_1,\tau_1,g_1)\in \nrFr_{n_1}(X_0,X_1)$ and 
$(Z_2,\tau_2,g_2)\in \nrFr_{n_2}(X_1,X_2)$ given by
$(Z_{12}, \tau_{12},g_{12} )\in \nrFr_{n_1+n_2}(X_0,X_2)$, where
% \[
$Z_{12}= Z_1\times_{X_1}Z_2$, 
$(\tau_{12},g_{12})\colon Z(\mathcal I^2_{Z_{12}}) \to Z(\mathcal I^2_{Z_{2}}) \to X_2$, 
% ,\]
the first morphism is induced by $(\tau_2,g_2)$, and the second one is induced by $(\tau_1,g_1)$.

(3)
% \begin{definition}
    Using \Cref{def:Cpair}, %In particular, 
we define categories $\Sm^\pair_S$ and $\Fr^\pair_+(S)$.
Denote $\nrFr_+(\Sm^\pair_S)=(\nrFr_+(S))^\pair$.
\end{definition}
% Then 
% Then according to the above definition 
\begin{remark}
$\nrFr_n((X,U),(Y,V))$ is the pointed set that is the quotient 
of the set of the elements $(Z,\tau,g)\in \nrFr_n(X,Y)$ such that
\[
((X\setminus U)\times_X Z)_\mathrm{red}\supset (Z\times_Y (Y\setminus V))_\mathrm{red},
\]
with respect to the subset of elements $(Z,\tau,g)\in \nrFr_n(X,Y)$ such that
$\emptyset = Z\times_Y (Y\setminus V)$.    
\end{remark}
% \end{definition}

% \subsection{First-order framed correspondences}\label{sect:Fr1th}

% We recall the following notion from \cite{SmModelSpectrumTP}.

% \begin{definition}
% % Denote the pointed set of such $(Z,\tau,g)$ by $\nrFr_n(X,Y)$. %painted at the element with $Z=\emptyset$.
% \end{definition}

\begin{example}
    $\sigma_X=(0\times X,t,g)\in\Fr_1(X,X)$, where 
    % $(t,g)\colon\A^1_X\to \A^1_S\times_S X$ is the canonical isomorphism.
    $t$ is the coordinate function on $\A^1_X$, and $g$ is the projection to $X$.
    For a matrix $G\in \GL_n(X)$, we define \[\sigma^G_{X}=(0\times V, G(t_1,\dots,t_n), \mathrm{pr}),\] 
    where $\mathrm{pr}\colon \A^n_V\to V$ is the canonical projection.
    Denote $\sigma^Df=f\circ \sigma^D\in\Fr_n(X,Y)$ for a given $S$-schemes $X$ and $Y$.
\end{example}

% \begin{definition}
A \emph{framed presheaf} 
% or
% a \emph{first-order framed presheaf}
is
a functor $\Fr_+(S)^{\mathrm{op}}\to \Set$.
A framed presheaf is \emph{quasi-stable}, if the morphisms $\sigma_X^*$ are auto-morphisms of $F(X)$ for all $X$. 
% 
% Similarly for 
% \emph{first-order framed presheaves}
% and $\nrFr_+(S)$.
% A first-order framed presheaf is a presheaf on $\nrFr_+(S)$ or $\nrFr_+(\Sm_S^\pair)$.
% \end{definition}
There is the functor 
$\Fr_+(S)\to \nrFr_+(S)$
% which preserves morphisms of the form $\sigma_X$, $X\in \Sm_S$,
defined by the natural map
\begin{equation}\label{eq:Fr->nrFr}
\Fr_n(X,Y)\to \nrFr_n(X,Y)\colon
(Z,V,\varphi,g)\mapsto (Z, \tau,g\big|_Z),
\end{equation}
where $\tau$ is induced by $\varphi\big|_{Z(\mathcal I^2_{\A^n_X}(Z))}$,
and similar definitions apply to
\emph{first-order framed presheaves}.
% and $\nrFr_+(S)$.

\begin{lemma}\label{lm:Fr<--nrFr}
% Suppose the base scheme $S$ is affine. 
Any $\A^1$-invariant quasi-stable framed presheaf $F$ on 
$\Aff_S$
% the category of affine schemes over 
for
an affine scheme $S$
is canonically equipped with the structure of a quasi-stable first-order framed presheaf.
\end{lemma}
\begin{proof}
Let $c,c^\prime\in \Fr_n(X,Y)$ be elements that images in $\nrFr_n(X,Y)$ are equal.
Denote $Z_{[2]}=Z(\mathcal I^2(Z))$.
Denote by $g_{Z_{[2]}}$ the morphism $g\big|_{Z_{[2]}}=g^\prime\big|_{Z_{[2]}}\colon Z_{[2]}\to Y$. 
% and by $g_{-\times Z_2}$ the composite with the projection $-\times Z_2\to Z_2$.
% and by $g_{U\times Z_2}$ the composite with the projection $U\times Z_2\to Z_2$.
% for each $U$.
% $U\times Z_2\to Z_2\xrightarrow{g_{Z_2}}Y$.
% define similarly $g_{\A^1\times Z_2}$, $g_{0\times Z_2}$, $g_{1\times Z_2}$.
% Note that $g\big|_{Z_2}=g^\prime\big|_{Z_2}\colon Z_2\to Y$.
% Denote this function by $g_{Z_2}$, and define similarly $g_{\A^1\times Z_2}$, $g_{0\times Z_{[2]}}$, $g_{1\times Z_{[2]}}$.
Consider the map 
\begin{gather}
\check{g}\colon 
E\to Y,\label{eq:amalggGlassY}\\
\begin{array}{llccccc}
E&=&(0\times (\A^n_X)^h_Z)&\amalg_{0\times Z_{[2]}}&(\A^1\times Z_{[2]})&\amalg_{1\times Z_{[2]}}&(1\times (\A^n_X)^h_Z)\\
\check{g}&=&g&\amalg_{g_{0\times Z_{[2]}}}& g_{\A^1\times Z_{[2]}}&\amalg_{g_{1\times Z_{[2]}}}& g^\prime
,\end{array}\nonumber
\end{gather}
where $g_{-\times Z_{[2]}}$ is the composite $(-\times Z_{[2]})\to Z_{[2]}\xrightarrow{g_2} Y$.
Choose a lifting of \eqref{eq:amalggGlassY}
\[\tilde g\colon (\A^1\times\A^n_X)^h_{\A^1\times Z}\to Y,\]
that exists by \cite[Theorem I.8]{Gru} see also \cite{Elkiksoleqhens}, \cite[Lemma 3.11]{FrRigidSmAffpairs},
% and \cite{Elkiksoleqhens}, %\cite[Lemma 3.6]{FrRigidSmAffpairs}todo add
since the canonical morphism $E\to (\A^1\times\A^n_X)^h_{\A^1\times Z}$ is an affine henseleian pair.
Consider 
% the framed homotopy
\[h = (\A^1\times Z,\lambda \varphi+(1-\lambda) \varphi^\prime, \tilde g)\in \Fr_n(\A^1\times X, Y).\]
Then
$h\circ i_0=c$, $h\circ i_1=c^\prime$, and since $F$ is $\A^1$-invariant,
it follows that $c^*={c^\prime}^*$.
Thus 
% the claim is proven.
the induced morphisms $c^*,{c^\prime}^*\colon F(Y)\to F(X)$ are equal.
Hence 
the %first 
claim of the lemma follows,
because the map \eqref{eq:Fr->nrFr} is surjective
and
% The second claim follows, because
% \eqref{eq:Fr->nrFr} 
preserves correspondences of the form $\sigma_X$, $X\in\Sm_S$.
% $\sigma_X\in \Fr_1(X,X)$ goes to $\sigma_X\in \nrFr_n(X,X)$ along morphism \eqref{eq:Fr->nrFr}.
\end{proof}

\subsection{Linearised correspondences, %framed presheaves %, 
and Nisnevich sheafification}\label{sect:ZFr}
% The canonical functor $\mathrm{Set}_\bullet\to \Ab$ form the category of pointed sets to the category of abelian groups induces  
% For a category $\mathcal C$ enriched over pointed sets, 
% denote by $\mathbb Z\mathcal C$ the category enriched over abelian groups that 
% with the same objects and morphisms given by the abelian groups associated with the pointed sets.
% have the same objects and hom-groups being the quotient of the free abelian groups associated with the hom-sets in $\mathcal C$ such that the pointed element equals zero.

% \begin{definition}\label{def:ZFstar}
% Let $\mathcal C$ be a category enriched over pointed sets. The \emph{additivisation} of $\mathcal C$ is the universal quotient category of $\mathbb Z\mathcal C$ such that $\mathcal C(X_1\amalg X_2,Y) = \mathcal C(X_1,Y)\oplus \mathcal C(X_2,Y)$.
% \begin{definition}
% Let
% $\ZF_*(S)$ and $\ZF_*^\pair(S)$ denote 
% the linear additivisations 
% of the categories $\Fr_+(S)$ and $\Fr_+^\pair(S)$, see \Cref{def:ZFstar}.
% \end{definition}

% \begin{definition}
% A \emph{framed presheaf} is
% a functor $\Fr_+(S)\to \Set$.
% A framed presheaf is \emph{quasi-stable}, if the morphisms $\sigma_X^*$ are auto-morphisms of $F(X)$ for all $X$.
% A first-order framed presheaf is a presheaf on $\nrFr_+(S)$ or $\nrFr_+(\Sm_S^\pair)$.
% % The same definitions apply to the case of the category of open pairs $\Fr_+^\pair(S)$, $\nrFr_+(\Sm^\pair_S)$.
% \end{definition}
\begin{definition}\label{def:ZFovZF}
(1) Let
$\ZF_*(S)$ denote 
the linear additivisations 
of the category $\Fr_+(S)$, see \Cref{def:ZFstar}.
%
% $\ZPre(\ZF_*(S))$ is the category of \emph{additive framed presheaves of abelian groups} 
% or a \emph{framed linear presheaves}
% over $S$.
$F\in \ZPre(\ZF_*(S))$ is called 
% an \emph{additive framed presheaf of abelian groups} 
% or 
a \emph{framed linear presheaf}
over $S$.
% is a functor $F\colon \Fr_+(S)\to \Ab$ such that $F(X_1\amalg X_2)\simeq F(X_1)\oplus F(X_2)$.
% Equivalently we can say that a linear framed presheaf defines a functor of categories enriched over abelian groups
% \[F\colon \ZF_*(S)\to \Ab.\]
(2)
Define
$\ovZF_n(X,Y)=\coker(\ZF_n(X\times\A^1,Y)\xrightarrow{i^*_0-i_1^*}\ZF_n(X,Y))$.
Similarly for $\Fr_+^\pair(S)$.
% The same definitions apply to %the case of 
% $\Fr_+^\pair(S)$.
\end{definition}

% \begin{definition}
% An \emph{additive framed presheaf of abelian groups} 
% or a \emph{linear framed presheaf}
% over $S$
% % over a base scheme $S$
% is a functor $F\colon \Fr_+(S)\to \Ab$ such that $F(X_1\amalg X_2)\simeq F(X_1)\oplus F(X_2)$.
% % Equivalently we can say that a linear framed presheaf defines a functor of categories enriched over abelian groups
% % \[F\colon \ZF_*(S)\to \Ab.\]
% The same definition applies to %the case of 
% $\Fr_+^\pair(S)$.
% \end{definition}

\begin{proposition}
% Consider the inverse and direct image functors %left Kan extension and the forgetful functors 
% $\gamma^*\colon \ZShS\to \ZShfrS$, $\gamma_*\colon \ZShfrS\to \ZShS$,
% see \Cref{sect:not:PreShaeves}.
% Then for any $E\in \ZShS$, and $F\in \ZShfrS$,
% \[\Ext^l_{\ZShS}(E,\gamma_* F)\simeq
% \Ext^l_{\ZShfrS}(\gamma^* E,F).\]
Consider the inverse and direct image functors %left Kan extension and the forgetful functors 
$\gamma^*\colon \ZShS\to \ZShfrS$, $\gamma_*\colon \ZShfrS\to \ZShS$,
see \Cref{sect:not:PreShaeves}.
Then for any $E\in \ZShS$ and $F\in \ZShfrS$,
\[\Ext^l_{\ZShS}(E,\gamma_* F)\simeq
\Ext^l_{\ZShfrS}(\gamma^* E,F).\]
\end{proposition}
\begin{proof}
Since the category $\ZShS$ is generated via colimits by the representable sheaves $\mathbb Z(U)$, $U\in \Sm_S$, it is enough to prove the claim for $E= \mathbb Z(U)$. 
Since
$\gamma^* E = \ZF_*(U)_\nis$,
the claim follows from \cite[Proposition 16.2]{hty-inv}.
\end{proof}

\begin{definition}\label{def:CohSupTransfer}
% For $S$-schemes $U$ and $X$, and a framed correspondence $c\in \Fr_n(U,X)$,
Let $U,X\in\Sch_S$, and $c\in \Fr_n(U,X)$.
% For any closed subscheme $Z$ in $X$, and the closed subscheme $Y=c^{-1}(Z)$ in $U$, 
For any closed subscheme $Z$ in $X$, and $Y=c^{-1}(Z)$, see \eqref{sect:corrpreimZ}, 
% \ZSh(\ZF_*(S))
for any $F\in\ZShfrS$,
% linear presheaf with framed transfers $F$ over $S$,
define the canonical morphism
on reduced cohomologies with support
% \[
$
c_*\colon\underline{H}^*_Z(X,F)\to \underline{H}^*_Y(U,F)
$
% \]
as the composite morphism
\begin{multline*}
H^l_{Z}(X,F)=
\Ext^l_{\ZShS}(\mathbb Z(X)/\mathbb Z(X-Z),F)\simeq\\
\Ext^l_{\ZShfrS}(\ZF(X)/\ZF(X-Z),F)\to\\
\Ext^l_{\ZShfrS}(\ZF(U)/\ZF(U-Y),F)\simeq\\
\Ext^l_{\ZShS}(\mathbb Z(U)/\mathbb Z(U-Y),F)=
H^l_{Y}(U,F).
\end{multline*}
\end{definition}

\section{Topology over a quotient-space $X/(X-Z)$} \label{sect:(un)NisXZ}sect:FrCorPresh
 % and cohomologies
\subsection{Topology $\Nis^X_Z$} %{Topologies.} 

Recall that $\Nis^X$ denotes the Nisnevich topology on $\Sch_X$ over a scheme $X$.
\begin{definition}\label{def:NisXZ}
Let $Z$ be a closed subscheme of $X$.
Define the topology $\Nis^X_Z$ on $\Sch_{X}$ as 
the strongest subtopology of $\Nis_X$ that restriction to %on 
$\Sch_{X-Z}$ is trivial.
\end{definition}
\begin{lemma}\label{lm:NisZXNisZpoints}
(1) An \'etale morphism $\widetilde U\to U$ is a $\Nis^X_Z$-covering %or $\underline{\Nis}^X_Z$-covering 
if 
% for any scheme $V$ of the forms \begin{equation}\label{eq:NisZpoints}U^h_v, \; v\in U,\quad U\times_X (X-Z),\end{equation}
% or \eqref{eq:XNisZpoints}, 
there is a dashed arrow in the triangle
\[\xymatrix{
&\widetilde U\ar[d]\\
V\ar@{-->}[ru]\ar[r]& U
}\]
for any scheme $V$ of the forms 
\begin{equation}\label{eq:NisZpoints}V=U^h_v, \; v\in U,\quad V=U\times_X (X-Z).\end{equation}
(2) The topology $\Nis^X_Z$ %and $\underline{\Nis}^X_Z$ 
has enough set of points given by the schemes \eqref{eq:NisZpoints}
% and
% \begin{equation}\label{eq:XNisZpoints}U\times_X X^h_{p(v)}, \quad U\times_X (X-Z),\end{equation}
for all $U\in\Sch_X$.
% and $v\in U$.
% where $p\colon U\to X$ is the canonical morphism.
% for all $U\in \Sm_X$.
\end{lemma}
\begin{proof}
% By \Cref{def:NisXZ}, %it follows that 
% $\Nis^X_Z$ is the 
% strongest topology on $\Sch_{X}$ such that 
% any $\Nis^X_Z$-covering $\widetilde U\to U$
% is 
% a Nisnevich coverings 
% and there is a dashed arrow in the diagram
% \begin{equation}\label{eq:liftNisXZ}
% \xymatrix{
% &\widetilde U\ar[d] \\
% U\times_X(X-Z)\ar@{-->}[ru]\ar[r]&U
% .}\end{equation}
% Then $\Nis^X_Z$ is the weakest topology on $\Sch_{X}$ 
% such that all
% Nisnevich coverings 
% $\widetilde U\to U$
% admitting lifting \eqref{eq:liftNisXZ}
% are $\Nis^X_Z$-coverings.
It follows from \Cref{def:NisXZ}, %it follows that 
that 
% $\Nis^X_Z$ is 
% the strongest %topology 
% and 
% the weakest topology 
% on $\Sch_{X}$ 
% such that 
a $\Nis^X_Z$-covering $\widetilde U\to U$
is 
a Nisnevich coverings 
such that there is a dashed arrow in the diagram
\begin{equation}\label{eq:liftNisXZ}
\xymatrix{
&\widetilde U\ar[d] \\
U\times_X(X-Z)\ar@{-->}[ru]\ar[r]&U
.}\end{equation}
% and all such morphisms are $\Nis^X_Z$-coverings.
% such that 
% any $\Nis^X_Z$-covering $\widetilde U\to U$
% is 
% a Nisnevich coverings 
% and there is a dashed arrow in the diagram
% \begin{equation}\label{eq:liftNisXZ}
% \xymatrix{
% &\widetilde U\ar[d] \\
% U\times_X(X-Z)\ar@{-->}[ru]\ar[r]&U
% ,}\end{equation}
% and all such morphisms are $\Nis^X_Z$-coverings.
% respectively.
% Then $\Nis^X_Z$ is the weakest topology on $\Sch_{X}$ 
% such that all
% Nisnevich coverings 
% $\widetilde U\to U$
% admitting lifting \eqref{eq:liftNisXZ}
% are $\Nis^X_Z$-coverings.
Then both points (1) and (2) follow.
% To prove the case of $\underline{\Nis}^X_Z$ one we note firstly that by the definition all the schemes \eqref{eq:XNisZpoints} are $\underline{\Nis}^X_Z$-paints.
% Then for any $\underline{\Nis}^X_Z$-local isomorphism $F\to G$ in $\Pre(X)$, see \Cref{sect:not:PreShaeves},  the morphisms of stalks $F(V)\to G(V)$ are isomorphisms for $V$ of the form \eqref{eq:XNisZpoints}. 
% Conversely,
% let  $F\to G$ be a morphism
% such that $F(U\times_X X^h_{p(v)})\simeq G(U\times_X X^h_{p(v)})$ for each $v\in U$.
% Then for each $v\in U$, there is an \'etale neighbourhood $\tilde X_v$ of $p(v)$ in $X$ such that $F(\tilde X_v)\simeq G(\tilde X_v)$.
% Define the $\underline{\Nis}^X_Z$-covering of $U$ $\widetilde X= \coprod_{v\in U}\tilde X_v\amalg (X-Z)\to U$.
% Then $F(\widetilde X)\simeq G(\widetilde X)$.
% So if both $F$ and $G$ are $\underline{\Nis}^X_Z$-sheaves then $F\cong G$.
%(2) follows from Point (1).
\end{proof}

\subsection{Framed transfers for cohomologies}\label{sect:fFrTransfersCohoutofinfinity}

% Consider the subcategories of sheaves in $\mathbb Z\Pre(X)$ with respect to topologies $\Nis^X_Z$ and $\underline{\Nis}^X_Z$. 
% Denote the sheafifications of $F\in\mathbb Z\Pre(X)$ by $F_{\Nis^X_Z}$, $F_{\underline{\Nis}^X_Z}$.
% Denote by $F_{\Nis^X_Z}$, $F_{\underline{\Nis}^X_Z}$
% the sheafifications of $F\in\mathbb Z\Pre(X)$ with respect to topologies $\Nis^X_Z$ and $\underline{\Nis}^X_Z$.

% Recall from \Cref{sect:not:PreShaeves}
% that $F_T=p^*(F)$ denotes the inverse image of 
% a presheaf $F$ over a base scheme $S$
% % $F$ 
% along structure morphisms $p\colon T\to S$
% for 
% % a presheaf $F$ over a base scheme $S$ and 
% an $S$-scheme $T$.
% We let us write $H_{\Nis^X_Z}^*(X,F)$ for $H_{\Nis^X_Z}^*(X,F_X)$.
Recall %from \Cref{sect:not:PreShaeves}
that $F_T=p^*(F)$ is 
the inverse image of 
$F\in\Pre(\Sch_X)$ along a given morphism $p\colon T\to X$,
and $H_{\nu}^*(T,F)=H_{\nu}^*(T,F_T)$,
see \Cref{sect:not:PreShaeves}.
% We let us write $H_{\nu}^*(T,F)$ for $H_{\nu}^*(T,F_T)$.

\begin{lemma}\label{lm:XNisZNisZ(UfinX)cong}
% Given $F\in\mathbb Z\Pre(X)$,
% denote by $F_{\Nis^X_Z}$ the ${\Nis^X_Z}$-sheafification.
For any $T\in \Sch_X$ finite over $X$
and any additive presheaf of abelian groups $F\in \mathbb Z\Pre(X)$, 
% there is a canonical isomorphism
% \[H^i_{\Nis^X_Z}(X,p_*(p^*(F)))\cong H^i_{\Nis^X_Z}(T,F),\]
\[
% $
H^i_{\Nis^X_Z}(X,p_*(F_T))\cong H^i_{\Nis^X_Z}(T,F)
% $
\]
for all $i\in\mathbb Z$,
where $p\colon T\to X$ is the canonical morphism.
% $F_{\Nis^X_Z}(U)\cong F_{\underline{\Nis}^X_Z}(U)$.
% $F_{\Nis^X_Z}(U)\cong F_{\underline{\Nis}^X_Z}(U)$.
\end{lemma}
\begin{proof}
%The embedding of topologies $\underline{\Nis}^X_Z\to \Nis^X_Z$ leads to the spectral sequence
% \[
Consider the spectral sequence
$
H^i_{\Nis^X_Z}(X,Rp_*^j(F_T))\Rightarrow H_{\Nis^X_Z}^{i+j}(T,F)
$, %\]
where $Rp_*^j(F_T)$ 
is the $\Nis^X_Z$-sheafifiticaton of the presheaf
$H^j_{\Nis^X_Z}(-\times_X T,F)$
on $\Sch_X$.
It follows by \Cref{lm:NisZXNisZpoints}(2) %, see \eqref{eq:NisZpoints},
that 
$\underline{H}^j_{\Nis^X_Z}(T\times_{X} U^h_v,F)=0$
for any $U\in\Sch_X$, $v\in U$,
because
\[T\times_{X} U^h_v=\coprod_{y\in T\times_X v}(T\times_X U)^h_y\]
since $T$ is finite over $X$. 
On the other hand, 
$\underline{H}^j_{\Nis^X_Z}(T\times_{X} (X-Z),F)=0$
since $\Nis^X_Z$ is trivial on $\Sch_{X-Z}$. 
Thus by \Cref{lm:NisZXNisZpoints}(2), %, see \eqref{eq:XNisZpoints},
\[\begin{cases}
    Rp_*^j(p^*(F))=0& \text{ for }j>0,\\
    Rp_*^j(p^*(F))=p_*(p^*(F))& \text{ for }j=0.
\end{cases}\tag*{\qedhere}\] 
% $Rp_*^j(p^*(F))=0$ for $j>0$, 
% $Rp_*^j(p^*(F))=p_*(p^*(F))$ for $j=0$.
\end{proof}

\begin{lemma}\label{lm:XNisZFr}
% (1)
For a closed subscheme $Z$ in a scheme $X$, 
for any quasi-stable linear (first-order) framed presheaf $F$ over $X$, 
the presheaf 
% $\Sch_X\to \Ab$, $T\mapsto H^*_{\Nis^X_Z}(X,p_*(F_T))$,
% \[\Sch_X\to \Ab;\quad (p\colon T\to X)\mapsto H^*_{\Nis^X_Z}(X,p_*(F_T))\] 
\[\Sch_X\to \Ab;\quad T\mapsto H^*_{\Nis^X_Z}(X,p_*(F_T)),\quad p\colon T\to X,\] 
% \[\Sch_X\to \Ab;\quad T\mapsto H^*_{\Nis^X_Z}(X,p_*(F_T))\] 
% where $p\colon T\to X$ is the structure morphism,
has a canonical structure of a quasi-stable (first-order) framed presheaf,
% (2)
such that
for 
% a closed subscheme $Z$ in a scheme $X$, 
% and 
a pair of closed subschemes $Z_1\subset Z_2$
in a scheme $X$, the canonical 
natural transformations %morphism 
\begin{equation}\label{eq:HNisXZ1XpFTtoHNisXZ2XpFTtoHNisTF}
    H^*_{\Nis^X_{Z_1}}(X,p_*(F_T))\to  H^*_{\Nis^X_{Z_2}}(X,p_*(F_T)) \to H^*_{\Nis}(T,F)
\end{equation}
are morphisms of framed presheaves. 
% \[\begin{array}{lcl}
%     H^*_{\Nis^X_{Z_1}}(X,p_*(F_T))&\to & H^*_{\Nis^X_{Z_2}}(X,p_*(F_T)) \\
%     H^*_{\Nis^X_Z}(X,p_*(p^*(F))) &\to &H^*_{\Nis}(T,F)
% \end{array}\]
% agree with the framed transfers structure on the above presheaves. 
Similarly for the reduced cohomologies \eqref{eq:redcohunHtau}.
% $\underline H^*_\nu(-,F_\nu)$ 
% and the morphism agrees too.
% The same applies to the reduced cohomologies \eqref{eq:redcohunHtau}. %$\underline H^*_\nu(-,F_\nu)$.
\end{lemma}
\begin{proof}
For any $c\in \Fr_n(U_1, U_2)$ %or $c\in \nrFr_n(U_1, U_2)$ 
over $X$ 
and a morphism $\widetilde X\to X$, 
the base change along 
% the morphism
a covering 
$\widetilde X^{\times l}\to X$
for each $l\in\mathbb Z_{\geq 0}$
%there is the canonical framed correspondences 
induces 
$\widetilde c^l\in \Fr_n(U_1\times_X \widetilde X^{\times l}, U_2\times_X \widetilde X^{\times l})$.
% and
% $\widetilde c^l\in \nrFr_n(U_1\times_X \widetilde X, U_2\times_X \widetilde X^{\times l})$.
% Summarising all together for $l\in\mathbb Z_{\geq 0}$ Hence for any 
For any $F\in\ZPre(\ZF_*(X))$,
correspondences $\widetilde c^l$ % quasi-stable 
% linear framed or first-order framed presheaf $F$, 
% there is 
induce the morphism of \v{C}ech complexes 
with respect to the coverings of $U_1$ and $U_2$.
% $U_1\times_X\widetilde X\to U_1$, and $U_2\times_X \widetilde X\to U_2$.
This defines a structure of %a 
%Hence the presheaf is 
a quasi-stable linear (first-order) framed presheaf on cohomologies
natural with respect to morphisms \eqref{eq:HNisXZ1XpFTtoHNisXZ2XpFTtoHNisTF}.
% Thus the first claim is proven. The second claim follows because of the naturality of the above construction.
% For any $c\in \Fr_n(U_1, U_2)$ or $c\in \nrFr_n(U_1, U_2)$ over $X$, and a morphism $\widetilde X\to X$ 
% for each $l$, there is the canonical framed correspondences $\widetilde c\in \Fr_n(U_1\times_X \widetilde X, U_2\times_X \widetilde X^{\times l})$ and
%  $\widetilde c\in \nrFr_n(U_1\times_X \widetilde X, U_2\times_X \widetilde X^{\times l})$ given by the base change.
% Hence for any quasi-stable linear framed or first-order framed presheaf $F$, 
% there is the morphism of \v{C}ech complexes with respect to the coverings $U_1\times_X\widetilde X\to U_1$, and $U_2\times_X \widetilde X\to U_2$.
% Hence the presheaf is a quasi-stable linear (first-order) framed presheaf.
% Thus the first claim is proven. The second claim follows because of the naturality of the above construction.
% The case of $\underline H^*_\nu(-,F_\nu)$ is similar.
The claim for $\nrFr_+(S)$ and reduced cohomologies follows similarly.
% $H^*_*()$
\end{proof}

% \Cref{lm:XNisZFr} 
% % is proven
% % for cohomologies $H^*_\nu(-,F_\nu)$,
% % and
% holds 
% for the reduced cohomologies \eqref{eq:redcohunHtau}
% $\underline H^*_\nu(-,F_\nu)$ 
% because of the same argument.

% Recall that $F_X$ denote the base changes of $F$ along structure morphisms $X\to S$,
% for a presheaf $F$ over a base scheme $S$, and an $S$-scheme $X$,
% and
% We let us write $H_{\Nis^X_Z}^*(X,F)$ for $H_{\Nis^X_Z}^*(X,F_X)$,
% see \Cref{sect:not:PreShaeves}.

% \begin{lemma}\label{lm:restrNisZtopology}
% For a morphism of schemes $T\to X$, and a closed subscheme $Z$ in $X$, 
% and $F\in\mathbb Z\Pre(X)$,
% there is a canonical isomorphism
% $H_{\Nis^T_{T\times_X Z}}^*(T,F)\cong H_{\Nis^X_Z}^*(T,F)$.
% \end{lemma}
% \begin{proof}
% The claim follows since the topology $\Nis^X_Z$ restricted on $\Sch_T$ equals %the topology 
% $\Nis^T_{T\times_X Z}$.
% \end{proof}
For a morphism of schemes $T\to X$, and a closed subscheme $Z$ in $X$, 
the topology $\Nis^X_Z$ restricted on $\Sch_T$ equals %the topology 
$\Nis^T_{T\times_X Z}$.
So for each $F\in\mathbb Z\Pre(X)$,
there is a canonical isomorphism
\begin{equation}\label{lm:restrNisZtopology}
    H_{\Nis^T_{T\times_X Z}}^*(T,F)\cong H_{\Nis^X_Z}^*(T,F).
\end{equation}

\begin{definition}\label{def:pushforwardtopcahngeinvimagefinhens}
Let $c=(T,\tau,g)\in \nrFr_n((X,U), (Y,V))$ for $(X,U), (Y,V)\in \Sm^\pair_S$ over a base scheme $S$.
Let $W=Y\setminus V$, $Z= X\setminus U$.
% Recall that $F_X$ and $F_Y$ denote the base changes of $F$ along structure morphisms $X\to S$ and $Y\to S$,
% see \Cref{sect:not:PreShaeves}.
Define the homomorphism
% \[c^*\colon H^*_{\Nis^Y_{W}}(Y,F_Y)\to H^*_{\Nis^X_Z}(X,F_X)\]
\[c^*\colon H^*_{\Nis^Y_{W}}(Y,F):=H^*_{\Nis^Y_{W}}(Y,F_Y)\to H^*_{\Nis^X_Z}(X,F_X)=:H^*_{\Nis^X_Z}(X,F)\]
as the composite homomorphism
% \begin{equation}\label{diag:pushforwardtopcahngeinvimagefinhens}\xymatrix{
% &H^i_{\Nis^{Y}_{W}}(T,F_Y)\ar[d]_{}& H^i_{\Nis^{Y}_{W}}(Y,F_Y)\ar[l]^{g^*}\\
% H^i_{\Nis^X_Z}(X,p_*(p^*(F_X)))\ar[d]_{\text{Lm}\,\ref{lm:XNisZFr}}\ar[r]^{\simeq}_{\text{Lm}\,\ref{lm:XNisZNisZ(UfinX)cong}} & 
% H_{\Nis^{X}_{Z}}(T,F_X)\\
% H_{\Nis^{X}_{Z}}(X,F_X)%\ar[l]^{\simeq}_{{\text{Lm}\,\ref{lm:XNisZNisZ(UfinX)cong}}}
% ,}\end{equation}
\begin{equation}\label{diag:pushforwardtopcahngeinvimagefinhens}\xymatrix{
&H^i_{\Nis^{Y}_{W}}(T,F_Y)\ar[d]_{}& H^i_{\Nis^{Y}_{W}}(Y,F_Y)\ar[l]^{g^*}\\
H^i_{\Nis^X_Z}(X,p_*(F_T))\ar[d]_{\text{Lm}\,\ref{lm:XNisZFr}}\ar[r]^{\simeq}_{\text{Lm}\,\ref{lm:XNisZNisZ(UfinX)cong}} & 
H_{\Nis^{X}_{Z}}(T,F_X)\\
H_{\Nis^{X}_{Z}}(X,F_X)%\ar[l]^{\simeq}_{{\text{Lm}\,\ref{lm:XNisZNisZ(UfinX)cong}}}
,}\end{equation}
where $H_\nu(-)=H^*_\nu(-,F_\nu)$,
and the second morphism in the sequence is provided by \eqref{lm:restrNisZtopology} and the equality $c^{-1}(W)=Z$.
The same applies to $c^*\colon\underline{H}^*_{\Nis^X_Z}(X,F)\to \underline{H}^*_{\Nis^Y_W}(Y,F)$.
\end{definition}

Let 
$\nrFr_{+,\mathrm{p}}(S)$ denote 
the subcategory in $\nrFr_+(S)$ 
formed by morphisms such that 
% that 
% morphisms 
% are 
% % correspondences $c=(T,\tau,g)$ 
% such ones that 
$g\cong\id_{T_{[2]}}$,
where $T_2=Z(\mathcal I^2(T))$,
in \eqref{eq:taugTtwoAY}.
% Let
% $\Fr_0(S)$ 
% be the subcategory of framed correspondences of level zero,
% which are given by morphisms of schemes.
For any $c=(T,\tau,g)\in \nrFr_n(X,Y)$, 
there is a decomposition 
\begin{equation}\label{eq:c=gcirccp}
c=g \circ c_p, 
\quad
g\in\Fr_0(T_{[2]},Y),
\quad
c_p=(T,\tau,\id_{T_{[2]})})\in\nrFr_n(X,T_{[2]})),
\end{equation}
and similarly for $\nrFrpair_n(S)$.

\begin{proposition}\label{prop:NisZFr}
For any quasi-stable first-order framed linear presheaf $F$ over a base scheme $S$,
the presheaf 
% \[
$(X,X-Z)\mapsto H^*_{\Nis^X_Z}(X,F)$
% \]
% where $\Nis^X_Z$ is the topology on $\Sch_X$ from \Cref{def:NisXZ},
has the canonical structure of the quasi-stable first-order framed presheaf.
The same applies to the presheaf $(X,X-Z)\mapsto \underline{H}^*_{\Nis^X_Z}(X,F)$
and also to $\A^1$-invariant framed presheaves over a base scheme $S$.
\end{proposition}
\begin{proof}
To prove the claim we need to show 
isomorphisms
%dispaly%\[
$(c_2 \circ c_1)^*\cong c_1^*\circ c_2^*$.
% \] 
% for
% $c_1\in \nrFr_{n_1}((X_0,X_0-Z_0),(X_1,X_1-Z_1))$,
% $c_2\in \nrFr_{n_2}((X_1,X_1-Z_1),(X_2,X_2-Z_2))$.
% composible framed correspondences $c_1,c_2$.
% Let 
% $\nrFr_{+,\mathrm{p}}(S)$ denote 
% the subcategory in $\nrFr_+(S)$ 
% that morphisms 
% are 
% % correspondences $c=(T,\tau,g)$ 
% such ones that $g\cong\id_{T_{[2]})}$,
% where $T_2=Z(\mathcal I^2(T))$,
% in \eqrer{eq:taugTtwoAY}.
% Let
% $\Fr_0(S)$ 
% be the subcategory of framed correspondences of level zero.
% % given by morphisms of schemes
% Since for any morphism $c=(T,\tau,g)\in \nrFr_n(X,Y)$ in $\nrFr_+(S)$, there is a decomposition 
% \[c=g \circ c_p, \quad
% c_p=(T,\tau,\id_{T_{[2]})})\in\nrFr_n(X,T_{[2]})),
% \]
% and similarly, for $\nrFr^\pair_n(S)$,
% \nrFr_n((X,X-Z),(T_2,T_2-T_2\times_X Z_2))
% Because of this 
% So
Because of \eqref{eq:c=gcirccp} it is enough to consider four cases: % of types of $c_1$ and $c_2$. 
\par (0)
For $c_1,c_2\in \nrFr_0(S)$,
the claim holds by the functoriality of the inverse images %pushforwards 
for regular morphisms of schemes.
\par (1)
For $c_1,c_2\in \nrFr_{+,\mathrm{p}}(S)$, the claim follows by \Cref{lm:XNisZFr}.
\par (2)
% The case of 
% For
% $c_1=(T,\tau,\id_{Z(\mathcal I^2(T)})\in \nrFr_n(X_1,Z(\mathcal I^2(T))$, $c_2=g\in \Fr_0(Z(\mathcal I^2(T),Y)=\Sch_S(Z(\mathcal I^2(T)),Y)$,
% the claim
% holds by \Cref{def:pushforwardtopcahngeinvimagefinhens}.
For
$c_1=(T,\tau,\id_{Z(\mathcal I^2(T)})$, $c_2=g$,
the claim
holds by \Cref{def:pushforwardtopcahngeinvimagefinhens}.
\par (3)
For 
$c_1=f\in \Sch^\pair_S((X,X-Z),(X^\prime,X^\prime-Z^\prime))$,
$c_2=(Y,\tau,\id_Y)\in \nrFr_n((X^\prime,X^\prime-Z^\prime),(Y,Y-W))$,
% % Let 
% % \[\begin{array}{ll}
% % c_1=f&\in \Fr_0((X,Z),(X^\prime,Z^\prime))=\Sch^\pair_S((X,Z),(X^\prime,Z^\prime)), \\
% % c_2=(Y,\tau,\id_Y)&\in \nrFr_n((X^\prime,Z^\prime),(Y,W))
% % .\end{array}\]
% Then the 
% The claim follows because of 
there is
the commutative diagram 
\[\xymatrix{
H^*_{\Nis^Y_W}(T,F)\ar[d] & H^*_{\Nis^Y_W}(Y,F)\ar[d]\ar[l] \\
H^*_{\Nis^X_Z}(T,F) & H^*_{\Nis^{X^\prime}_{Z^\prime}}(Y,F)\ar[l] \\
H^*_{\Nis^X_Z}(X,p^T_*(F_T))\ar[u]^{\simeq}_{\text{Lm}\,\ref{lm:XNisZNisZ(UfinX)cong}}\ar[d]^{\text{Lm}\,\ref{lm:XNisZFr}} & H^*_{\Nis^{X^\prime}_{Z^\prime}}(X,p^Y_*(F_Y))\ar[u]^{\simeq}_{\text{Lm}\,\ref{lm:XNisZNisZ(UfinX)cong}}\ar[l]\ar[d]^{\text{Lm}\,\ref{lm:XNisZFr}} \\
H^*_{\Nis^X_Z}(X,F) & H^*_{\Nis^{X^\prime}_{Z^\prime}}(X^\prime,F)\ar[l]^{f^*} 
,}\]
where $T=Y\times_{X^\prime} X$,
$p^T\colon T\to S$, $p^Y\colon Y\to S$.

The case of reduced cohomologies \eqref{eq:redcohunHtau} is similar.
The claim for $\A^1$-invariant framed presheaves follows
by \eqref{eq:Fr->nrFr} and
\Cref{lm:Fr<--nrFr}.
\end{proof}

% \begin{lemma}\label{prop:CommCohSupTrans}
% Let $V$ be a scheme over $S$, and let $U$ be a $V$-scheme.
% Let $X$ be an $S$-scheme, and $c\in \Fr_n(U,X)$ be a framed correspondence over $S$.
% Then for any linear presheaf with framed transfers $F$ over $S$,
% there is the commutative square
% \[\xymatrix{
% \underline{H}^*_{\Nis_{Y}^{U}}(U,F)\ar[d] & \underline{H}^*_{\Nis_{Z}^{X}}(X,F)\ar[l]^{c^*}\ar[d] \\ 
% \underline{H}^*_\Nis(U,F)&\underline{H}^*_{\Nis}(X,F)\ar[l]^{c^*}
% ,}\]
% where $Y=c^{-1}(Z)$.
% \end{lemma}
\begin{lemma}\label{prop:CommCohSupTrans}
% Let $U,X\in\Sch_S$, 
% and 
% $c\in \Fr_n(U,X)$ be a framed correspondence over $S$.
% Then 
Under the assumptions of \Cref{prop:NisZFr},
the homomorphism \[\underline{H}^*_{\Nis_{Z}^{X}}(X,F)\to \underline{H}^*_{\Nis}(X,F)\]defines the homomorphism of presheaves on $\Fr_+^\mathrm{pair}(S)$.
% framed presheaves on $\Sch_S^\mathrm{pair}$.
% for any closed subscheme $Z$ of $X$,
% % and $Y=c^{-1}(Z)$,
% for any linear presheaf with framed transfers $F$ over $S$,
% there is the commutative square
% \[\xymatrix{
% \underline{H}^*_{\Nis_{Y}^{U}}(U,F)\ar[d] & \underline{H}^*_{\Nis_{Z}^{X}}(X,F)\ar[l]^{c^*}\ar[d] \\ 
% \underline{H}^*_\Nis(U,F)&\underline{H}^*_{\Nis}(X,F)\ar[l]^{c^*}
% ,}\]
% where $Y=c^{-1}(Z)$.
\end{lemma}
% \begin{lemma}\label{prop:CommCohSupTrans}
% Let $U,X\in\Sch_S$, 
% and 
% $c\in \Fr_n(U,X)$ be a framed correspondence over $S$.
% Then 
% for any closed subscheme $Z$ of $X$,
% % and $Y=c^{-1}(Z)$,
% for any linear presheaf with framed transfers $F$ over $S$,
% there is the commutative square
% \[\xymatrix{
% \underline{H}^*_{\Nis_{Y}^{U}}(U,F)\ar[d] & \underline{H}^*_{\Nis_{Z}^{X}}(X,F)\ar[l]^{c^*}\ar[d] \\ 
% \underline{H}^*_\Nis(U,F)&\underline{H}^*_{\Nis}(X,F)\ar[l]^{c^*}
% ,}\]
% where $Y=c^{-1}(Z)$.
% \end{lemma}
\begin{proof}
The claim follows because 
% morphisms of 
natural transformations in
% diagram 
\eqref{diag:pushforwardtopcahngeinvimagefinhens} 
commute with the morphisms $H_\nu^*(-,F_\nu)\to H_{\Nis}^*(-,F_\Nis)$ for $\nu$ as in \Cref{def:pushforwardtopcahngeinvimagefinhens}.
\end{proof}

\section{Cohomologies out of infinity on $\A^1_V$.}\label{sect:outofinffinitesupSHI}

In this section, we show 
% strict homotopy invariance for ``cohomologies out of infinity'',
% namely, 
triviality of ``cohomologies out of infinity'' on the relative affine line $\A^1_V$ over a local henselian scheme $V$,
where
``cohomologies out of infinity''
means
% are defined as 
``the universal cohomologies trivial at infinity''
defined as cohomologies 
with respect to the topology $\Fin_V$, see \Cref{def:FinS}.
% on $\A^1_V$.
% in the sense of cohomologies on $\A^1_V$ with respect to the topology $\Fin_V$.
%
% In this section, we show a partial strict homotopy invariance result.
%
% In this section, we show a strict homotopy invariance for ``cohomologies out of infinity''.
% % for $\A^1$-invariant quasi-stable linear framed presheaves.
% Namely, 
% \Cref{th:FinSuppLocHensbaseafflinetriviality},
% proves
% %we show 
% that 
% ``cohomologies out of infinity''
% on the relative affine line $\A^1_V$ over a local henselian scheme $V$ are trivial.
% Here
% ``cohomologies out of infinity''
% means
% % are defined as 
% ``the most universal cohomologies trivial at infinity''
% in the following sense.
% Namely, we define ``cohomologies out of the infinity'' as ``the universal cohomologies that are trivial at infinity'',
% and 
% show that such cohomologies on the relative affine line $\A^1_V$ over a local henselian scheme $V$ are trivial, \Cref{th:FinSuppLocHensbaseafflinetriviality}.
% By the definition ``cohomologies out of infinity'' %are ``the universal cohomologies that are trivial at infinity'', 
% are cohomologies with respect to the union of subtopologies $\Nis^{\A^1\times V}_Z$, defined in \Cref{sect:(un)NisXZ} for all closed subschemes $Z$ finite over $V$.

\begin{definition}\label{def:PropS}
Let $X$ be a scheme over a scheme $S$. Denote by $\mathrm{Prop}_S^X$ 
the 
minimal
topology on $\Sch_{X}$
% generated by 
that contains
topologies $\Nis_{Z}^{X}$ for all closed subschemes $Z$ in $X$ proper over $S$, see \Cref{def:NisXZ}.
\end{definition}

\begin{example}\label{def:FinS}
Let $S$ be a scheme. Define $\Fin_S=\mathrm{Prop}^{\A^1\times S}_S$, which is
% \] 
the minimal topology on $\Sch_{\A^1\times S}$
that contains
% generated by 
topologies $\Nis_{Z}^{\A^1\times S}$ for all closed subschemes $Z$ in $\A^1\times S$ finite over $S$.
\end{example}

% By \Cref{th:FinSuppLocHensbaseafflinetriviality}
% % In \Cref{th:FinSuppLocHensbaseafflinetriviality}, we prove that
% for a local henselian scheme $V$,
% \begin{equation}\label{eq:cohNisFinunNis}
%     \underline{H}^*_{\Fin_V}(V\times\A^1,F)\simeq 0.
% \end{equation}

\begin{remark}\label{remark:ExLpArg}
The basic principles of the arguments 
in this section
up to ``cohomolgies out of infinity"
ware introduced in \cite{Voe-hty-inv} 
% came originally from \cite{Voe-hty-inv} 
% 
% being 
and 
reworked 
% combined 
%later 
with 
% the interpretation
% in terms of 
the use of open pairs categories
% combined with the use of open pairs 
% correpondences 
suggested in \cite{hiWt} %,phiWtshv,shvNhiWt
and %the technique based on 
Serre's theorem \cite[III, Corollary 10.7]{Hartshorne-AlG} 
% and the technique of sections of linear bundles 
suggested to the author by I. Panin.
% Up to % the concidering of
% ``cohomolgies out of infinity"
% % i.e. % the topology $\Fin_S$,
% Up to ``cohomolgies out of infinity"
% the basic principles of the arguments 
% came originally from \cite{Voe-hty-inv} 
% being 
% combined %later 
% with the use of open pairs 
% % combined with the use of open pairs 
% % correpondences 
% suggested in \cite{hiWt} %,phiWtshv,shvNhiWt
% and the technique of sections of linear bundles suggested by I. Panin.
% Up to ``cohomolgies out of infinity"
% the basic principles of the arguments 
% came originally from \cite{Voe-hty-inv} 
% being 
% interpreted in %later 
% terms of open pairs categories
% % combined with the use of open pairs 
% % correpondences 
% suggested in \cite{hiWt} %,phiWtshv,shvNhiWt
% and 
% the technique of sections of linear bundles suggested by I. Panin.
% The %proof of 
% arguments of 
% \Cref{lm:pairZarExA1Vhalffin,lm:ZarExA1VhalffinF} 
% summarises 
% various ones 
% from \cite{Framed} and from \cite{GWStrHomInv,ccorrs},
% % Though the argument of \Cref{lm:pairZarExA1Vhalffin} is not like in \cite{Framed}, it is like in \cite{GWStrHomInv,ccorrs},
% and a reader familiar with the latter sources would not meet anything %something %essentially 
% novel. 
% \Cref{lm:pairZarExA1Vhalffin}
% and 
% \Cref{lm:ZarExA1VhalffinF}

The 
argument strategy 
for excision isomorphisms
used here 
came from \cite{GWStrHomInv,ccorrs},
and constructions of correspondences from \cite{Framed}.

% The proof of excision isomorphisms
% % combines 
% uses
% constructions of correspondences from \cite{Framed} and 
% the 
% %strategy from 
% % scheme of the 
% argument
% strategy 
% from \cite{GWStrHomInv,ccorrs}.
%% Though the argument of \Cref{lm:pairZarExA1Vhalffin} is not like in \cite{Framed}, it is like in \cite{GWStrHomInv,ccorrs},
% A reader familiar with the latter sources would not meet anything %something %essentially 
% novel. 
\end{remark}

\subsection{Excision isomorphisms}

\begin{proposition}[\'Etale excision \protect{\cite[Theorem 3.10]{surj-etale-exc}}]\label{prop:EtExcision}
Let $\widetilde X$ be an \'etale neighbourhood of a closed subscheme $Z$ in an essentially smooth local henselian scheme $X$ over a field $k$.
Then for any quasi-stable framed linear presheaf $F$, 
the following square is Cartesian and coCartesian
\[\xymatrix{
F(\widetilde X-Z)& F(\widetilde X)\ar[l]\\
F(X-Z)\ar[u]& F(X).\ar[l]\ar[u]
}\]
% the square 
% \[\xymatrix{
% F(\widetilde X-Z)& F(\widetilde X)\ar[l]\\
% F(X-Z)\ar[u]& F(X)\ar[l]\ar[u]
% }\]
% is Cartesian and coCartesian.
\end{proposition}

\begin{proposition}(Zariski excision with finite support)\label{lm:pairZarExA1Vhalffin}
Let $V$ be a local scheme. Let $Z_0$, $Z_1$ be closed subschemes in $\A^1_V$, 
such that
$Z_0\cap Z_1=\emptyset$,
% Let $Z_0,Z_1\to \A^1_V$ be closed immersions, 
% $Z_0\times_{\A^1_V}Z_1=\emptyset$.
and $Z_0$ is finite over $V$.
Then for any quasi-stable framed linear presheaf $G\colon \ZF_*^\pair(V)\to \Ab$, 
% Then for any quasi-stable linear framed presheaf $G\colon \ZF_*^\pair(V)\to \Ab$, see 
% % \Cref{def:Smpair} and \Cref{def:Frpair_+} 
% \Cref{def:explFr},
% \Cref{def:ZFstar},
% \Cref{sect:ZFr}
% and \Cref{def:Cpair}
% for the category $\ZF_*^\pair(V)$, %of open pairs,
the morphism 
\begin{equation}\label{eq:j*GpairsA1V}j^*\colon G(\A^1_V,\A^1_V-Z_0)\to G(\A^1_V-Z_1,\A^1_V-Z_1-Z_0),\end{equation}
induced by the canonical open immersion of pairs
%display%
\[j\colon (\A^1_V-Z_1,\A^1_V-Z_1-Z_0)\to (\A^1_V, \A^1_V-Z_0)\]
is an isomorphism.
See 
\Cref{def:explFr},
\Cref{def:ZFstar},
\Cref{sect:ZFr}
and \Cref{def:Cpair}
for the category $\ZF_*^\pair(V)$.
\end{proposition}
\begin{proof}
Denote by $\calO(1)$ the canonical ample bundle on $\PP^1$ and let $t_\infty$  
denote the section that vanishing locus is $\infty$. 
% Denote by the same symbols the inverse images on $\PP^1_V$, and schemes of the form $-\times\PP^1_V$.
% 
Denote $\infty_V=\infty\times V=\PP^1_V\setminus \A^1_V$, 
and let $\ovZ_1$ be the closure of $Z_1$ in $\PP^1_V$.
Since $Z_0$ is finite over $V$ it is a closed subscheme in $\PP^1_V$.
So $Z_0\cap \ovZ_1 = Z_0\cap Z_1 = \emptyset$.
Denote by $\delta\in \Gamma(\A^1_V\times_V \PP^1_V,\mathcal O(1))$ a section that vanishing locus is the graph of the canonical immersion $\A^1_V\to \PP^1_V$. 
%
% Since $Z_0$ and $\infty_V\cup \ovZ_1$ are finite over a local scheme, 
% they are semi-local, 
% and 
By
\Cref{lm:semilocfinlocLinearbundle}
for any
% ll large enough 
$l\in \mathbb Z$,
there are 
invertible sections
\[\gamma_{1,\infty}\in \Gamma(\infty_V\cup\ovZ_1,\mathcal O(l)), \gamma_{Z_0}\in \Gamma(Z_0,\mathcal O(l-1)).\] 

We are going to show that \eqref{eq:j*GpairsA1V} has left and right inverse morphisms.

Part 1) 
Firstly, we construct the left inverse.
By \cite[III, Corollary 10.7]{Hartshorne-AlG}
for large enough $l\in \mathbb Z$,
there exist sections 
\begin{gather*}
\begin{array}{ll}
s_1\in \Gamma(\A^1\times \PP^1_V,\mathcal O(l)),&
s_0^\prime\in \Gamma(\A^1\times \PP^1_V,\mathcal O(l-1)),\\
s_1\big|_{\A^1\times Z_0}=\gamma_{Z_0}\delta\big|_{\A^1\times Z_0}, &
s_0^\prime\big|_{\A^1\times Z_0}=\gamma_{Z_0},\\
s_1\big|_{\A^1\times (\infty_V\cup\ovZ_1)}=\gamma_{1,\infty}, &
s_0^\prime\big|_{\A^1\times \infty_V}=\gamma_{1,\infty}\delta^{-1}\big|_{\A^1\times \infty_V},\\
& s_0^\prime\big|_{Z(\delta)}=t_\infty^{l-1}.
\end{array}
\end{gather*}
% \[s_1\in \Gamma(\A^1\times \PP^1_V,\mathcal O(l)),\;
% s_0^\prime\in \Gamma(\A^1\times \PP^1_V,\mathcal O(l-1))\]
% such that
% \[\begin{array}{ll}
% s_1\big|_{\A^1\times Z_0}=\gamma_{Z_0}\delta\big|_{\A^1\times Z_0}, &
% s_0^\prime\big|_{\A^1\times Z_0}=\gamma_{Z_0},\\
% s_1\big|_{\A^1\times (\infty_V\cup\ovZ_1)}=\gamma_{1,\infty}, &
% s_0^\prime\big|_{\A^1\times \infty_V}=\gamma_{1,\infty}\delta^{-1}\big|_{\A^1\times (\infty_V)},\\
% & s_0^\prime\big|_{Z(\delta)}=t_\infty^{l-1}.
% \end{array}\]
Define the section $s=s_1\lambda+s_0(1-\lambda)$ and morphisms in $\Fr_+^\mathrm{pair}(V)$
% framed correspondences of open pairs over $V$
\[\begin{array}{lll}
c &=& (Z(s),\A^1\times \A^1_V\times_V \A^1_V, s/t_\infty^l, g),\\&\in& 
\Fr_1^\mathrm{pair}( \A^1\times (\A^1_V,\A^1_V-Z_0) , (\A^1_V, \A^1_V-Z_0) ),\\ 
c^\prime_0&=& (Z(s^\prime),0\times \A^1_V\times_V (\A^1_V-Z_0), s_0^\prime/t_\infty^l, g_0^\prime)\\&\in& 
\Fr_1^\mathrm{pair}(0\times (\A^1_V,\A^1_V-Z_0), (\A^1_V-Z_0, \A^1_V-Z_0),\\
\tilde c_1&=& (Z(s_1),1\times \A^1_V\times_V (\A^1_V-Z_1), s_1/t_\infty^l, \tilde g_1)\\&\in& 
\Fr_1^\mathrm{pair}(1\times (\A^1_V,\A^1_V-Z_0), (\A^1_V-Z_1,\A^1_V-Z_1-Z_0)),
\end{array}\]
where $g$, $g_0^\prime$ and $\tilde g_1$ are given by the projections to the second multiplicand.
Define
\[\begin{array}{lll}
\id^\tau &=& 
(Z(\delta),0\times\A^1_V\times_V\A^1_V, s_0/t_\infty^l,g_0)\\
&\in& \Fr_1^\mathrm{pair}(0\times(\A^1_V,\A^1_V-Z_0),(\A^1_V,\A^1_V-Z_0)),
\end{array}\]
$g_0$ is induced by the projection to the second multiplicand.
Then the following equalities hold
\[c\circ i_1 = j \circ \tilde c_1, c\circ i_0 = \id^\tau + e\circ c_0^\prime\]
in $\ZF_*^\pair(V)$,
where 
\begin{gather*}
\begin{array}{lll}
i_0,i_1\colon (\A^1_V,\A^1_V-Z_0)&\to& \A^1\times (\A^1_V,\A^1_V-Z_0),
\end{array}\\
% \begin{array}{lll}
% % i_*\colon *\times (\A^1_V,\A^1_V-Z_0)&\to& \A^1\times (\A^1_V,\A^1_V-Z_0),\\
% i_0\colon 0\times (\A^1_V,\A^1_V-Z_0)&\to& \A^1\times (\A^1_V,\A^1_V-Z_0),\\
% i_1\colon 1\times (\A^1_V,\A^1_V-Z_0)&\to& \A^1\times (\A^1_V,\A^1_V-Z_0),
% \end{array}\\
e\colon (\A^1_V-Z_0, \A^1_V-Z_0)\to (\A^1_V, \A^1_V-Z_0).
\end{gather*}
Since $G$ is $\A^1$-invariant
the framed correspondence $\id^\tau$ induces the same endomorphism on $G(\A^1_V,\A^1_V-Z_0)$
as $\sigma_{(\A^1_V,\A^1_V-Z_0)}$ because of the $\A^1$-homotopy
\begin{multline}\label{eq:idtauhomotopysigma}
    (Z(\delta),(s_0\lambda + \delta t_0^{l-1})/t_\infty^l,g_0\circ \pr)\in \\
    \Fr_1^\mathrm{pair}(\A^1\times(\A^1_V,\A^1_V-Z_0),(\A^1_V,\A^1_V-Z_0)),
\end{multline}
% \[(Z(\delta),(s_0\lambda + \delta t_0^{l-1})/t_\infty^l,g_0\circ \pr)\in \Fr_1^\mathrm{pair}(\A^1\times(\A^1_V,\A^1_V-Z_0),(\A^1_V,\A^1_V-Z_0)),\]
where $\pr\colon \A^1\times(\A^1_V,\A^1_V-Z_0)\to (\A^1_V,\A^1_V-Z_0)$.
Since $G$ is additive, we have $i_0^* c^* = (\id^\tau)^*+ (e c_0^\prime)^*$.
The morphism $c^\prime$ is trivial morphism in the category of pairs $\Fr_+^\pair(V)$.
% , 
% and induces consequently the trivial morphism $(c^\prime)^*$ on the values of the presheaf $G$.
Hence %Thus
\[c_1^* j^* = (\id^\tau)^*+ (e c_0^\prime)^*=(\id^\tau)^*=\sigma^*_{(\A^1_V,\A^1_V-Z_0)}.\]
Since $G$ is quasi-stable,
the morphism \eqref{eq:j*GpairsA1V} has the left inverse.
% Since the latter morphism is an automorphism, because $G$ is quasi-stable,
% the morphism \eqref{eq:j*GpairsA1V} has the left inverse.

Part 2)
Next, we construct the right inverse to \eqref{eq:j*GpairsA1V}.
Define $V_1=\A^1_V - Z_1$.
By \cite[III, Corollary 10.7]{Hartshorne-AlG},
for large enough $l\in \mathbb Z$,
there exist sections 
\begin{gather*}
\begin{array}{ll}
s_1\in \Gamma(\A^1\times \PP^1_V,\mathcal O(l)),&
s_0^\prime\in \Gamma(V_1\times \PP^1_V,\mathcal O(l-1)),\\
s_1\big|_{\A^1\times_V Z_0}=\gamma_{Z_0}\delta\big|_{\A^1\times Z_0}, &
s_0^\prime\big|_{V_1\times_V Z_0}=\gamma_{Z_0},\\
s_1\big|_{\A^1\times (\infty_V\cup\ovZ_1)}=\gamma_{1,\infty}, &
s_0^\prime\big|_{V_1\times_V (\infty_V\cup\ovZ_1)}=\gamma_{1,\infty}\delta^{-1}\big|_{V_1\times_V (\infty_V\cup\ovZ_1)},\\
& s_0^\prime\big|_{Z(\delta)}=t_\infty^{l-1}.
\end{array}
\end{gather*}
% \[s_1\in \Gamma(\A^1\times \PP^1_V,\mathcal O(l)),\;
% s_0^\prime\in \Gamma(V_1\times \PP^1_V,\mathcal O(l-1))\]
% such that
% \[\begin{array}{ll}
% s_1\big|_{\A^1\times_V Z_0}=\gamma_{Z_0}\delta\big|_{\A^1\times Z_0}, &
% s_0^\prime\big|_{V_1\times_V Z_0}=\gamma_{Z_0},\\
% s_1\big|_{\A^1\times (\infty_V\cup\ovZ_1)}=\gamma_{1,\infty}, &
% s_0^\prime\big|_{V_1\times_V (\infty_V\cup\ovZ_1)}=\gamma_{1,\infty}\delta^{-1}\big|_{V_1\times_V (\infty_V\cup\ovZ_1))},\\
% & s_0^\prime\big|_{Z(\delta)}=t_\infty^{l-1}.
% \end{array}\]
Define sections 
and
framed correspondences of open pairs over $V$
\begin{gather*}
s_0 = \delta s_0^\prime,\quad s=s_1\lambda+s_0(1-\lambda),
\\ 
\begin{array}{lll}
c &=& (Z(s),\A^1\times V_1\times_V V_1, s/t_\infty^l, g)\\ &\in& \Fr_1^\mathrm{pair}( \A^1\times (V_1,V_1-Z_0) , (V_1, V_1-Z_0) ),\\ 
c^\prime_0&=& (Z(s^\prime),0\times V_1\times_V (V_1-Z_0), s_0^\prime/t_\infty^l, g_0^\prime)\\ &\in& \Fr_1^\mathrm{pair}(0\times (V_1,V_1-Z_0), (V_1-Z_0, V_1-Z_0),\\
\tilde c_1&=& (Z(s_1),1\times \A^1_V\times_V V_1, s_1/t_\infty^l, g)\\ &\in& \Fr_1^\mathrm{pair}(1\times (\A^1_V,\A^1_V-Z_0), (V_1,V_1-Z_0)),
\end{array}
\end{gather*}
% \[s_0 = \delta s_0^\prime,\quad s=s_1\lambda+s_0(1-\lambda),\] 
% and
% framed correspondences of open pairs over $V$
% \[\begin{array}{lll}
% c &=& (Z(s),\A^1\times V_1\times_V V_1, s/t_\infty^l, g)\\ &\in& \Fr_1( \A^1\times (V_1,V_1-Z_0) , (V_1, V_1-Z_0) ),\\ 
% c^\prime_0&=& (Z(s^\prime),0\times V_1\times_V (V_1-Z_0), s_0^\prime/t_\infty^l, g_0^\prime)\\ &\in& \Fr_1(0\times (V_1,V_1-Z_0), (V_1-Z_0, V_1-Z_0),\\
% \tilde c_1&=& (Z(s_1),1\times \A^1_V\times_V V_1, s_1/t_\infty^l, g)\\ &\in& \Fr_1(1\times (\A^1_V,\A^1_V-Z_0), (V_1,V_1-Z_0)),
% \end{array}\]
where $g$, $g_0^\prime$ and $\tilde g_1$ are 
%given 
% induced
% by the 
the morphisms 
%projections on 
to the right side multiplicands $V_1$ and $V_1-Z_0$.
Define
\[\begin{array}{lll}
\id^\tau &=& (Z(\delta),0\times V_1\times_V V_1, s_0,g_0)\\
&\in& \Fr_1^\mathrm{pair}(1\times(V_1,V_1-Z_0),(V_1,V_1-Z_0)),
\end{array}\]
% $
% \id^\tau = (Z(\delta),0\times V_1\times_V V_1, s_0,g_0)\in \Fr_1(1\times(V_1,V_1-Z_0),(V_1,V_1-Z_0)),
% $
where
$g_0$ is the morphism %induced by the projection 
to the right side %multiplicand 
$V_1$.
% Then the following equalities hold in $\ZF_*^\pair(V)$
% \[c\circ i_1 = j \circ \tilde c_1,\; c\circ i_0 = \id^\tau + e c_0^\prime\]
Then $c\circ i_1 = j \circ \tilde c_1$ and $c\circ i_0 = \id^\tau + e c_0^\prime$ in $\ZF_*^\pair(V)$,
% where 
\[\begin{array}{lll}
i_0,i_1\colon (V,V_1-Z_0)&\to& \A^1\times (V_1,V_1-Z_0),\\
e\colon (V_1-Z_0, V_1-Z_0)&\to& (V_1,V_1-Z_0).
\end{array}\]
% \[\begin{array}{lll}
% i_0\colon 0\times (V,V_1-Z_0)&\to& \A^1\times (V_1,V_1-Z_0),\\
% i_1\colon 1\times (V,V_1-Z_0)&\to& \A^1\times (V_1,V_1-Z_0),\\
% e\colon (V_1-Z_0, V_1-Z_0)&\to& (V_1,V_1-Z_0).
% \end{array}\]
% Since $G$ is $\A^1$-invariant,
% %t
The framed correspondences 
$\id^\tau$ and $\sigma_{(\A^1_V,\A^1_V-Z_0)}$
induce the same endomorphism on $G(V_1,V_1-Z_0)$
since $G$ is $\A^1$-invariant,
and 
%the framed correspondence 
because of the $\A^1$-homotopy
% $\id^\tau$ 
% induces the same endomorphism on $G(V_1,V_1-Z_0)$
% as $\sigma_{(\A^1_V,\A^1_V-Z_0)}$ because of the $\A^1$-homotopy
\[(Z(\delta),(s_0\lambda + \delta t_0^{l-1})/t_\infty^l,g_0\circ \pr)\in \Fr_1^\mathrm{pair}(\A^1\times(V_1,V_1-Z_0),(V_1,V_1-Z_0)),\]
where $\pr\colon \A^1\times(V_1,V_1-Z_0)\to (V_1,V_1-Z_0)$.
Since $G$ is additive, we have $i_0^* c^* = (\id^\tau)^*+ (e c_0^\prime)^*$.
The object $(V_1-Z_0,V_1-Z_0)$ is trivial in the category of pairs, so the morphism $c^\prime_0$ is trivial, and the induced morphism $(c^\prime_0)^*$ is trivial too.
Thus 
\begin{equation}\label{eq:juscousidtauussigmaus}
    j^* c_1^* = (\id^\tau)^* + (e c_0^\prime)^* = (\id^\tau)^*=\sigma^*_{(V_1,V_1-Z_0)}
,\end{equation}
where the right side equality follows because of an
$\A^1$-homotopy of framed correspondences similar to 
\eqref{eq:idtauhomotopysigma}.
Since $G$ is quasi-stable, 
% the latter morphism above 
% the morphism 
\eqref{eq:juscousidtauussigmaus}
is an automorphism, 
and consequently,
% hence 
% the morphism 
\eqref{eq:j*GpairsA1V} has a right inverse.
\end{proof}

\begin{corollary}\label{lm:ZarExA1Vhalffin}
Let $Z\subset \A^1_V$ be finite over 
a local scheme
$V$,
and $X^\prime\subset X$ be a pair of
open neighbourhoods of $Z$.
% 
% Let $V$ be a local scheme.
% Let $Z\subset \A^1_V$ be finite over $V$,
% and $X^\prime\subset X$ be a pair of
% open neighbourhoods of $Z$.
% and $X^\prime$ and $X$ be open neighbourhoods of $Z$, $X^\prime\subset X$.
% Let $V$ be a local scheme. 
% Let $Z$ be a closed subscheme in $\A^1_V$ finite over $V$.
% Let $X^\prime$ and $X$ be open neighbourhoods of $Z$ in $\A^1_V$ such that $X^\prime\subset X$.
%
Then for any $\A^1$-invariant quasi-stable framed linear presheaf 
%display% \[G\colon \ZF_*^\pair(V)\to \Ab,\]
%
$G\colon \ZF_*^\pair(V)\to \Ab$,
the canonical morphism 
\begin{equation*}G(X^\prime,X^\prime-Z)\to G(X,X-Z),\end{equation*}
is an isomorphism,
see 
% \Cref{def:Smpair} and \Cref{def:Frpair_+,def:ZFstar} 
\Cref{def:explFr},
\Cref{def:ZFstar},
\Cref{sect:ZFr}
and \Cref{def:Cpair}.
% for the category %of open pairs 
% $\ZF_*^\pair(V)$.
\end{corollary}
\begin{proof}
By \Cref{lm:pairZarExA1Vhalffin}, the claim holds for $X=\A^1_V$.
For an arbitrary open neighbourhood $X$ of $Z$ the claim 
follows because of the composition
%display% \begin{equation*}G(X^\prime,X^\prime-Z)\simeq G(\A^1_V,\A^1_V-Z_0)\simeq G(X,X-Z).\qedhere\end{equation*}
$G(X^\prime,X^\prime-Z)\simeq G(\A^1_V,\A^1_V-Z_0)\simeq G(X,X-Z)$.
% The case of an arbitrary open neighbourhood $X$ of $Z$ follows because of the composition
% %display% \begin{equation*}G(X^\prime,X^\prime-Z)\simeq G(\A^1_V,\A^1_V-Z_0)\simeq G(X,X-Z).\qedhere\end{equation*}
% $G(X^\prime,X^\prime-Z)\simeq G(\A^1_V,\A^1_V-Z_0)\simeq G(X,X-Z)$.
\end{proof}
\begin{corollary}\label{lm:ZarExA1VhalffinF}
Under the assumptions of \Cref{lm:ZarExA1Vhalffin},
for any $\A^1$-invariant quasi-stable framed linear presheaf $F$ on $\Sm_V$,
the following
square 
is Cartesian and coCartesian
\[\xymatrix{
F(X^\prime-Z)& F(X^\prime)\ar[l]\\
F(X-Z)\ar[u]& F(X).\ar[l]\ar[u]
}\]

% the square 
% \[\xymatrix{
% F(X^\prime-Z)& F(X^\prime)\ar[l]\\
% F(X-Z)\ar[u]& F(X)\ar[l]\ar[u]
% }\]
% is Cartesian and coCartesian.
\end{corollary}
\begin{proof}
Consider the framed presheaves on 
$\ZF_*^\pair(V)$
% the category of open pairs over $V$ 
given by
\begin{equation}\label{eq:KerCokerF(X->U)}
(X,U)\mapsto\Ker(F(X)\to F(U)),\quad
(X,U)\mapsto\Coker(F(X)\to F(U))
\end{equation}
%display
% \begin{equation}\label{eq:KerCokerF(X->U)}\begin{array}{lcl}
% (X,U)&\mapsto&\Ker(F(X)\to F(U)),\\
% (X,U)&\mapsto&\Coker(F(X)\to F(U))
% \end{array}\end{equation}
By \Cref{lm:ZarExA1Vhalffin}, since 
the morphism $j^*$ induced by
the canonical morphism % open immersion of pairs
\[j\colon (\A^1_V-Z_1,\A^1_V-Z_1-Z_0)\to (\A^1_V, \A^1_V-Z_0)\]
is an isomorphism 
% on the values of $G$
%display% \begin{equation*}j^*\colon G(\A^1_V,\A^1_V-Z_0)\to G(\A^1_V-Z_1,\A^1_V-Z_1-Z_0)\end{equation*}
for any of presheaves %$G$ 
% from 
in \eqref{eq:KerCokerF(X->U)}.
So the claim follows.
% for $G$ being anyone from \eqref{eq:KerCokerF(X->U)}.
\end{proof}

\begin{corollary}\label{lm:ZarExlocA1Vhalffin}
% Let 
% $V$ be a local scheme, 
% $X$ be an open subscheme of $\A^1_V$,
% $Z$ be a closed subscheme in $X$ finite over $V$ and local
% with the closed point $z$.
Let 
$X$ be an open subscheme of $\A^1_V$
over a local scheme $V$.
Let $Z\subset X$ be finite over $V$
with the unique closed point $z$.
% $Z$ be a closed subscheme in $X$ finite over $V$ and local
% with the closed point $z$.
%
Then for any $\A^1$-invariant quasi-stable framed linear presheaf $F$ over $V$,
the square 
\[\xymatrix{
F(X_z-Z)& F(X_z)\ar[l]\\
F(X-Z)\ar[u]& F(X)\ar[l]\ar[u]
}\]
is Cartesian and coCartesian,
where $X_z$ denotes the local scheme of $X$ at $z$.
\end{corollary}
\begin{proof}
The claim follows from \Cref{lm:ZarExA1VhalffinF}.
\end{proof}

% \begin{remark}\label{remark:ExLpArg}
% The basic principles of the arguments in the section 
% came from \cite{Voe-hty-inv}. 
% % Further developed in
% % \cite{hiWt,phiWtshv,shvNhiWt},
% % . The arguments were later developed for other types of transfers, and at the same time, the scheme of some arguments and the techniques were modified.
% The use of category of open pairs appeared in works \cite{hiWt,phiWtshv,shvNhiWt} on Witt correspondences.
% The use of sections of linear bundles and Serre's theorem was suggested to the author by I. Panin.
% % The latter works did not contain the argument for excision on $\A^1_U$ over a local henselian $U$.
% % \cite{Framed} used some of these principles too, namely the category of pairs and the use of Serre's theorem for the choice sections, % and moreover \cite{Framed} used them for the proof of some excision isomorphism on the relative affine line, 
% %
% Proofs of excision isomorphisms on $\A^1\times U$ in the neighbourhood of $0
% \times U$ %for different types of transfers
% were presented 
% in
% \cite{Framed}
% for framed transfers,
% and
% \cite{GWStrHomInv} and \cite{ccorrs}
% for other types of transfers.
% The present argument 
%  % for framed transfers
% %for this case 
% is different
% from the one in \cite{Framed} and 
% coincides with the ones in 
% \cite{GWStrHomInv} and \cite{ccorrs}.
% \end{remark}

\subsection{Coverings and cohomologies}
\begin{theorem}\label{th:FinSuppLocHensbaseafflinetriviality}
Let $V$ be a local henselian scheme, and $F$ be an $\A^1$-invariant quasi-stable framed linear presheaf over $V$.
Let $X$ be an open subscheme of $V\times\A^1$.
Then \[\underline{H}^*_{\Nis^X_Z}(X,F)=0\]
for any closed subscheme $Z$ in $X$ finite over $V$.
\end{theorem}
\begin{proof}
Since $Z$ is finite over a local henselian scheme, the connected components of $Z$ are local henselian schemes.
We prove the claim by induction on the number of the connected components of $Z$.
If $Z=\emptyset$, 
% i.e. the number is zero, 
the claim is a tautology, since the topology $\Nis^X_Z$ is trivial.
Let $Z= Z_0\amalg Z_1$, where $Z_0$ is non-empty local henselian. %, and $Z_1$ has less number of connected components than $Z$.
Let $z$ be the closed point of $Z_0$.

Consider the long exact sequence
\begin{multline}\label{eq:HnisZXX-Z0X-Z1X-Z}
\dots\to \underline{H}^{n-1}_{\Nis^X_Z}(X-Z,F)\to \underline{H}^n_{\Nis^X_Z}(X,F)\to \\ 
\underline{H}^n_{\Nis^X_Z}(X-Z_0,F)\oplus \underline{H}^n_{\Nis^X_Z}(X-Z_1,F)\to\\ \underline{H}^n_{\Nis^X_Z}(X-Z,F)\to\dots
.\end{multline}
Since any $\Nis^X_Z$-covering over $X-Z$ admits a left inverse morphism, $\underline{H}^{n-1}_{\Nis^X_Z}(X-Z,F)=0$.
Next, 
% \[%display
$\underline{H}^n_{\Nis^X_Z}(X-Z_0,F)\cong\underline{H}^n_{\Nis^X_{Z_1}}(X-Z_0,F)=0$,
% \]
% where the first isomorphism holds because $(X-Z_0)\cap Z= Z_1$, 
% and the second isomorphism holds by the inductive assumption applied to the topology $\Nis^{X}_{Z_1}$.
because $(X-Z_0)\cap Z= Z_1$, 
and by the inductive assumption applied to the topology $\Nis^{X}_{Z_1}$.
% For $Y$ being $X-Z$, or $X-Z_0$, 
% \[\underline{H}^n_{\Nis^X_Z}(Y,F)=\underline{H}^n_{\Nis^X_{Z_1}}(Y,F)=0,\]
% where the first isomorphism holds because $Y\cap Z= Z_1$, and the second isomorphism holds because by the inductive assumption the claim holds for the topology $\Nis^{X}_{Z_1}$.
% Note that 
% $\underline{H}^n_{\Nis^X_Z}(X-Z_1,F)=\underline{H}^n_{\Nis^X_{Z_0}}(X-Z_1,F)$.
% Similarly, since $(X-Z_1)\cap Z= Z_0$, 
% %display
% $
% \underline{H}^n_{\Nis^X_Z}(X-Z_1,F)\cong\underline{H}^n_{\Nis^X_{Z_0}}(X-Z_1,F)
% $.
% Let $z$ be the closed point of $Z_0$.
%
Further, 
\begin{multline*}
    \underline{H}^n_{\Nis^X_Z}(X-Z_1,F)\stackrel{(0)}{\cong}
    \underline{H}^n_{\Nis^X_{Z_0}}(X-Z_1,F)\stackrel{(1)}{\cong}
    \underline{H}^n_{\Nis^X_{Z_0},Z_0}(X-Z_1,F)\stackrel{(2)}{\cong}\\
    \underline{H}^n_{\Nis^X_{Z_0},Z_0}(X_z,F)\stackrel{(3)}{\cong}
    \underline{H}^n_{\Nis^X_{Z_0},Z_0}(X^h_z,F)\stackrel{(4)}{\cong}
    \underline{H}^n_{\Nis^X_{Z_0}}(X^h_z,F)\stackrel{(5)}{=}0.
\end{multline*}
Here 
(0) holds since $(X-Z_1)\cap Z= Z_0$,
(1) holds since $\underline{H}^n_{\Nis^X_{Z_0}}(X-Z_1-Z_0,F)\cong 0$ because $X-Z_1-Z_0$ is a $\Nis^X_{Z_0}$-point,
% the topology $\Nis^X_{Z_0}$ is trivial on $X-Z_1-Z_0$,
(2) and (3) follow by %the use of 
\Cref{lm:ZarExlocA1Vhalffin}
and %the use of 
\Cref{prop:EtExcision} 
respectively
both applied with substitutions $Z=Z_0$ 
because each $\Nis^X_{Z_0}$-covering of $X-Z_1$ and $X_z$ 
has a shrinking given by  one Nisnevich square,
and
(4), (5) hold because $X^h_z$ and $X^h_z-Z_0$ are $\Nis^X_{Z_0}$-points.
\end{proof}
\begin{definition}\label{def:underNisV}
% Let $S$ be a scheme.
% Define the topology $\underline{\Nis}^S$ on $\Sch_S$ as the topology
% with coverings being %\'etale 
% morphisms of the form $\widetilde U\times_S V\to U\times_S V$
% for a Nisnevich covering $\widetilde U\to U$, $U\in\Et_S$.

Let $S$ be a scheme.
Denote by $\underline{\Nis}^S$ the topology
on $\Sch_S$ 
that is 
% the inverse image of 
% the Nisnevich topology
% on $V\in \Et_S$,
% i.e. 
the weakest topology
such that the functor 
$\Et_S\to\Sch_S$
is continuous.
% Define the topology $\underline{\Nis}^S$ on $\Sch_S$ as the weakest topology
% % subtopology of $\Nis^X$
% such that  
% % all 
% morphisms $\widetilde U\to U$
% % for a Nisnevich covering $\widetilde U\to U$, $U\in\Et_S$.
% such that for each $v\in V$, $V\in \Et_S$, there is a dashed arrow in the triangle
% \begin{equation}\label{eq:UXhpvlifting}\xymatrix{
% &\widetilde U\ar[d] \\
% U\times_S V_v^h\ar@{-->}[ru]\ar[r]&U,
% }\end{equation}
% where $p\colon U\to S$ is the structure morphism,
% are $\underline{\Nis}^S$-coverings.
\end{definition}
% It is follows from the definition that $\underline{\Nis}^X$ is a subtopology of $\Nis^X$.

% % \begin{definition}
% % Define \emph{the topology $\underline{\Nis}^X_Z$} on $\Sm_X$ as the intersection $\underline{\Nis}^X \cap \Nis^X_Z$.\end{definition}

\begin{corollary}\label{cor:NisFalongA1hominv}
Let $S$ be a scheme.
Consider the topology $\underline\Nis^S\cup\Fin_S$ on $\Sm_{S\times\A^1}$, see \Cref{def:underNisV} and \Cref{def:FinS}.
% Consider the topology $\underline\Nis^S\cup\Fin_S$ on $\Sm_{S\times\A^1}$ generated by $\underline\Nis^S$ and $\Fin_S$, see \Cref{def:underNisV} and \Cref{def:FinS}.
Then the projection $S\times\A^1\to S$ induces a natural isomorphism 
\begin{equation}\label{eq:cohNisFinunNis}
    H^*_{\underline\Nis^S\cup\Fin_S}(S\times\A^1,F)\simeq H^*_\Nis(S,F)
\end{equation}
for any $\A^1$-invariant quasi-stable framed linear presheaf $F$ on $\Sm_S$.
%todounreaded check
\end{corollary}
\begin{proof}
Consider 
the presheaves on the small \'etale site over $S$
% and
% the canonical morphism 
\[
U\mapsto h^q_\mathrm{fnis}(U)= H^q_{\Fin_U}(U\times\A^1,F),\quad
U\mapsto h^q_\mathrm{nis}(U)= H^q_{\Nis}(U\times\A^1,F).
\]
%display
% \begin{gather}
% \nonumber
% \begin{array}{lclcl}
% U&\mapsto& h^q_\mathrm{fnis}(U)&=& H^q_{\Fin_U}(U\times\A^1,F),\\%nisF
% U&\mapsto&h^q_\mathrm{nis}(U)&=& H^q_{\Nis}(U\times\A^1,F).\end{array}
% \end{gather}
%display
% \begin{gather}
% \nonumber
% \begin{array}{lclcl}
% U&\mapsto& h^q_\mathrm{fnis}(U)&=& H^q_{\Fin_U}(U\times\A^1,F)\\%nisF
% U&\mapsto&h^q_\mathrm{nis}(U)&=& H^q_{\Nis}(U\times\A^1,F),\end{array}
% \\
% \label{eq:hqfnisnis}
%     h^q_\mathrm{fnis}\to h^q_\mathrm{nis}
% .
% \end{gather}
% \[\begin{array}{lclcl}
% U&\mapsto& h^q_\mathrm{fnis}(U)&=& H^q_{\Fin_U}(U\times\A^1,F)\\%nisF
% U&\mapsto&h^q_\mathrm{nis}(U)&=& H^q_{\Nis}(U\times\A^1,F),\end{array}\]
% and
% the canonical morphism of presheaves 
% \begin{equation}\label{eq:hqfnisnis}
%     h^q_\mathrm{fnis}\to h^q_\mathrm{nis}
% .\end{equation}
% induced by the embedding $\Fin_S\to\Nis$,
The canonical morphism $h^q_\mathrm{fnis}\to h^q_\mathrm{nis}$ is a Nisnevich local isomorphism
by \Cref{th:FinSuppLocHensbaseafflinetriviality} applied to the schemes $X=\A^1\times U^h_v$ for all $v\in U$, $U\in\Et_S$. 
Then 
the induced
morphism of 
spectral sequences
\[
H^p_{\nis}(S, h^q_\mathrm{fnis}) \Rightarrow H^{p+q}_{\underline\Nis^S\cup\Fin_S}(S\times\A^1,F), \quad
H^p_{\nis}(S, h^q_\mathrm{nis}) \Rightarrow H^{p+q}_{\Nis}(S\times\A^1,F)
\]
%display
% \[
% \begin{array}{lcl}
% H^p_{\nis}(S, h^q_\mathrm{fnis}) & \Rightarrow& H^{p+q}_{\underline\Nis^S\cup\Fin_S}(S\times\A^1,F), \\
% H^p_{\nis}(S, h^q_\mathrm{nis}) & \Rightarrow& H^{p+q}_{\Nis}(S\times\A^1,F)
% \end{array}
% \]
is an isomorphism.
% induces the claimed isomorphism \eqref{eq:cohNisFinunNis}.
\end{proof}

\section{Finite support injectivity criterion}\label{sect:InjThbyFinSupHom} 

In this section, 
% we apply the result of the previous section
% to formulate some enough criterion used in the next sections
% The 
we 
formulate 
% get 
an enough criterion for the injectivity 
of the restriction homomorphism
\[H^*_\nis(X,F_\nis)\hookrightarrow H^*_\nis(X^{(0)},F_\nis)\]
for an $\A^1$-invariant quasi-stable framed linear presheaf $F$.
% and
% an irreducible scheme $X$.
% with the generic point $\eta$
% of the Nisnevich cohomologies of an $\A^1$-invariant quasi-stable linear framed presheaf $F$
% from an irreducible scheme $X$ to the generic point $\eta$. 
%
%We formulate and use it in the forms of 
% \Cref{cor:preimFinFrCintr->Inj}
% and \Cref{cor:preimFinCovFrA1hom}.

\begin{proposition}\label{prop:FrA1HomFinCovPreimage}
Let $S$ be a scheme over a field $k$, $X\in \Sm_k$, and $Z$ be a closed subscheme in $X$.
% Assume 
Suppose
there is a framed correspondence $c\in \Fr_N(S\times\A^1, X)$ such that $c^{-1}(Z)$ is finite over $S$.
Then the morphisms 
\[c_0^*, c_1^*\colon \underline{H}^i_{\Nis_{Z}^{X}}(X,F)\to \underline{H}^i_\nis(S,F),\quad i\geq 0,
\]
where
$c_0=c\circ i_0$, $c_1=c\circ i_1$, where $i_0\colon S\times 0\to S\times\A^1$, $i_1\colon S\times 1\to S\times\A^1$,
are equal for any $\A^1$-invariant quasi-stable framed linear presheaf $F$.
\end{proposition}
\begin{proof}
The claim follows because of the commutative diagram
\[
\xymatrix{
& \underline{H}^*_{\Nis_{Z}^{X}}(X,F)\ar[dl]\ar[d]\ar[dr] & \\
\underline{H}^*_{\Nis_{c^{-1}_0(Z)}^{S\times 0}} (S\times 0, F)\ar[d] & 
\underline{H}^*_{\Nis_{c^{-1}(Z)}^{S\times\A^1}} (S\times\A^1, F)\ar[l]\ar[r]\ar[d] &
\underline{H}^*_{\Nis_{c^{-1}_1(Z)}^{S\times 1}} (S\times 1, F)\ar[d]\\
\underline{H}^*_\nis (S\times 0, F) & 
\underline{H}^*_{\underline{\Nis}^S\cup\Fin_S} (S\times\A^1, F)\ar[l]_{\cong}\ar[r]^{\cong} &
\underline{H}^*_\nis (S\times 1, F)
,}
\]
see \Cref{def:NisXZ}
for the topologies in the second row, 
\Cref{def:underNisV}
for
$\underline{\Nis}^S$, and 
\Cref{def:FinS}
for
$\Fin_S$;
the arrows between the first two rows are given by \Cref{def:CohSupTransfer}; 
the arrows from the second row to the third one are inverse images for the morphisms of sites.
The commutativity of the diagram
is provided by \Cref{prop:NisZFr} and \Cref{prop:CommCohSupTrans}.
The isomorphisms in the bottom row hold by \Cref{cor:NisFalongA1hominv}.
\end{proof}

\begin{corollary}\label{cor:preimFinCovFrA1hom}
Let $S$ be a scheme over a field $k$, and $X\in \Sm_k$, and $Z$ be a closed subscheme in $X$.
Suppose that 
\begin{itemize}%[leftmargin = 10pt]
\item[(1)] 
there is a framed correspondence $c\in \ZF_N(S\times\A^1, X)$ given by the difference of explicit framed correspondences $c=c^+-c^-$
such that schemes $(c^+)^{-1}(Z)$ and $(c^-)^{-1}(Z)$ are finite over $S$. 
\item[(2)]
there is $\widetilde c_1\in \ZF_N(S,X-Z)$ such that
$c_1=j\circ \widetilde c_1$,
where 
$c_1=c\circ i_1$, $i_1\colon S\to S\times\A^1$ is the unit section,
$j\colon X-Z\to X$ is the canonical immersion.
\end{itemize}
Then for any $\A^1$-invariant quasi-stable framed linear presheaf $F$ on $\Sm_k$
the morphisms
\[c_0^*\colon \underline{H}^i_{\Nis^X_Z}(X,F)\to \underline{H}^i_\nis(X,F)\]
induced by $c_0$ vanish for all $i\geq 0$.
\end{corollary}
\begin{proof}
The claim follows since the homomorphism
$j^*\colon \underline{H}^i_{\Nis^X_Z}(X,F)\to \underline{H}^i_{\nis}(X-Z,F)$ given by the inverse image vanishes,
and $c_0^*=\widetilde c_1^*j^*$ by \Cref{prop:FrA1HomFinCovPreimage}.
\end{proof}
\begin{corollary}\label{cor:preimFinFrCintr->Inj}
Let $X\in \Sm_k$.
If the assumption of \Cref{cor:preimFinCovFrA1hom} holds for any $Z$ in $X$ of positive codimension, then \[\underline{H}^i_\nis(X,F)\cong 0 \quad \forall i\geq 0\]
for any $\A^1$-invariant quasi-stable framed linear presheaf $F$ on $\Sm_k$.
 % for all $i\geq 0$
\end{corollary}
\begin{proof}
Recall that $\underline{H}^i_{\nis}(X,F_{\nis})=\underline{H}^i_{\Nis^X}(X,F_{\Nis^X})$.
Let $\XTheta$ denote the scheme that is the union of generic points of $X$.

Consider the topology $\Nis^{X/\XTheta}:=\bigcup_{Z}\Nis^X_Z$, where $Z$ runes over all  closed subschemes $Z$ in $X$ of positive codimension.
Then $\underline{H}^i_{\Nis^{X/\XTheta}}(X,F_{\Nis^{X/\XTheta}})\cong \varinjlim_Z \underline{H}^i_{\Nis^X_Z}(X,F_{\Nis^X_Z})$.
By \Cref{cor:preimFinCovFrA1hom}, 
% we have that 
the canonical morphisms 
$\underline{H}^i_{\Nis^X_Z}(X,F_{\Nis^X_Z})\to \underline{H}^i_{\Nis^X}(X,F_{\Nis^X})$
are trivial for all closed subschemes $Z$ of positive codimension.
Hence  the canonical morphism \begin{equation}\label{eq:HNisXThetaHNisX}\underline{H}^i_{\Nis^{X/\XTheta}}(X,F_{\Nis^X_Z})\to \underline{H}^{i}_{\Nis^X}(X,F_{\Nis^X})\end{equation} is trivial for each $i\in\mathbb Z$.

On the other side, 
% $\Nis^{X/\XTheta}$ is
% the strongest subtopology of $\Nis^X$ 
% such that 
% for any covering $\widetilde U\to U$, $U\in \Sch_X$, 
% there is a lifting $U\times_X\XTheta\to\widetilde U$, 
% because such a lifting exists is and only if there is a lifting $U\times_X (X-Z)\to \widetilde U$ for some closed subscheme $Z$. 
% % Zariksi neighbourhood $U$ of $\XTheta$.
% Hence 
$\Nis^{X/\XTheta}$ has enough set of points that is %given by 
the union of sets of Nisnevich points and schemes of the form $U\times_X\XTheta$, $U\in \Sch_X$.
Consider 
the restrictions of $\Nis^{X/\XTheta}$ and $\Nis^{X}$ 
on %the small \'etale site 
$\Et_X$, 
and the morphism of sites given by the embedding of the topologies $(\Nis^{X/\XTheta})\big|_{\Et_X}\to (\Nis^X)\big|_{\Et_X}$.
For any $U\in \Et_X$, the scheme $U\times_X\XTheta$ has Krull dimension zero. 
Hence for any additive presheaf $F$, there is the isomorphism of presheaves $H^*_{\Nis^{X/\XTheta}}(-,F_{\Nis^{X/\XTheta}})\simeq H^*_{\Nis^{X}}(-,F_{\Nis^{X}})$ on $\Et_X$.
So the canonical morphism \eqref{eq:HNisXThetaHNisX}
is an isomorphism.

Since \eqref{eq:HNisXThetaHNisX} is trivial and is an isomorphism, $\underline{H}^i_{\Nis^X}(X,F)\cong 0$.
% Thus \eqref{eq:HNisXThetaHNisX} is trivial and is an isomorphism at the same time; 
% so %the isomorphism 
% $\underline{H}^i_{\Nis^X}(X,F)\cong 0$% follows
% .
\end{proof}

\section{Compactified correspondences over one-dimensional base schemes}\label{sect:homotpiesatinfinity}
 % and cohomoligies at $\infty$ on $\A^1_V$.
% In this section,
% we prove 
% triviality of ``cohomologies at infinity'' 
% on %the relative affine line 
% $\A^1_V$ over 
% a local essentially smooth scheme $V$ over a field $k$,
% \Cref{th:GenFifLocEssSmSchemeonedimbasetriviality},
% applying \Cref{cor:preimFinFrCintr->Inj},
% to a framed $\A^1$-homotopy ``with finite supports''
% constructed in \Cref{cor:contractinggenpointhomoveretaSmX}.

% In \Cref{cor:contractinggenpointhomoveretaSmX},
% % this section,
% we construct
% a framed $\A^1$-homotopy ``with finite supports''
% required to apply 
% \Cref{cor:preimFinFrCintr->Inj},
% to
% the generic fibres of essentially smooth local schemes over one-dimensional
% base schemes.

% \Cref{cor:contractinggenpointhomoveretaSmX}.
In this section,
we prove triviality of cohomologies  
on %the relative affine line 
the generic fibres of essentially smooth local schemes over one-dimensional
base schemes
% $\A^1_V$ over 
% a local essentially smooth scheme $V$ over a field $k$,
% \Cref{th:GenFifLocEssSmSchemeonedimbasetriviality},
applying \Cref{cor:preimFinFrCintr->Inj},
to a framed $\A^1$-homotopy ``with finite supports''
constructed in \Cref{cor:contractinggenpointhomoveretaSmX},
which moves cohomology classes to the trivial ones.
For this purpose 
we develop the technique of compactified correspondences over one-dimensional base schemes.
%
% In detail, 
% we construct a homotopy that move
% a fibre of an essentially smooth scheme to a given dense open subscheme
% that sends infinity point away from the closed complement
% to the open sub 
% \begin{remark}
%     This section concentrated the most complicated part, and is entirely new.     
% \end{remark}
% In this section, we prove the triviality of the ``cohomologies at infinity'' on the relative affine line $\A^1_V$ over a local essentially smooth scheme $V$, \Cref{th:GenFifLocEssSmSchemeonedimbasetriviality}. To do this in \Cref{cor:contractinggenpointhomoveretaSmX} we construct a framed $\A^1$-homotopy ``with finite supports'' that allows to apply \Cref{cor:preimFinFrCintr->Inj}.
% This section concentrates the

\subsection{Compactified correspondences}\label{subsect:dcopmactifiedframedcorrespondece} 

\begin{definition}\label{def:reldimFr}
A \emph{$d$-dimensional framed correspondence} $\Phi$ over a scheme $U$ 
from an affine $U$-scheme $X$ to a $U$-scheme $Y$ 
is a set of data $(S,\varphi,g)$, where
\begin{itemize}
\item[(1)]
$S$ is an $X$-scheme of relative dimension $d$ over $X$ equipped with a closed embedding $S\to \A^n_X$,
\item[(2)]
$\varphi\in \mathcal O((\A^n_X)^h_{S})^{\oplus n-d}$ is a set of regular functions such that $S= Z(\varphi)$, %and
\item[(3)]
$g\colon (\A^n_X)^h_{S}\to Y$ is a morphism of $U$-schemes.
\end{itemize}
\end{definition}

\begin{definition}\label{def:compactifiedfocusedFr}
Let $B$ a local irreducible scheme, $\dim B=1$,
$z\in B^{(1)}$, $\eta\in B^{(0)}$.
Let $V\in \Sch_B$, $V_{\underline{\eta}}=V\times_B \eta$, $T\in \Sch_\eta$. 
% Let $B$ a local irreducible one-dimensional scheme, 
% $z\in B^{(1)}$, $\eta\in B^{(0)}$.
% Let $B$ a scheme, $\dim B=1$,
% $z\in B^{(1)}$, $\eta\in B^{(0)}$,
% $B=\overline{\eta}$, $\eta=B-z$.
% $z\in B$ be the closed point, $\eta\in B$ be the generic point. 
%
% Let $V\in \Sch_B$, $T\in \Sch_\eta$, and $V_{\underline{\eta}}=V\times_B \eta$. 
Assume we are given with 
\begin{itemize}
\item[(1)] 
$\ovX\in\Sch_V$, $\ovX$ is projective, $\dim_V X=d$, see \Cref{sect:not:Schemes}, %\eqref{eq:dimSX=d},
% of pure relative dimension $d$,
an ample bundle $\mathcal O(1)$ on $\ovX$, 
a section $t_\infty\in \Gamma(\ovX, \mathcal O(1))$
such that the closed subscheme $X_\infty=Z(t_\infty)$ in $\ovX$ has positive relative codimension over $V$;
% a projective equidimensional scheme $\ovX$ over $V$ of pure relative dimension $d$,
% and an ample bundle $\mathcal O(1)$ on $\ovX$ with a section $t_\infty\in \Gamma(\ovX, \mathcal O(1))$
% such that the closed subscheme $X_\infty=Z(t_\infty)$ in $\ovX$ has positive relative codimension over $V$;
\item[(2)]
a $d$-dimensional framed correspondence 
$\Phi=(X_{\underline{\eta}}, \varphi, g)$ from $V_{\underline{\eta}}$ to $T$ over $\eta$,
\item[(3)]
an isomorphism
$X_{\underline{\eta}}\cong(\ovX-X_\infty)\times_B \eta$;
\end{itemize}
The set of data
$(\ovX,X_\infty,\Phi)$ is called 
\emph{$V$-focused compactified $d$-dimensional $\eta$-framed correspondence} from $V_{\underline{\eta}}$ to $T$.
\end{definition}

\begin{definition}\label{def:smoothcontain}
Let $B$ a local irreducible one-dimensional scheme, 
$z\in B^{(1)}$, $\eta\in B^{(0)}$.
% $z\in B$ be the closed point, $\eta\in B$ be the generic point. 
Let $V\in \Sch_B$, $T\in \Sch_\eta$, and $V_{\underline{\eta}}=V\times_B \eta$. 
Let $r\colon V_{\underline{\eta}}\to T$ be a morphism of schemes.
Let $(\ovX,X_\infty,\Phi)$ be
$V$-focused compactified $d$-dimensional $\eta$-framed correspondence from $V_{\underline{\eta}}$ to $T$, $\Phi = (X_{\underline{\eta}},\varphi,g)$.

% n

We say that $(\ovX,X_\infty,\Phi)$ \emph{$V$-smoothly contains} $r$ if
there is a closed subscheme $\Gamma$ in $X=\ovX-X_\infty$
such that \begin{itemize}
\item[(1)] the canonical projection $X\to V$ is smooth over $\Gamma$ and induces an isomorphism $\gamma\colon \Gamma\cong V$,
\item[(2)] $g \circ \gamma^{-1}_{\underline{\eta}}=r$, where $\gamma_{\underline{\eta}}\colon \Gamma\times_B \eta\cong V_{\underline{\eta}}$.
\end{itemize}
\end{definition}

\subsection{Universal endo-correspondence} 
% \subsection{Universal endo-correspondence and contracting correspondence of dimension one.} 

% The main result of this subsection is \Cref{lm:DimDim1},
% that relates to the the second part of the title, namely, 
% it provides the correspondence of dimension one that is a 'curve'-homotopy that contracts some cohomology classes.

% \Cref{prop:eqXsubcurve} is used in the proof of \Cref{lm:DimDim1} together with
% \Cref{lm:XetafDimFrS}, 
% that provides a $d$-dimensional framed correspondences from a $d$-dimensioanl scheme, 
% called by the universal endo-correspondence above.
% % a $d$-dimensional framed correspondences, that is

An appropriate compactification of a smooth scheme $X$ over $B$ % $X\in\Sm_B$
provides
a focused compactified framed correspondence of dimension $d=\dim X$ from the local scheme of $X\in\Sm_B$ at $x\in X$ to $X$ itself
as
shown in
the following lemma.
% The following lemma constructs a 
% focused compactified framed correspondence of dimension $d=\dim X$ 
% from the local scheme of $X\in\Sm_B$ at $x\in X$ to $X$ itself
% provided by an appropriate compactification of $X$ over $B$.

\begin{lemma}\label{lm:XetafDimFrS}
Let $B$ be a one-dimensional irreducible local scheme, 
$z\in B^{(1)}$, $\eta\in B^{(0)}$.
% $z\in B$ be the closed point, $\eta\in B$ be the generic point.
Let $X\in \Sm_B$, $x\in X_{\underline{z}}=X\times_B z$. 
Let $d=\dim_{B}^{x} X$ be the relative dimension of $X$ at $x$ over $B$.
Let $Z$ be a closed subscheme of positive codimension in $X_{\underline{\eta}}=X\times_B \eta$.

Then there is a 
$d$-dimensional 
$\Xx$-focused compactified 
$\eta$-framed correspondence
$(\ovS,S_\infty,\Phi)$
from $(\Xx)_{\underline{\eta}}=\Xx\times_B \eta$ to $X_{\underline{\eta}}=X\times_n \eta$,
where $\Phi=(S_{\underline{\eta}}, \varphi, \val)$, 
$S_{\underline{\eta}} = S\times_B \eta$, $S = S-S_\infty$,
such that 
\begin{itemize}
\item[(1)]
$(\ovS,S_\infty,\Phi)$ $\Xx$-smoothly contains 
$\can_{\underline{\eta}}\colon(\Xx)_{\underline{\eta}}\to X_{\underline{\eta}}$;
% the canonical morphism $(\Xx)_{\underline{\eta}}\to X_{\underline{\eta}}$;
\item[(2)]
% the closure $%Y
% \ovcalZ$ %=\overline{\val^{-1}(Z)}
% of $\val^{-1}(Z)$ in $\ovS$ 
% is of positive relative codimension, 
the scheme
$\ovcalZ=\mathrm{Cl}_{\ovS}(\val^{-1}(Z))$,
% i.e. the closure of $\val^{-1}(Z)$ in $\ovS$, 
see \eqref{eq:Cl=overline},
is of positive relative codimension, 
% see\Cref{eq:Cl=overline}
% in $\ovS$ %\ovX_{\underline{\eta}}
\item[(3)]
% \begin{equation}\label{eq:codimzovYXinfclaim}
$\codim_{\ovS\times_\Xx x}(\mathrm{Cl}_{\ovS}((\ovcalZ\cap S_\infty)\times_B \eta)\times_\Xx x)\geq 2$.
% \end{equation}
% %\ovYXinf
% for
% % \[
% $F = \mathrm{Cl}_{\ovS}((\ovcalZ\cap S_\infty)\times_B \eta)$,
% % $F = \overline{(\ovcalZ\cap S_\infty)\times_B \eta}\subset\ovS$,
% % % \]
% % being the closure of %(\overline{\val^{-1}(Z)}\cap S_\infty)\times_B \eta
% % %\ovYXinf =
% % %Y
% % % where 
% % %Y
% % $(\ovcalZ\cap S_\infty)\times_B \eta$ in $\ovS$, 
% there is an inequality
% \begin{equation}\label{eq:codimzovYXinfclaim}\codim_{\ovS\times_\Xx x}(F\times_\Xx x)\geq 2.\end{equation}
\end{itemize}
\end{lemma}
\begin{proof}
Suppose $d=0$,
then
$Z=\emptyset$,
and 
the claim is trivial.
Suppose $d>0$.

Step 0)
% We can 
% assume without loss of generality that
% $X\in \SmAff_B$, $X$ is irreducible, $\dim_B X=d$, $T_{X_{\underline{\eta}}/\eta}\simeq \mathbb 1_{X_{\underline{\eta}}}^d$, where $X_{\underline{\eta}} =X\times_B\eta$, because we can
% shrink the given scheme $X$ to a Zariski neighbourhood of $x$ with such properties. 
% Shrinking $X$ to a Zariski neighbourhood of $x$, 
% we can assume without loss of generality that
% $X\in \SmAff_B$, $X$ is irreducible, $\dim_B X=d$, $T_{X_{\underline{\eta}}/\eta}\simeq \mathbb 1_{X_{\underline{\eta}}}^d$, where $X_{\underline{\eta}} =X\times_B\eta$.
% Note additionally that, 
% $X_{\underline{\eta}}$ is dense in $X$ since $X\in \Sm_B$.
We shrink $X$ to an 
irreducible smooth affine
Zariski neighbourhood of $x$
of dimension $d$
such that
$T_{X/B}\simeq \mathbb 1_{X}^d$.
% $T_{X_{\underline{\eta}}/\eta}\simeq \mathbb 1_{X_{\underline{\eta}}}^d$.
% , where $X_{\underline{\eta}} =X\times_B\eta$.
Note that $X_{\underline{\eta}}$ is dense in $X$ since $X\in \Sm_B$.
% Step 1)
We consider a closed immersion $X\to \A^{N^\prime}_B$, $N^\prime\in \mathbb Z$.
Let $\ovX=\mathrm{Cl}_{\PP^{N^\prime}_B}(X)$, 
% denote the closure of $X$ in $\PP^{N^\prime}_B$,
\[\ovX_{\underline{\eta}}=\ovX\times_B \eta, \quad \ovX_{\underline{z}}=\ovX\times_B z.\] 
Then $\ovX$ is irreducible projective $B$-scheme, $X_{\underline{\eta}}$ is a dense open subscheme.
Since $\dim B=1$, 
the scheme $\ovX$ is equidimensional over $B$
by \Cref{cor:OnedimBaseEquidimSch}.
Consider %the closed complement
$X_{\infty,\underline{\eta}}=\ovX_{\underline{\eta}} \setminus X_{\underline{\eta}}$, and
$\Xhinf=\mathrm{Cl}_{\ovX}(X_{\infty,\underline{\eta}})$.
% the closure $\Xhinf$ of $X_{\infty,\underline{\eta}}$ in $\ovX$.
Note that 
\[\Xhinf\times_B\eta = X_{\infty,\underline{\eta}} ,\quad \Xhinf\subset \ovX\setminus X.\]
Since $\dim B=1$,
%it follows 
by \Cref{cor:OnedimBaseEquidimSch}
%, that 
$\Xhinf$ is equidimensional over $B$, and %\[
$\dim_z (\Xhinf \times_B z)=\dim_\eta (\Xhinf \times_B \eta)<\dim_\eta X_\eta=d$. %\]
So 
$\codim_{\ovX/B}(\Xhinf)$.
% $\Xhinf$ is of positive relative codimension in $\ovX$ over $B$.

% Step 1)
% % Step 2)
% Since $\ovX$ is a projective scheme over $B$, there is an ample bundle $\calO(1)$ on $\ovX$.
% Let $F\subset \ovX_{\underline{z}}$ be a finite set of 
% closed points that intersects non-emptily every irreducible component of 
% $\ovX_{\underline{z}}-((\Xhinf\times_B z) \cup x)$. 
% Then for a large enough $d$, there is a section $t_\infty\in \Gamma(\ovX,\calO(d))$ 
% \[t_\infty\big|_{\Xhinf}=0,\quad t_\infty\big|_S\in \Gamma(F,\calO(d)^\times).\]
% %
% Since $\dim_B ((\Xhinf\times_B z)\cup x)<\dim_B X$, it follows that 
% $\ovX_{\underline{z}}$ equals the closure of $\ovX_{\underline{z}}-((\Xhinf\times_B z) \cup x)$. 
% Hence $F$ intersects each irreducible component of $\ovX_{\underline{z}}$.
% Then \[Z(t_\infty)\supset \Xhinf,\quad \dim_B Z(t_\infty) = n-1,\quad n=\dim_B X.\] 
% Let us
% redenote $\calO(1):=\calO(d)$
% and put \[\Xhpinf=Z(t_\infty),\quad X^\pri=\ovX-\Xhpinf.\]
Step 1)
Since $\dim_z ((\Xhinf\times_B z)\cup x)<\dim_z \ovX_{\underline{z}}$, 
it follows that 
$\ovX_{\underline{z}}-((\Xhinf\times_B z) \cup x)$
is dense in
$\ovX_{\underline{z}}$. 
% So $F$ intersects each irreducible component of $\ovX_{\underline{z}}$.
Let $F\subset \ovX_{\underline{z}}-((\Xhinf\times_B z) \cup x)$ be a finite set of 
closed points that intersects non-emptily every irreducible component of 
$\ovX_{\underline{z}}$. 
Since $\ovX$ is a projective scheme over $B$, there is 
an ample bundle $\calO(1)$ on $\ovX$
and $t_\infty\in \Gamma(\ovX,\calO(1))$ such that
\[t_\infty\big|_{\Xhinf}=0,\quad t_\infty\big|_F\in \Gamma(F,\calO(1)^\times).\]
%
% Since $\dim_B ((\Xhinf\times_B z)\cup x)<\dim_B X$, it follows that 
% $\ovX_{\underline{z}}$ equals the closure of $\ovX_{\underline{z}}-((\Xhinf\times_B z) \cup x)$. 
% Hence $F$ intersects each irreducible component of $\ovX_{\underline{z}}$.
Then 
%\[Z(t_\infty)\supset \Xhinf,\quad \dim_B Z(t_\infty) = n-1,\quad n=\dim_B X.\] 
$Z(t_\infty)\supset \Xhinf$, $\dim_B Z(t_\infty) = \dim_B X-1$.
% , $d=\dim_B X$.
Denote 
% $\calO(1):=\calO(d)$,
$\Xhpinf=Z(t_\infty)$,
$X^\pri=\ovX-\Xhpinf$.
% \]
% \[\Xhpinf=Z(t_\infty),\quad X^\pri=\ovX-\Xhpinf.\]

Step 2)
% Step 3)
% By Step (2)
% $X^\pri$ is an affine Zariski neighbourhood of $x$ in the projective $B$-scheme $\ovX$, 
% the complement
% $\Xhpinf$ is of positive relaitve codimension in $\ovX$ over $B$, 
% and $\Xhpinf\supset \Xhinf$. 
% Consequently, 
% $X^\pri$ is smooth at $x$ over $B$, 
% and 
% there is the canonical open immersion  
% \[j\colon X^\pri_{\underline{\eta}}\to X_{\underline{\eta}},\]
% where $X^\pri_{\underline{\eta}} =X^\pri\times_B \eta$. 
% By Step (1)
%
% Since
% $X^\pri$ is an affine Zariski neighbourhood of $x$ in $\ovX$, and $\Xhpinf\supset \Xhinf$,
% % Consequently, 
% $X^\pri$ is smooth at $x$ over $B$, 
Since
$\Xhpinf\supset \Xhinf$
% $X^\pri\subset X$,
and $x\in X^\pri$,
% is an affine Zariski neighbourhood of $x$ in $\ovX$, and $\Xhpinf\supset \Xhinf$,
% Consequently, 
there is the canonical open immersion  
\[
% $
% j\colon X^\pri_{\underline{\eta}}\to X_{\underline{\eta}}
j\colon X^\pri_{\underline{\eta}}=X^\pri\times_B \eta\to X_{\underline{\eta}}
% $,
,\]
and 
$X^\pri$ is smooth at $x$ over $B$.
% where $X^\pri_{\underline{\eta}} =X^\pri\times_B \eta$. 
Note %additionally 
that
$\codim_{\ovX/B}(\Xhpinf)>0$.
% $\Xhpinf$ is of positive relative codimension in $\ovX$ over $B$.
%
% Step 3)
% Step 4)
% Since the tangent bundle of $X_{\underline{\eta}}$ is trivial, 
Since $T_{X_{\underline{\eta}}/\eta}\cong\mathbb 1^{d}_{X_{\underline{\eta}}}$, 
we get
% \[
% X^\pri_{\underline{\eta}}\in \SmAff_{\underline{\eta}}, \quad
% T_{X^\pri_{\underline{\eta}}}\simeq \mathbb 1^d_{X^\pri_{\underline{\eta}}}.
% \]
$X^\pri_{\underline{\eta}}\in \SmAff_\eta$,
$T_{X^\pri_{\underline{\eta}}/\eta}\simeq \mathbb 1^d_{X^\pri_{\underline{\eta}}/\eta}$.
Then 
by \Cref{lm:framedscheme},
there is a 
$d$-dimensional framed correspondence
%display%\[
$
(X^\pri_{\underline{\eta}})^\fr=(X^\pri_{\underline{\eta}},f_{d+1},\dots, f_{N}, \tilde j)
$
% \]
where 
$\tilde j\colon (\A^N_\eta)^h_{X^\pri_{\underline{\eta}}}\to X_{\underline{\eta}}$ is such that
$\tilde j\big|_{X^\pri_{\underline{\eta}}}=j$, see \cite[Theorem I.8]{Gru}, see also %or 
\cite{Elkiksoleqhens}, \cite[Lemma 3.11]{FrRigidSmAffpairs}.
% $\tilde j\big|_{X^\pri_{\underline{\eta}}}=j\colon X^\pri_{\underline{\eta}}\to X_{\underline{\eta}}$
% is the canonical immersion.
%
So \[
(\ovX, \Xhpinf, (X^\pri_{\underline{\eta}})^{\fr} )\]
is a $d$-dimensional 
$\Xx$-focused compactified $\eta$-framed correspondence
from $\eta$ to $X_{\underline{\eta}}$.
The base change along the map $\Xx\to B$ gives 
the $d$-dimensional $\Xx$-focused compactified $\eta$-framed correspondence
$(\ovS,S_\infty, \Phi)$
from $\Xx$ to $X$,
where $\ovS=\ovX\times_B \Xx$, $S_\infty=\Xhpinf\times_B \Xx$, 
and $\Phi$ is the base change of $(X^\pri_{\underline{\eta}})^{\fr}$.
% This gives us $(\ovS,S_\infty,\Phi)$.
% proves the initial claim of the lemma. 
% In what follows, we proceed with claims (1), (2) and (3).

Step 3)
% Step 4)
% Step 5)
The composite morphism
$\Xx\to X\to \ovX$
induces 
% the section 
$\Delta\colon \Xx\to\ovS$.
Since 
% the morphism $X\to B$ is smooth,
$X\in\Sm_B$
by assumption,
the morphism $\ovS\to \Xx$ is smooth over $\Delta(\Xx)$.
Thus $(\ovS,S_\infty, \Phi)$ 
$\Xx$-smoothly contains the canonical morphism %canonical map 
$(\Xx)_{\underline{\eta}}\to X_{\underline{\eta}}$.
This proves claim (1).

Step 4)
% Step 5)
% Step 6)
%
%%%%%%%%%%%%%%%%%%%%%%%%%%%%%%%%%%%%%%%%%%%%%%%%%%%%%%%%%%%%%%%%%%%%%%%%%%%%%%%%%%%%%%%%%%%%%
% \todo{is being written}
We have 
$j^{-1}(Z)=(Z\cap X^\pri_{\underline{\eta}})$,
% We use notation $\overline{j^{-1}(Z)}=\mathrm{Cl}_{\ovX}(j^{-1}(Z))$.
%$\overline{j^{-1}(Z)}_{\underline{\eta}}=$
%display
$
\mathrm{Cl}_{\ovX}(j^{-1}(Z))\times_B\eta=\mathrm{Cl}_{\ovX_{\underline{\eta}}}(j^{-1}(Z))
$,
\begin{equation}\label{ovcalZsupsetClovSvalpreimZ}
\ovcalZ=\mathrm{Cl}_{\ovS}(\val^{-1}(Z))
\subset \overline{j^{-1}(Z)}\times_B \Xx,
\end{equation}
where $\overline{j^{-1}(Z)}=\mathrm{Cl}_{\ovX}(j^{-1}(Z))$.
Since  
$\codim_{X^\pri_{\underline{\eta}}}(j^{-1}(Z))=\codim_{X_{\underline{\eta}}}(Z)>0$,
% and 
% it follows that
%\[
$
\codim_{\ovX_{\underline{\eta}}}(\overline{j^{-1}(Z)}_{\underline{\eta}})>0
$.
%\]
Thus by \Cref{cor:OnedimBaseEquidimSch}, 
\[\codim_{\ovS/\Xx}(\ovcalZ)=\codim_{\ovX/B}(\overline{j^{-1}(Z)})\stackrel{\text{Cor. \ref{cor:OnedimBaseEquidimSch}}}{>}0.\]
% $\codim_{\ovX/B}(\overline{j^{-1}(Z)})>0$,
% % is of positive relative codimension in $\ovX$,
% and 
% %Y
% % $\ovcalZ$ is of positive codimension in $\ovS$.
% $\codim_{\ovS/\Xx}(\ovcalZ)>0$.
% is of positive codimension in $\ovS$.
This proves claim (2).
% \[\val^{-1}(Z)=(\Xx)_{\underline{\eta}}\times_B j^{-1}(Z)\]
% is of positive codimension in $S_{\underline{\eta}}$.

Step 5)
% Step 6)
% Step 7)
% Finally, we prove inequality \eqref{eq:codimzovYXinfclaim} from claim (3).
Since
% \[
$\codim_{\ovX_{\underline{\eta}}}((\overline{j^{-1}(Z)}\cap \Xhpinf)\times_B \eta)=
\codim_{\ovX_{\underline{\eta}}}((\overline{j^{-1}(Z)}_{\underline{\eta}}\cap X_{\infty,\underline{\eta}})\leq 2
$,
% \]
% and 
% \[(\overline{j^{-1}(Z)}\cap \Xhpinf)\supset\mathrm{Cl}_{\ovX}((\overline{j^{-1}(Z)}\cap \Xhpinf)\times_B \eta),\]
by \Cref{cor:OnedimBaseEquidimSch},
\[\codim_{\ovX_{\underline{z}}}(\mathrm{Cl}_{\ovX}((\overline{j^{-1}(Z)}\cap \Xhpinf)\times_B \eta)\times_B z)\leq 2.\]
Then 
% the inequality \eqref{eq:codimzovYXinfclaim} from 
claim (3)
% claim 
follows because by \eqref{ovcalZsupsetClovSvalpreimZ}
\[\ovcalZ\cap S_\infty\subset(\overline{j^{-1}(Z)}\cap \Xhpinf)\times_B \Xx.\qedhere\]
% % Step 6)
% % Step 7)
% Finally, we prove inequality \eqref{eq:codimzovYXinfclaim} from claim (3). We note that 
% % $\overline{\val^{-1}(Z)}\cap S_\infty=(\overline Z\cap \Xhpinf)\times_B (B-\eta)\times_B\Xx$.
% %Y
% \[\ovcalZ\cap S_\infty\supset(\overline{j^{-1}(Z)}\cap \Xhpinf)\times_B \Xx.\]
% Since
% % \[
% $\codim_{\ovX_{\underline{\eta}}}((\overline{j^{-1}(Z)}\cap \Xhpinf)\times_B \eta)\leq 2
% $,
% % \]
% and 
% \[(\overline{j^{-1}(Z)}\cap \Xhpinf)\supset\mathrm{Cl}_{\ovX}((\overline{j^{-1}(Z)}\cap \Xhpinf)\times_B \eta),\]
% % $\overline{j^{-1}(Z)}\cap \Xhpinf$ is the closure of $(\overline{j^{-1}(Z)}\cap \Xhpinf)\times_B \eta$ in $\ovX$,
% by \Cref{cor:OnedimBaseEquidimSch},
% \[\codim_{\ovX_{\underline{z}}}(\mathrm{Cl}_{\ovX}((\overline{Z}\cap \Xhpinf)\times_B \eta)\times_B z)\leq 2.\qedhere\]
%% \[\codim_{\ovX_{\underline{z}}}(\overline{(\overline{Z}\cap \Xhpinf)\times_B \eta)}\times_Bz\leq 2.\qedhere\]
\end{proof}

\subsection{Contracting correspondence of dimension one.} 
% \subsection{Universal endo-correspondence and contracting correspondence of dimension one.} 

% The main result here %of this subsection 
% is 
% \Cref{prop:DimDim1} 
%%, %\Cref{lm:DimDim1},
% which 
% provides 
In \Cref{prop:DimDim1},
we construct
a one-dimensional correspondence %of dimension one 
that is a ``curve''-homotopy moving
% contracting 
% curtain 
% considered 
cohomology classes on
$(\Xx)_{\underline{\eta}}$
to the trivial one.
\Cref{prop:eqXsubcurve} is used to deduce \Cref{prop:DimDim1} %\Cref{lm:DimDim1} 
from \Cref{lm:XetafDimFrS}.

% The main result of this subsection is \Cref{prop:DimDim1}, %\Cref{lm:DimDim1},
% that provides the correspondence of dimension one that is a 'curve'-homotopy that contracts some cohomology classes.
% \Cref{prop:eqXsubcurve} is used to deduce \Cref{prop:DimDim1} %\Cref{lm:DimDim1} 
% from \Cref{lm:XetafDimFrS}.
% a $d$-dimensional framed correspondences, that is
% The first part of the title reflects 
% \Cref{lm:XetafDimFrS}.
% The second part relates to \Cref{lm:DimDim1},
% that is the main result.
% The following proposition is used in the proof.

\begin{proposition}\label{prop:eqXsubcurve}
Let $B$ be a local scheme, $z$ be the closed point.
Let $Z$ be a closed subscheme, $U=B-Z$. %be the open complement.
Let $\ovX$ be a projective $B$-scheme, and $\calO(1)$ denote an ample bundle on $\ovX$. 

Let $Y\to \ovX$ and $D\to \ovX$ be closed immersions, 
% $F$ be the closure 
% of $F_U=(Y\cap D)\times_B U$
% in $\ovX$. 
$F=\mathrm{Cl}_{\ovX}(F_U)$,  
$F_U=(Y\cap D)\times_B U$. 
Let $\Delta\colon B\to X$ be a $B$-morphism. 
Suppose that
\begin{itemize}
\item[(a1)] $\ovX$ is equidimensional of pure relative dimension $d$ over $B$;
\item[(a2)]
$\codim_{\ovX\times_B z}(Y\times_B z)>0$ and $\codim_{\ovX\times_B z}(D\times_B z)>0$;
\item[(a3)] $\codim_{\ovX\times_B z}(F\times_B z)\geq 2$;
%$F\times_B Z$ is of codimension at least 2 in $\ovX\times_B Z$; %the closed fibre
\item[(a4)]
%the $B$-scheme 
$\ovX$ is $B$-smooth over the subscheme $\Delta(B)$, and $\Delta(B)\cap D=\emptyset$.
\end{itemize}

Then for large enough $l_i\in\mathbb Z$, there are sections
$s_i\in \Gamma(\ovX,\calO(l_i))$, $i=2,\dots ,d$, such that 
\begin{itemize}
\item[(c1)]
the vanishing locus $\ovC=Z(s_2,\dots, s_d)$ is of pure relative dimension one over $B$;
\item[(c2)]
$\Delta(B)\subset \ovC$, and
the morphism $\ovC\to B$ is smooth over $\Delta(B)$;
\item[(c3)]
$\ovC\cap D$ and $\ovC\cap Y$ are finite over $B$;
\item[(c4)]
$(\ovC\cap Y\cap D)\times_B U=\emptyset$.
\end{itemize}
\end{proposition}
\begin{proof}
% Step 0)
Consider the closed subscheme $\Delta(B)$ in $\ovX$
and 
% the first order thickening 
$\Delta(B)_{[2]}=Z(\calI^2(\Delta(B))$.
Since the morphism $\ovX\to B$ is smooth over $\Delta(B)$, and $\Delta(B)$ is local, 
% there are trivialisations 
$\calO(1)\big|_{\Delta(B)}\simeq \calO(\Delta(B))$, $N_{\Delta(B)/(\ovX)}\simeq \mathbb 1^d_{\Delta(B)}$.
Hence for any $d$, there are sections 
$\gamma_1,\dots, \gamma_d\in \Gamma(\Delta(B)_{[2]},\calO(d))$ such that 
$Z(\gamma_1,\dots,\gamma_d)=\Delta(B)$.

By induction 
for $i=d,\dots, 2$,
we construct sections $s_i\in\Gamma(\ovX,\calO(l_i))$
and define 
schemes 
$X_{i}$, $Y_{i}$, $D_{i}$, $F_{i}$, $S_i$, 
\begin{equation}\label{eq:s_nDelta(B)S}
s_i\in \Gamma(\ovX,\calO(l_i)),\quad s_i\big|_{\Delta(B)_{[2]}}=\gamma_i,\quad s_i\big|_{S_i}\in \calO(l_i)^\times,
\end{equation}
% Define 
%\[
$\ovX_n:=\ovX$,
$\ovX_{i-1}=Z(s_i)\subset\ovX_i$.
% \]
$Y_i=Y\cap\ovX_i$, $D_i=D\cap\ovX_i$, $F_{i,\underline{U}}=(Y\cap D\cap \ovX_i)\times_B U$, $F_i=\mathrm{Cl}_{\ovX_i}(F_{i,\underline{U}})$,
$S_i\subset \ovX_i\times_B z$ is a finite set of closed points
such that 
$S_i\cap N\neq\emptyset$
for any 
irreducible component $N$ of anyone of schemes
% Consider schemes 
\begin{equation*}\label{eq:ovXzYzDzYDz}
\ovX_i\times_B z, \quad \ovX_i\times_B z, \quad Y_i\times_B z, \quad D_i\times_B z, \quad F_i\times_B z.
\end{equation*}
% the intersection $N\cap S$ is not empty.
%
% Then \eqref{eq:s_nDelta(B)S}
% implies that
% $\dim_z (\ovX_{n-1}\times_B z)=n-1$. Hence %\begin{equation}\label{eq:dimBX(n-1)}
% $\dim_B \ovX_{n-1}=n-1$.
% % \end{equation}
% Moreover, by the same reason
% % the same condition implies that 
% \begin{equation*}\label{eq:codimYzDzYDz(n-1)}
% \begin{array}{lllll}
% \codim_{\ovX_{n-1}\times_B z }(Y\cap \ovX_{n-1})&\times_B& z&\geq& 1, \\
% \codim_{\ovX_{n-1}\times_B z} (D\cap \ovX_{n-1})&\times_B& z&\geq& 1, \\
% \codim_{\ovX_{n-1}\times_Bz } (F\cap \ovX_{n-1})&\times_B& z&\geq& 2,
% \end{array}
% \end{equation*}
% and consequently,
% % \begin{equation}\label{eq:codimYD(n-1)}
% $
% \codim_{\ovX_{n-1}\times_Bz } (F_{n-1}\cap \ovX_{n-1})\times_B z\geq 2,
% $
% % \end{equation}
% where $F_{n-1}= \overline{F_{U,n-1}}$, $F_{U,n-1}=(Y\cap D\cap \ovX_{n-1})\times_B U$.
% % By the first condition on $s_n$ in 
The middle equality in \eqref{eq:s_nDelta(B)S}
% Also we conclude that 
implies that
$\Delta(B)\subset \ovX_{i}$,
and
the morphism $\ovX_{i}\to B$ is smooth over $\Delta(B)$.
% $\Delta(B)\subset\ovX_i$, and $\ovX_{i}$ is $B$-smooth over $\Delta(B)$.
The right side condition in \eqref{eq:s_nDelta(B)S}
implies that
\begin{gather}
\label{eq:dimBX(i)}
\dim_B \ovX_{i}=i, \\
\label{eq:codimYz(i)}
\codim_{\ovX\times_Bz}(Y_i\times_B z)\geq 1,\\
\label{eq:codimDz(i)}
\codim_{\ovX\times_Bz}(D_i\times_b z)\geq 1, \\
\label{eq:codimYD(i)}
\codim_{\ovX\times_Bz}(F_i\times_B z)\geq 2,
\end{gather}
So the scheme $\ovC=Z(s_2,\dots, s_d)=\ovX_1$
satisfies claim (c2).
Claim (c1) holds by \eqref{eq:dimBX(i)}.
Claim (c3) follows by \eqref{eq:codimYz(i)}, and \eqref{eq:codimDz(i)}. 
Then $F_{1}$ is a finite over $B$, and 
since by \eqref{eq:codimYD(i)} $F_1\times_B z=\emptyset$, 
it follows by Nakayama's Lemma that $F_1=\emptyset$.
Whence 
\[(Y\cap D\cap\ovC)\times_B U=F_1\times_B U=\emptyset.\]
So claim (c4) holds.
\end{proof}

\begin{lemma}\label{lm:DimDim1}
Let $B$ be a one-dimensional irreducible local scheme, 
$z\in B^{(1)}$, $\eta\in B^{(0)}$.
% $z\in B$ be the closed point, $\eta\in B$ be the generic point.
Suppose 
there is a $d$-dimensional $\Xx$-focused 
compactified $\eta$-framed correspondence
$(\ovS,S_\infty,\Phi)$
from the $\eta$-scheme $\calV = (\Xx)_{\underline{\eta}}=\Xx\times_B \eta$ to the $\eta$-scheme $X_{\underline{\eta}}=X\times_B \eta$, where 
$\Phi= (S_{\underline{\eta}}, \varphi, \val_{S_{\underline{\eta}}})$,
$S_{\underline{\eta}}=(\ovS-S_\infty)\times_B\eta$,
 such that \begin{itemize}
\item[(a1)] $(\ovS,S_\infty,\Phi)$ smoothly contains the canonical morphism $(\Xx)_{\underline{\eta}}\to X_{\underline{\eta}}$,
\item[(a2)] 
% the closure $\ovcalZ=\overline{\val^{-1}_S(Z)}$ of the scheme $\calZ=\val^{-1}_S(Z)$ in $\ovS$ is of positive codimension over $\Xx$,
% the closure $F$ %\ovYXinf
% of $(\ovcalZ\cap S_\infty)\times_B \eta$ in $\ovS$ has codimension at least 2 in the closed fibre over $\Xx$, i.e. 
% \begin{equation}\label{eq:codimzovYXinf}\codim_{\ovS\times_\Xx x}(F\times_\Xx x)\geq 2.\end{equation} %\ovYXinf
$\codim_{\Xx}(\ovcalZ)>0$, where
$\ovcalZ=\mathrm{Cl}_{\ovS}(\calZ)$,
% \overline{\val^{-1}_S(Z)}
% of the scheme $\calZ=\val^{-1}_S(Z)$ in $\ovS$ 
$\calZ=\val^{-1}_S(Z)$,
and
\begin{equation}\label{eq:codimzovYXinf}\codim_{\ovS\times_\Xx x}(\mathrm{Cl}_{\ovS}((\ovcalZ\cap S_\infty)\times_B \eta)\times_\Xx x)\geq 2.\end{equation} %\ovYXinf
\end{itemize}

Then there is a 
one-dimensional $\Xx$-focused compactified $\eta$-framed correspondence
% \[
$(\ovC,C_\infty,\Psi)$
% \]
from $\calV = (\Xx)_{\underline{\eta}}$ to $X_{\underline{\eta}}$,
where $\Psi=(C_{\underline{\eta}}, \psi, \val_C)$, 
$C_{\underline{\eta}} = (\ovC-C_\infty)\times_B \eta$,
such that 
\begin{itemize}
\item[(c1)] $(\ovC,C_\infty,\Psi)$ $\Xx$-smoothly contains 
$\can_{\underline{\eta}}\colon(\Xx)_{\underline{\eta}}\to X_{\underline{\eta}}$,
% xrightarrow{}
% the canonical morphism $(\Xx)_{\underline{\eta}}\to X_{\underline{\eta}}$,
\item[(c2)] 
$\overline{\val^{-1}_C(Z)}=\mathrm{Cl}_{\ovC}(\val^{-1}_C(Z))$,
see \eqref{eq:Cl=overline},
is finite over $\Xx$, and 
\[
% $
(\overline{\val^{-1}_C(Z)}\cap C_\infty)\times_B \eta=\emptyset
% $
.
\]
% the closure $\overline{\val^{-1}_C(Z)}$ of $\val^{-1}_C(Z)$ in $\ovC$ is finite over $\Xx$, and 
% % \[
% $(\overline{\val^{-1}_C(Z)}\cap C_\infty)\times_B \eta=\emptyset
% $.
% % \]
\end{itemize}
\end{lemma}
\begin{proof}
Let $\calO(1)$ be an ample line bundle 
on the projective $B$-scheme $\ovS$, 
and and $t_\infty$ be a 
global
section 
% on the projective scheme $\ovS$ 
such that $Z(t_\infty) = S_\infty$ provided by \Cref{def:compactifiedfocusedFr}.
We are going to apply \Cref{prop:eqXsubcurve} 
% with 
to the scheme 
$\overline{X}:=\ovS$, and 
the subschemes 
\[D:=S_\infty, \quad
Y:=\ovcalZ, %=\overline{\val^{-1}_S(Z)}, 
\quad \Delta(B):=\Delta(\Xx).\]
where $B:=\Xx$, $Z:=\Xx\times_Bz$, $U:=\Xx\times_B\eta$.
Assumption (a1) in \Cref{prop:eqXsubcurve} holds by %the properties of the scheme $\ovS$ provided by 
\Cref{def:compactifiedfocusedFr};
(a4) in \Cref{prop:eqXsubcurve} holds by
(a1) in \Cref{lm:DimDim1} and 
\Cref{def:smoothcontain}; 
% The assumption 
(a2) in \Cref{prop:eqXsubcurve} follows from %the properties of $S_\infty$ provided by 
(a2) in \Cref{lm:DimDim1} and 
\Cref{def:compactifiedfocusedFr};
(a3) in \Cref{prop:eqXsubcurve} holds by 
\eqref{eq:codimzovYXinf}.
Then by \Cref{prop:eqXsubcurve},
there is a vector of sections $(s_2,\dots, s_d)$, where $s_i\in \Gamma(\ovS,\calO(l))$ for a large enough $l$, 
such that 
\begin{itemize}\item[(p1)]
$\dim_B\ovC=1$, where $\ovC=Z(s_2,\dots, s_d)\subset \ovS$,
% is of relative dimension one over $\Xx$, 
% the scheme $\ovC=Z(s_2,\dots, s_d)\subset \ovS$
% is of relative dimension one over $\Xx$, 
\item[(p2)]
% there is the inclusion 
$\Delta(\Xx)\subset \ovC$, and
the morphism $\ovC\to \Xx$ is smooth over $\Delta(\Xx)$,
\item[(p3)]
$\overline{\calZ}\cap \ovC$ is finite over $\Xx$,
$C_\infty=S_\infty\cap \ovC$ is finite over $\Xx$,
$(\overline{\calZ} \cap C_\infty)\times_B \eta=\emptyset$.
\end{itemize}
By assumption we are given with a closed embedding $S_{\underline{\eta}}\to \A^N_{\calV}$,
a %regular 
morphism $\val\colon (\A^N_{\calV})^h_{S_{\underline{\eta}}}\to X_{\underline{\eta}}$, 
and 
% a vector of regular functions 
$\varphi\in \mathcal O((\A^N_{\calV})^h_{S_{\underline{\eta}}})^{\oplus (N-d)}$ such that $Z(\varphi)=S_{\underline{\eta}}$.
By property (p1) %in the list above 
there are canonical embeddings 
$\ovC\to \ovS$, $C_{\underline{\eta}}\to S_{\underline{\eta}}$,
% \[\ovC\to \ovS, \quad C_{\underline{\eta}}\to S_{\underline{\eta}},\]
where $C_{\underline{\eta}}=C\times_B\eta$, $C=\ovC-C_\infty$. 
Define $\val_{C_{\underline{\eta}}}\colon (\A^N_{\calV})^h_{C_{\underline{\eta}}}\to X_{\underline{\eta}}$ as the composite morphism
\[(\A^N_\calV)^h_{C_{\underline{\eta}}}\to (\A^N_\calV)^h_{S_{\underline{\eta}}}\xrightarrow{\val_{S_{\underline{\eta}}}} X_{\underline{\eta}}.\]
Choose liftings $\varphi_j\in \mathcal O(\A^N_\calV)$ of functions $s_j/t_\infty^l\in \mathcal O(S_{\underline{\eta}})$, $j=2,\dots ,d$,
and denote by the same symbols their inverse images on $(\A^N_\calV)^h_{C_{\underline{\eta}}}$.
Then $\Psi=(C_{\underline{\eta}}, \varphi_2,\dots, \varphi_d, \varphi, \val_{C_{\underline{\eta}}})$ is a one-dimensional framed correspondence from $\calV$ to $X_{\underline{\eta}}$.
Moreover, $(\ovC,C_\infty, \Psi)$ is an $\Xx$-focused compactified $\eta$-framed correspondence because $\ovC$ is of pure relative dimension one and $C_\infty$ is finite over $\Xx$ by property (p1) and the second part of (p3) in the list above respectively.

Claim (c1) of the lemma 
follows by property (p2).
Claim (c2) follows by property (p3),
because $\overline{\val^{-1}_C(Z)}\subset \overline{\calZ}\cap\ovC$ is finite over $\Xx$,
and 
the fiber of $\overline{\val^{-1}_C(Z)}\cap C_\infty\subset \overline{\calZ}\cap C_\infty$ over $\eta$
% the fiber of $\overline{\val^{-1}_C(Z)}\cap C_\infty\subset \overline{\val^{-1}_S(Z)}\cap C_\infty$ over $\eta$
is empty.
\end{proof}

\begin{proposition}\label{prop:DimDim1}
Let $B$ be a one-dimensional irreducible local scheme, 
$z\in B^{(1)}$, $\eta\in B^{(0)}$.
% $z\in B$ be the closed point, $\eta\in B$ be the generic point.
Let $X\in \Sm_B$, $x\in X_{\underline{z}}=X\times_B z$. 
Let $d=\dim_{B}^{x} X$. %be the relative dimension of $X$ at $x$ over $B$.
Let $Z$ be a closed subscheme of positive codimension in $X_{\underline{\eta}}=X\times_B \eta$.

Then there is a 
one-dimensional $\Xx$-focused compactified $\eta$-framed correspondence
$(\ovC,C_\infty,\Psi)$
from $\calV = (\Xx)_{\underline{\eta}}$ to $X_{\underline{\eta}}$, 
where $\Psi=(C_{\underline{\eta}}, \psi, \val_C)$, $C_{\underline{\eta}} = (\ovC-C_\infty)\times_B \eta$,
such that 
\begin{itemize}
\item[(1)] $\Xx$-smoothly contains the canonical morphism $(\Xx)_{\underline{\eta}}\to X_{\underline{\eta}}$,
\item[(2)]  the closure $\overline{\val^{-1}_C(Z)}$ of $\val^{-1}_C(Z)$ in $\ovC$ is finite over $\Xx$, and 
$(\overline{\val^{-1}_C(Z)}\cap C_\infty)\times_B \eta=\emptyset$.
% \[(\overline{\val^{-1}_C(Z)}\cap C_\infty)\times_B \eta=\emptyset.\]
\end{itemize}
\end{proposition}
\begin{proof}
The claim follows by \Cref{lm:XetafDimFrS} and \Cref{lm:DimDim1}. 
\end{proof}

% \subsection{Finite correspondence homotopy, and the vanishing of cohomoligies at infinity.}
\subsection{Finite correspondence homotopy.}

\begin{proposition}\label{cor:contractinggenpointhomoveretaSmX}
Let $B$ be one-dimensional scheme, 
$z\in B^{(1)}$, $\eta\in B^{(0)}$.
% with a closed point $z\in B$ and a generic point $\eta\in B$.
Let $X\in \Sm_B$, $x\in X_{\underline{z}}=X\times_B z$.

Define $U=\Xx$ and $\calV=U\times_B \eta$.
% Then for any closed subscheme $Z$ in $X_{\underline{\eta}}=X\times_B \eta$, there are 
% Denote $\calV=(\Xx)_{\underline{\eta}}=\Xx\times_B \eta$.
Then for any closed subscheme $Z$ in $X_{\underline{\eta}}=X\times_B \eta$, there are framed correspondences 
$c\in \ZF_N(\calV\times\A^1, X_{\underline{\eta}})$, $c^\prime\in \ZF_N(\calV\times\A^1, X_{\underline{\eta}}-Z)$
such that \begin{itemize}
\item[(1)] $c \circ i_0=\sigma^N\can_{\underline{\eta}}$, $c\circ i_0 = j\circ c^\prime$, %\sigma^N_{\calV}
where $i_0,i_1\colon \calV\to \calV\times\A^1$ are the zero and unit sections, $\can_{\underline{\eta}}\colon \calV\to X_{\underline{\eta}}$, $j\colon X_{\underline{\eta}}-Z\to X_{\underline{\eta}}$ are the canonical morphisms,
\item[(2)] $c^{-1}(Z)$ is finite over $\calV$, see \eqref{sect:corrpreimZ} for $c^{-1}(Z)$. 
\end{itemize}
\end{proposition}
\begin{proof}
% Using the base change to the local scheme at the point $z$ of the irreducible component of $B$ containing $\eta$,
Using the base change along
the morphism $(\overline{\eta})_z\to B$,
we 
reduce the claim to %the one over 
local irreducible base scheme $B$.
By \Cref{prop:DimDim1}, 
there is a one-dimensional $\Xx$-focused compactified framed correspondence 
$(\ovC,C_\infty,\Psi)$
from the $\eta$-scheme $\calV$ to the $\eta$-scheme $X_{\underline{\eta}}$,
where $\Psi=(C_{\underline{\eta}}, \varphi, \val)$, $C_{\underline{\eta}}= (\ovC-C_\infty)\times_B \eta$,
such that \begin{itemize}
\item[(a1)]
$(\ovC,C_\infty,\Psi)$ smoothly contains $\can_{\underline{\eta}}$,
\item[(a2)]
the closure $\ovcalZ$ of $\calZ=\val^{-1}(Z)$ in $\ovC$ is finite over $\Xx$,
and $(\ovcalZ\cap C_\infty)\times_B\eta=\emptyset$,
% where $\Psi=(C_{\underline{\eta}}, \varphi, \val)$, $C_{\underline{\eta}}= (\ovC-C_\infty)\times_B \eta$.
\end{itemize}
% The claim follows by the combination of 
% \Cref{prop:DimDim1} and \Cref{lm:CurveFrCor}.
% \end{proof}
% \begin{lemma}\label{lm:CurveFrCor}
% % Let $B$ be a one-dimensional local irreducible scheme, $z\in B$ be the closed point, $\eta\in B$ be the generic point, $X\in Sm_B$, $Z$ be a closed subscheme in $X_{\underline{\eta}}=X\times_B \eta$.
% %
% %Suppose 
% %
% def:ZFovZF
We are going to construct 
framed correspondences 
$\tilde c\in \Fr_N(\calV\times\A^1, X_{\underline{\eta}})$, $c_0^+,c_1\in \Fr_N(\calV\times\A^1, X_{\underline{\eta}}-Z)$, $N\in \mathbb Z$, such that \begin{itemize}
\item[(c1)] $[\tilde c \circ i_0]=[\sigma^N\can_{\underline{\eta}}+ c^+_0]$, $[\tilde c \circ i_1]=[j\circ c_1]$ in $\ovZF_N(\calV,X_{\underline{\eta}})$
see
\Cref{def:ZFovZF},
% where $i_0,i_1\colon \calV\to \calV\times\A^1$ are the zero and unit sections, $j\colon X_{\underline{\eta}}-Z\to X_{\underline{\eta}}$ is the open immersion,
\item[(c2)] ${\tilde c}^{-1}(Z)$ is finite over $\calV$.
\end{itemize}
Then $c=\tilde c-(c_0^+\circ\pr)$, where $\pr\colon \calV\times\A^1\to\calV$ is the canonical projection, satisfies properties (1)-(2) in the lemma above.
% \end{lemma}
% \begin{proof}

Step 1)
% \item[Step 1)]
Firstly, assume the residue field $\calO(x)$ has at least three elements. %, $\#\mathcal O_z(z)>2$.
According to \Cref{def:reldimFr,def:compactifiedfocusedFr,def:smoothcontain} 
we are given with 
\begin{itemize}
\item[(d1)]
the projective equidimensional scheme $\ovC$ over $\Xx$ of relative dimension one, 
a closed subscheme $C_{\infty}$ finite over $\Xx$, 
a closed embedding $C_{\underline{\eta}}\to \A^N_{\calV}$ for some $N\in \mathbb Z$;
\item[(d2)]
functions $e_i\in \calO(\mathcal E)$, $i=2,\dots, N$, and a regular map $\val\colon \mathcal E\to X_{\underline{\eta}}$,
such that $C_{\underline{\eta}} = Z(e_2,\dots,e_N)$,
where $\mathcal E$ is the henselisation $(\A^N_{\calV})^h_{C_{\underline{\eta}}}$;
\item[(d3)]
a map $\Delta\colon \Xx\to C$, 
such that $\val\circ\Delta=\can\colon \Xx\to X$ is the canonical map,
and the morphism $\ovC\to \Xx$ is smooth over $\Delta(\Xx)$.
\end{itemize}

Here 
$\ovC_{\underline{\eta}}=\ovC\times_B \eta$, %=\ovC\times_{\Xx} \calV
and denote
$C_{\infty,\underline{\eta}}= C_{\infty}\times_{B} \eta$, %C_{\infty}\times_{\Xx} \calV
% \[\ovC_{\underline{\eta}}=\ovC\times_B \eta=\ovC\times_{\Xx} \calV,\quad C_{\infty,\underline{\eta}}= C_{\infty}\times_{\Xx} \calV,\]
and note that $\calZ=\val^{-1}(Z)=\ovcalZ\times_B\eta$. 
Then 
$\ovC_{\underline{\eta}}$ is a projective scheme of relative dimension one over $\calV$,
and $C_{\infty,\underline{\eta}}$ and $\calZ$ are closed subschemes finite over $\calV$.

We let us write $\Delta\subset\ovC$ for the subscheme $\Delta(\Xx)$.
By (d3) the morphism $\ovC\to \Xx$ is smooth over $\Delta$.
It follows that there is a line bundle $\calL(\Delta)$ on $\ovC$ with a section $\delta$ 
such that $Z(\delta)=\Delta$.
% that vanishing locus $Z(\delta)$ equals $\Delta(\Xx)$.
We write $\calL(\Delta)$ for 
%both 
% $\calL(\Delta(\Xx))$ and 
$\calL(\Delta)\big|_{\ovC_{\underline{\eta}}}$.
% its restriction on $\ovC_{\underline{\eta}}$.
%
% By (d3) the morphism $\ovC\to \Xx$ is smooth over $\Delta(\Xx)$.
% It follows that there is a line bundle $\calL(\Delta(\Xx))$ on $\ovC$ with a section $\delta$ 
% such that $Z(\delta)=\Delta(\Xx)$.
% % that vanishing locus $Z(\delta)$ equals $\Delta(\Xx)$.
% We write $\calL(\Delta)$ for 
% %both 
% $\calL(\Delta(\Xx))$ and 
% $\calL(\Delta(\Xx))\big|_{\ovC_{\underline{\eta}}}$.
% % its restriction on $\ovC_{\underline{\eta}}$.
%
% Since $\calZ$ and $\Delta(\Xx)$ are finite over the local scheme $\Xx$, 
% by \Cref{lm:semilocfinlocLinearbundle}, the restrictions
% $\calO(1)\big|_{\calZ}$ and $\calL(\Delta)\big|_{\calZ}$ are trivial.
% Since $\#\mathcal O_z(z)>2$, 
% by \Cref{lm:b0b1binfty},
% it follows that for any $d\in \mathbb Z_{\geq 0}$, there are invertible sections
% $\gamma_\infty, \gamma^+_0,\gamma_1\in \Gamma(\calZ, \mathcal O(d)^\times)$,  
% $\gamma_\infty= \gamma_1-\gamma_0\delta\big|_\calZ$.
% % \[\gamma_\infty, \gamma^+_0,\gamma_1\in \Gamma(\calZ, \mathcal O(1)^\times), \quad 
% % \gamma_\infty= \gamma_1-\gamma_0\delta\big|_\calZ.\]
% Since $C_{\infty,\underline{\eta}}$ is finite over $\Xx$ there is an invertible section 
% % \[
% $
% \beta_0^+\in \Gamma(C_{\infty,\underline{\eta}},\mathcal O(1)^\times)
% $.
% % \]
Since $\calZ$, $\Delta$ and $C_{\infty,\underline{\eta}}$ are finite over $\Xx$, 
by \Cref{lm:semilocfinlocLinearbundle}, 
the restrictions
$\calO(1)\big|_{\calZ\cup\Delta}$, 
$\calL(\Delta)\big|_{\calZ\cup\Delta}$, 
$\mathcal O(1)\big|_{C_{\infty,\underline{\eta}}}$
are trivial.
Since $\#\mathcal O_z(z)>2$, 
by \Cref{lm:b0b1binfty}, %it follows that 
for any $l\in \mathbb Z_{\geq 0}$, there are invertible sections
$\gamma_\infty, \gamma^+_0,\gamma_1\in \Gamma(\calZ\cup\Delta, \mathcal O(l)^\times)$,  
$\gamma_\infty= \gamma_1-\gamma^+_0\delta\big|_{\calZ\cup\Delta}$,
% \[\gamma_\infty, \gamma^+_0,\gamma_1\in \Gamma(\calZ, \mathcal O(1)^\times), \quad 
% \gamma_\infty= \gamma_1-\gamma_0\delta\big|_\calZ.\]
and
let
% Let
% \[
$
\beta_0^+\in \Gamma(C_{\infty,\underline{\eta}},\mathcal O(1))
$
% \]
%^\times
be
an invertible section.
% Also, there 
% is an invertible section 
% $
% \beta_0^+\in \Gamma(C_{\infty,\underline{\eta}},\mathcal O(1)^\times)
% $

Since $\calZ\cap C_{\infty,\underline{\eta}}= \emptyset$ 
by (a2),
% by the assumption (2) of the lemma, 
by Serre's theorem %on ample bundles 
\cite[III, Corollary 10.7]{Hartshorne-AlG}, 
for a large enough $l\in\mathbb Z$, 
there are sections 
\[\begin{array}{ll}
s_0^+\in \Gamma(\ovC_{\underline{\eta}}, \calL(\Delta(\calV))(l))),&
s_1\in \Gamma(\ovC_{\underline{\eta}}, \calO(l))\\ 
s_\infty\big|_{\calZ}=\gamma_\infty, &
s_0^+\big|_{\calZ\cup\Delta}= \gamma^+_0, \\
s_\infty\big|_{C_{\infty,\underline{\eta}}}=0, &
s_0^+\big|_{C_{\infty,\underline{\eta}}}=(\beta_0^+)^d
.\end{array}\]
% there are sections 
% $s_0^+\in \Gamma(\ovC_{\underline{\eta}}, \calL(\Delta(\calV))(d)))$ and $s_1\in \Gamma(\ovC_{\underline{\eta}}, \calO(d))$ 
% such that
% \[\begin{array}{ll}
% s_\infty\big|_{\calZ}=\gamma_\infty, &
% s_0^+\big|_{\calZ}= \gamma^+_0, \\
% s_\infty\big|_{C_{\infty,\underline{\eta}}}=0, &
% s_0^+\big|_{C_{\infty,\underline{\eta}}}=(\beta_0^+)^d
% .\end{array}\]
Define 
% \begin{gather*}
% s_0=\delta s_0^+, \quad s_1=s_0+s_\infty,\\
% s_\lambda = s_0(1-\lambda)+s_1\lambda.
% \end{gather*}
\[
s_0=\delta s_0^+, \quad s_1=s_0+s_\infty, \quad s_\lambda = s_0(1-\lambda)+s_1\lambda.
\]
Let
$g\colon (\A^N_{\calV\times\A^1})^h_{Z(s_\lambda)}\to X_{\underline{\eta}}$ be 
a lifting of the morphism 
%display%
\[
% $
Z(s_\lambda)\to C_{\underline{\eta}}\times\A^1\to C_{\underline{\eta}}\to X_{\underline{\eta}}
% $
\]
along the canonical closed immersion, see \cite[Theorem I.8]{Gru}, see also %or 
\cite{Elkiksoleqhens}, \cite[Lemma 3.11]{FrRigidSmAffpairs},
and
$e_1\in \calO(\A^N_{\calV\times\A^1})$ be a lifting of 
$s_\lambda/t_\infty^d\in \calO(C_{\underline{\eta}}\times\A^1)$.
% Since 
% $s_0^+\big|_{\calZ}=\delta^{-1}\gamma_0$ and $s_1\big|_{\calZ}=\gamma_1$ are invertible,
% $g$ induces maps $\val_0\colon (\A^N_\calV)^h_{Z(s_0^+)}\to X_{\underline{\eta}}-Z$ and $\val_1\colon (\A^N_{\calV})^h_{Z(s_1)}\to X_{\underline{\eta}}-Z$.
Define
\[\begin{array}{lclcl}
\tilde c&=&(Z(s_\lambda), e_1, e_2,\dots, e_N, g)&\in& \Fr_N(\calV\times\A^1, X_{\underline{\eta}}),\\ 
c_0^+ &=& (Z(s_0^+), e_1\big|_{\A^N_{\calV\times 0}}, e_2,\dots, e_N, \val_{0})&\in& \Fr_N(\calV,X_{\underline{\eta}}-Z),\\
c_1 &=& (Z(s_1), s_1/t_\infty^d, e_2,\dots, e_N, \val_{1})&\in& \Fr_N(\calV,X_{\underline{\eta}}-Z).
\end{array}\]
where
$\val_i\colon (\A^N_\calV)^h_{Z(s_i^+)}\to X_{\underline{\eta}}-Z$, $i=0,1$, 
% and $\val_1\colon (\A^N_{\calV})^h_{Z(s_1)}\to X_{\underline{\eta}}-Z$
are
induced by
$g$, 
because
sections $s_0^+\big|_{\calZ}=\gamma^+_0$ and $s_1\big|_{\calZ}=\gamma_1$ are invertible.
% $s_0^+\big|_{\calZ}=\delta^{-1}\gamma_0$ and $s_1\big|_{\calZ}=\gamma_1$ are invertible.
% and the maps $g_0\colon (\A^N_\calV)^h_{Z(s_0^+)}\to X_{\underline{\eta}}-Z$ and $\val_1\colon (\A^N_{\calV})^h_{Z(s_1)}\to X_{\underline{\eta}}-Z$ 
% are induced by $g$, because $s_0^+\big|_{\calZ}=\delta^{-1}\gamma_0$ and $s_1\big|_{\calZ}=\gamma_1$ are invertible.
%
Then since the section $s^+_0\big|_{\Delta}=\gamma_0\big|_{\Delta}$ is invertible,
$\tilde c\circ i_0= \can^\prime_{\underline{\eta}}+ j\circ c_0^+$,
% \[\tilde c\circ i_0= \can^\prime_{\underline{\eta}}+ j\circ c_0^+, \quad \tilde c\circ i_1 = j\circ c_1,\]
where $\can^\prime_{\underline{\eta}} = (\Delta(\calV), e_1, e_2,\dots, e_N, g)=\sigma^D\can_{\underline{\eta}}$, %it was \sigma^N_D 
and $D\in \GL_N(\calV)$ is defined by the differentials of functions $e_i$.
We redenote 
$\tilde c:=\tilde c\circ \sigma^{D^{-1}}_{\calV\times\A^1}$, and similarly for $c_0$, $c_0^+$, $c_1$, 
then by \Cref{lm:ElemenatryMatrixesZaroSectionSupportframedhomotopy}, since $\Xx$ is local,
we conclude
$[\tilde c\circ i_0]= [\can_{\underline{\eta}}+ j\circ c_0^+]$,
$\tilde c\circ i_1 = j\circ c_1$.

Since the section
%display% \[
$
s_\lambda\big|_{C_{\infty,\underline{\eta}}\times\PP^1}= \delta\big|_{C_{\infty,\underline{\eta}}}(\beta_0^+)^d
$
% \]
is invertible,
$Z(s_\lambda)\subset C_{\underline{\eta}}\times\PP^1$,
and the morphism
$p\colon Z(s_\lambda)\to \calV\times\A^1$ is finite. %A
Hence the scheme 
\[c^{-1}(Z) = p(g^{-1}_{Z(s_\lambda)}(Z))=p(Z(s_\lambda\big|_{\calZ\times\A^1}) )\] is finite over $\calV$
since the section $s_\infty\big|_{\calZ}=\gamma_\infty$ is invertible.

% Since 
% %\[%todo%display
% $s_\lambda\big|_{\calZ\times\infty}=s_\infty\big|_{\calZ}=\gamma_\infty
% $
% % \]
% is invertible,
% $Z(s_\lambda\big|_{\calZ\times\PP^1})$ is contained in $\calZ\times\A^1$ and is finite over $\calZ$, and consequently, it is finite over $(\Xx)_{\underline{\eta}}=\calV$.
% Hence  
% \[c^{-1}(Z) = p(g^{-1}_{Z(s_\lambda)}(Z))=p(Z(s_\lambda\big|_{\calZ\times\A^1}) )\] is finite over $\calV$,
% where $p\colon Z(s_\lambda)\to \calV\times\A^1$ is the finite morphism induced by the projection $C_{\underline{\eta}}\times\A^1\to \calV\times\A^1$.
%
% Then
% \[s_\lambda\big|_{C_{\infty,\underline{\eta}}\times\PP^1}= \delta\big|_{C_{\infty,\underline{\eta}}}(\beta_0^+)^d\]
% is an invertible section.
% Hence $Z(s_\lambda)$ is finite over $\Xx\times\PP^1$ %A
% and is contained in $C_{\underline{\eta}}\times\PP^1$.
% Next, \[s_\lambda\big|_{\calZ\times\infty}=s_\infty\big|_{\calZ}=\gamma_\infty\]
% is invertible, and it follows that
% $Z(s_\lambda\big|_{\calZ\times\PP^1})$ is contained in $\calZ\times\A^1$ and is finite over $\calZ$, and consequently, it is finite over $(\Xx)_{\underline{\eta}}=\calV$.
% Hence  
% \[c^{-1}(Z) = p(g^{-1}_{Z(s_\lambda)}(Z))=p(Z(s_\lambda\big|_{\calZ\times\A^1}) )\] is finite over $\calV$,
% where $p\colon Z(s_\lambda)\to \calV\times\A^1$ is a finite morphism induced by the projection $C_{\underline{\eta}}\times\A^1\to \calV\times\A^1$.

% \item[Step 2)]
Step 2)
%Consider the case of an arbitrary residue field $k=\mathcal O_z(z)$.
For an arbitrary residue field $k=\mathcal O_z(z)$,
let $q_0,q_1\in \mathbb Z_{>0}$, 
%disp\[(q_0,q_1)=1, \quad (q_0,\chark k)=1, \quad (q_2,\chark k)=1.\]
$(q_0,q_1)=1$, $(q_0,\chark k)=1$, $(q_2,\chark k)=1$.
Define polynomials 
\[f_i= t^{q_i-1}+t^{q_i-2}+\dots +1=(t^{q_i}-1)/(t-1)\in R[t], R=\calO_B(B).\]
Then $f_0$ and $f_1$ are separable over the closed point $z$.
% Denote by $B_i=Z(f_i)$ the closed subschemes in $\A^1_B$, and let $p_i\colon B_i\to B$ denote the canonical projections.
% Then $p_i$ are finite \'etale morphisms of schemes.
Denote by $B_i=Z(f_i)$, $i=0,1$,
% the closed subschemes in $\A^1_B$, 
then 
the projection to $\A^1_B\to B$
induces finite \'etale morphisms $p_i\colon B_i\to B$,
and 
there are \'etale neighbourhoods $V_i$ of 
$B_i$ in $\A^1_B$
% $Z(f_i)$ 
with retractions
$r_i\colon V_i\to B_i$.
% Denote
% $p_i^{X_{\underline{\eta}}}=
% % p_i\times_B id_{X_{\underline{\eta}}}
% \colon B_i\times_B X_{\underline{\eta}}\to X_{\underline{\eta}}$.
Define
%display% \[
% $
% e_i=(Z(f_i)\times_B\calV,f_i,r_i)\in \Fr_1(\calV, B_i\times_B\calV)
% $
$
e_i=(Z(f_i),V_i,f_i,r_i)\in \Fr_1(B, B_i)
$
% \]
and
$c = \sum\limits_{i=0,1}
(-1)^i
(p_i\times_B id_{X_{\underline{\eta}}})\circ c^{B_i}\circ (e_i\times_B id_{\calV\times\A^1})\circ (b_i)_\varepsilon
% + 
% c^{B_1}\circ (e_1\times id_{\A^1})\circ (b_1)_\varepsilon
\in 
\ZF_{2+N}(\calV\times\A^1, X_{\underline{\eta}})$,
where $c^{B_i}\in\ZF_{N}(\calV\times_B B_i\times\A^1,X_{\underline{\eta}})$ 
from
% denote
% the correpondences constructed in 
(Step 1).
Then
\[\begin{array}{lcl}
c\circ i_0 &=& 
\sum\limits_{i=0,1}
(-1)^i
(p_i\times_B id_{X_{\underline{\eta}}})\circ
c^{B_i}_0\circ (e_i\times_B id_{\calV})\circ (b_i)_\varepsilon
\\&\stackrel{\text{Step 1}}{=}&
\sum\limits_{i=0,1}
(-1)^i
(p_i\times_B id_{X_{\underline{\eta}}})\circ\sigma^N\can^{B_i}\circ (e_i\times_B id_{\calV})\circ (b_i)_\varepsilon\\ &=& 
% \sum\limits_{i=0,1}
% (-1)^i
% \sigma^N \can\circ ((p_i\times_B id_{\calV})\circ (e_i\times_B id_{\calV})\circ (b_i)_\varepsilon.
\sigma^N \can\circ (\sum\limits_{i=0,1}(-1)^i(p_i\circ e_i\circ (b_i)_\varepsilon)\times_B id_{\calV})
\\&\stackrel{\text{\Cref{lm:nvarepsilonclass}}}{=}&
\sigma^N \can\circ (\sum\limits_{i=0,1}(-1)^i((q_i)_{\varepsilon}\circ (b_i)_\varepsilon)\times_B id_{\calV})
\\&\stackrel{\text{\Cref{lm:classtrivsigma}}}{=}&
\sigma^N \can\circ (\sigma^{2}\times_B id_{\calV})
\\&=&
\sigma^{2+N} \can\circ (id_{\calV})
% \\&
\in
% &
\ovZF_{2+N}(\calV, X_{\underline{\eta}}).
\end{array}\tag*{\qedhere}\] 
%% where $p_i^{\calV}\colon B_i\times_B\calV\to \calV$.
% Finally, since by \Cref{lm:nvarepsilonclass} %,%punctuation 
% $[p_i\circ e_i]=[(q_i)_{\varepsilon}]\in \ovZF_1(\calV, \calV)$,
% it follows by \Cref{lm:classtrivsigma} that
% % \[
% $
% [c\circ i_0] = [\sigma^{2+N}\can]\in\ovZF_{2+N}(\calV, X_{\underline{\eta}})
% $.
% \]
% and
% $c = c^{B_0}\circ e^{\A^1}_0 + c^{B_1}\circ e^{\A^1}_1\in \ZF_{2+N}(\calV\times\A^1, X_{\underline{\eta}})$.
% Then
% \begin{multline*}
% c\circ i_0 = 
% c^{B_0}_0\circ e_0\circ (b_0)_\varepsilon - c^{B_1}_0\circ e_1\circ (b_1)_\varepsilon = \\
% \sigma^N\can^{B_0}\circ e_0\circ (b_0)_\varepsilon - \sigma^N\can^{B_1}\circ e_1\circ (b_1)_\varepsilon = \\
% \sigma^N \can\circ (p_0\circ e_0\circ (b_0)_\varepsilon- p_1\circ e_1\circ (b_1)_\varepsilon),
% \end{multline*} 
% where $p_i^{\calV}\colon B_i\times_B\calV\to \calV$.
% Finally, since by \Cref{lm:nvarepsilonclass} %,%punctuation 
% $[p_i\circ e_i]=[(q_i)_{\varepsilon}]\in \ovZF_1(\calV, \calV)$,
% it follows by \Cref{lm:classtrivsigma} that
% % \[
% $
% [c\circ i_0] = [\sigma^{2+N}\can]\in\ovZF_{2+N}(\calV, X_{\underline{\eta}})
% $.
% % \]
\end{proof}

\subsection{Vanishing of cohomologies.} % at infinity
\begin{theorem}\label{th:GenFifLocEssSmSchemeonedimbasetriviality}
Let $X$ be a smooth scheme over a one-dimensional irreducible base scheme $B$. Let $x\in X$ be a point over a closed point $z\in B$.
Let $F$ be an $\A^1$-invariant quasi-stable framed linear presheaf over the generic point $\eta\in B$. 
Then 
$
% \[
\underline{H}^l_\nis(\Xx\times_B \eta, F)\cong 0
$.
% \]
\end{theorem}
\begin{proof}
The claim follows by the combination of \Cref{cor:preimFinFrCintr->Inj} with 
\Cref{cor:contractinggenpointhomoveretaSmX}.
\end{proof}

\section{Strict homotopy invariance}\label{sect:StrHomInvTheorem}

% In \Cref{th:GenFifLocEssSmSchemeonedimbasetriviality}
% % proves the va
% % In this section,
% we prove 
% triviality of ``cohomologies at infinity'' 
% on %the relative affine line 
% $\A^1_V$ over 
% a local essentially smooth scheme $V$ over a field $k$,
% \Cref{th:GenFifLocEssSmSchemeonedimbasetriviality}
% applying criterion form \Cref{sect:InjThbyFinSupHom}
% % \Cref{cor:preimFinFrCintr->Inj}
% and 
% the $\A^1$-homotopy 
% % result of 
% from
% \Cref{sect:homotpiesatinfinity}.
% % to a framed $\A^1$-homotopy ``with finite supports''
% % constructed in \Cref{cor:contractinggenpointhomoveretaSmX}.
% Then we deduce the strict homotopy invariance theorem, see \Cref{th:strhominv}, combining 
% \Cref{th:GenFifLocEssSmSchemeonedimbasetriviality}
% % the latter result 
% with 
% \Cref{cor:NisFalongA1hominv}.
% % from \Cref{sect:outofinffinitesupSHI}.

We prove the strict homotopy invariance theorem 
combining 
% the triviality of ``cohomologies at infinity'' on $\A^1\times X^h_x$ for $X\in\Sm_k$
% provided by
\Cref{th:GenFifLocEssSmSchemeonedimbasetriviality}
and 
% the triviality of ``cohomologies out infinity''
% provided by
% \Cref{cor:NisFalongA1hominv}.
% or 
\Cref{th:FinSuppLocHensbaseafflinetriviality}.
\begin{theorem}\label{th:strhominv}
Let $F$ be a quasi-stable homotopy invariant framed linear presheaf over $k$.

(1)
% The presheaves $h_\nis^i(F)= H^i_\nis(-,F)$, $i\geq 0$, 
% on the category $\Sm_k$ are $\A^1$-invariant.
The presheaves $h_\nis^i(F)\colon \Sm_k\to\Ab$, $X\mapsto H^i_\nis(X,F)$,
are $\A^1$-invariant for all $i\geq 0$.

(2)
% \label{cor:HniscongHzar}
There are canonical isomorphisms
$H^i_\nis(X,F_\nis)\cong H^i_\zar(X,F_\zar)$ for all $X\in \Sm_k$, $i\geq 0$.
% \label{cor:ZarStrictHomotopyInvariance}
% (3)
% The canonical projection induces the isomorphism
% $H^*_\zar(\A^1\times X,F_\zar)\simeq H^*_\zar(X,F_\zar)$
% for all $X\in \Sm_k$.
\end{theorem}
\begin{proof}
(1)
The claim is equivalent to
isomorphisms
\[\begin{cases} H^0_\nis(V\times\A^1, F)\cong F(V\times\A^1),  \\ H^i_\nis(V\times\A^1, F)=0, &i>0,\end{cases}\]
for
all local henselian essentially smooth schemes $V$
over $k$.
The morphism of sites given by the embedding of the topology $\Fin_{S}$ on $\Sm_{\A^1\times V}$ into the Nisnevich topology leads to the spectral sequence 
\begin{equation}\label{eq:H(Fin<-VNis)}
H^p_{\Fin_V} (\A^1\times V, H^q_{\nis}(F)_{\Fin_V} )\Rightarrow H^{r}_\nis(\A^1\times V, F_\nis),r=p+q,
\end{equation}
where ${H}^*_\nis(F)_{\Fin_V}$ denotes the sheafification with respct to $\Fin_V$ of the presheaf $H^*_\nis(-,F_\nis)$.
By \Cref{lm:NisZXNisZpoints}, the topology $\Fin_V$ has enough set of points given by schemes of the forms
% \[
$U^h_v$, 
for $v\in U$,
and
$U\times_{\PP^1\times V}(\PP^1\times V)_{(\infty,z)}$,
% \]
where $z$ is the closed point of $V$.
Since the schemes $U^h_v$ are Nisnevich points, $\underline H^*_\nis(U^h_v,F)\cong 0$.
By \Cref{th:GenFifLocEssSmSchemeonedimbasetriviality}, \[\underline H^*_\nis(U\times_{\PP^1\times V}(\PP^1\times V)_{(\infty,z)},F)\simeq 0.\]
Combining the latter two isomorphisms, we conclude
$\underline{H}^*_\nis(F)_{\Fin_V}\cong 0$, and consequently%
,
% \[
$
{H}^*_\nis(F)_{\Fin_V}\cong F_{\Fin_V}
$.
% \]
Thus \eqref{eq:H(Fin<-VNis)} implies \[H^p_{\Fin_V} (\A^1\times V, F_{\Fin_V} )\cong H^{p}_\nis(\A^1\times V, F_\nis);\]
and the claim, since $\underline{H}^*_{\Fin_V}(\A^1\times V, F)\cong 0$ by \Cref{th:FinSuppLocHensbaseafflinetriviality}.
% \end{proof}
% The following theorem reformulates the above result in an alternative standard form.
% \begin{theorem}\label{cor:strhominv}
% The presheaves $h_\nis^i(F)= H^i_\nis(-,F)$, $i\geq 0$, on the category $\Sm_k$ are $\A^1$-invariant.
% \end{theorem}
% \begin{proof}

(2)
We repeat the argument from \cite{Sus-nonperftalk}. 
% Since 
% the presheaves $H^*_\nis(-,F_\nis)$ are $\A^1$-invariant, 
By point (1) and \cite[Theorem 3.11]{hty-inv} the canonical morphism $H^*_\nis(U,F_\nis)\to H^*_\nis(\eta,F_\nis)$ is injective for any essentially smooth local scheme $U$ with the generic point $\eta$. Hence the presheaf $\underline{H}^*_\nis(-,F_\nis)$ is Zariski locally trivial. Then because of the spectral sequence $H^{p}_\zar(X, H^{q}_\nis(-,F_\nis)_\zar)\Rightarrow H^{p+q}_\nis(X,F_\nis)$ the claim follows.
\end{proof}

\section{Cousin complex}\label{sect:CousinExactness}

Let $B$ be a one-dimensional local irreducible noetherian scheme.
Let $E\in\SH(B)$.
Then by Voevodsky's Lemma,
see \cite[Construction 3.1]{Framed},
% by Voevodsky's lemma 
$E$ defines the functor 
$E\colon \Fr_+(B)\to \SH$
such that for any $X,X_0,X_1\in \Sm_B$, the morphisms 
%\[F(X)\xrightarrow[\simeq]{\sigma^*_X}F(X)), F(X_0\amalg X_1)\simeq F(X_0)\oplus F(X_1), F(X\times\A^1)\simeq F(X).\]
\[
F(X)\xrightarrow{\sigma^*_X}F(X), \quad 
F(X_0\amalg X_1)\to F(X_0)\oplus F(X_1), \quad 
F(X)\to F(X\times\A^1)
\]
are isomorphisms in $\SH$.
% the category $\mathcal Spt$ of $S^1$-spectra.
Consequently, since the stable homotopy category $\SH$ is additive, $F$ defines the functor
\begin{equation}\label{eq:ZF->SH}F\colon \ZF_*(B)\to \SH.\end{equation}

\begin{definition}
For $X\in \Sm_B$, and a closed subscheme $Z$ in $X$,
denote \[E^l_Z(X)=\pi_{-l}\coCone(E(X)\to E(X-Z)), l\in \mathbb Z,\]
% Denote by $\ZF_*(B)^\mathrm{pair}$ the category of morphisms in $\ZF_*(B)$ spanned by open immersions $X-Z\hookrightarrow X$ in $\Sm_B$, then there are 
% that 
% leads to 
%has the stucture of %defines %induces
and define 
functors \begin{equation}\label{eq:E(X/X-Z)}\ZF_*(B)^\mathrm{pair}\to \Ab; \quad ((X-Z)\hookrightarrow X)\to E^l_Z(X).\end{equation}
\end{definition}
\begin{lemma}\label{lm:A1invfibreparis}
Presheaves \eqref{eq:E(X/X-Z)} are $\A^1$-invariant.
\end{lemma}
\begin{proof}
The claim follows because \eqref{eq:ZF->SH} is $\A^1$-invariant.
% By \Cref{th:strhominv} the Nisnevich sheaves associated with the abelian presheaves $\pi_{-l} L_\nis F\colon \ZF_+(B)\to \Ab$, where $l\in \mathbb Z$, are strictly $\A^1$-homotopy invariant. 
% Hence $L_\nis F$ defines an $\A^1$-invariant functor $\ZF_*(B)\to \SH$.
% Hence the presheaves \eqref{eq:E(X/X-Z)} are $\A^1$-invariant.
\end{proof}

%Let $F$ be an $\A^1$-invariant additive presheaf $\mathrm{Corr}^\fr(k)\to \SH$.
%For $X\in \Sm_k$, and a closed subscheme $Z$ in $X$,
%denote $E^l_Z(X)=\pi_{-l}\operatorname{fib}(L_\nis F(X)\to L_\nis F(X-Z))$,
%where $L_\nis$ denote the associated Nisnevich local sheaf.
%Denote by $\mathrm{Corr}^\fr(k)^\mathrm{pair}$ the subcategory in the category of arrows of $\mathrm{Corr}(k)$ spanned by the morphisms %induced by open immersions $X-Z\to X$ in $\Sm_k$, then there is the functor %\begin{equation}\label{eq:E(X/X-Z)}\mathrm{Corr}^\fr(k)^\mathrm{pair}\to \SH; ((X-Z)\to X)\to E^l_Z(X).\end{equation}
%By \Cref{cor:strhominv} and \cite[Proposition 3.2.14]{five-authors}, 
%$L_\nis F$ defines $\A^1$-invariant functor $\mathrm{Corr}^\fr(k)\to \SH$.
%Hence the presheaves \eqref{eq:E(X/X-Z)} are $\A^1$-invariant.
%\begin{remark}
%Note that any $\A^1$-invariant stable framed Nisnevich sheaf of abelian groups $F\in \mathrm{HI}^\ZF_\nis(k)$ defines 
%the $\A^1$-invariant presheaf of Eilenberg-MacLane spectra $\mathrm{Corr}^\fr(k)\to \SH$ by %\Cref{}
%\cite[Theorem 5.14]{TowardsconservativityofGm–stabilization}. %Eilenber-Maklein spectra presheaf that 
%\end{remark}

\begin{lemma}\label{lm:extensionofsupperts}
Let $U$ be essentially smooth local scheme over $B$,
and $\calV=U\times_B \eta$, where $\eta$ is the generic point of $B$. 
Then for any closed subscheme $W$ of codimension $r$ in $\calV$, 
the morphism 
% \begin{equation}\label{eq:ElWcalVtovarinjcodVWpElWpcalV}
% E^l_W(\calV)\to \varinjlim_{\codimcalVWprimermo} E^l_{W^\prime}(\calV)
% \end{equation}
\begin{equation}\label{eq:ElYcalVtovarinjcodVWpElWpcalV}
E^l_Y(X\times_B\eta)\to \varinjlim_{\codimcalVWprimermo} E^l_{W^\prime}(\calV)
\end{equation}
is trivial for 
any $l\in\mathbb Z$,
where 
$X=U$, $Y=W$,
and
$W^\prime$ runs over the filtered set of closed subschemes of codimension $r-1$.
\end{lemma}
\begin{proof}
For any 
scheme $X\in \Sm_B$, 
point $x\in X$ and 
% $Y\hookrightarrow X$
closed subscheme $Y$ of $X_{\underline{\eta}}=X\times_B\eta$
such that 
$U=X_x$, 
$W=Y\times_X U$
% $W=Y_x$,
% by \cite[Theorem 2.11]{hty-inv}, if $\dim B=0$ $Z(f)=\emptyset$, and 
by \Cref{cor:contractinggenpointhomoveretaSmX},
% , if $Z(f)\neq\emptyset$,
there are 
\[c^\prime\in \ZF_N(\calV,X_{\underline{\eta}}-Y), \quad c=(Z,V,\varphi,g)\in \ZF(\calV\times\A^1, X_{\underline{\eta}})\] such that 
$c\circ i_0=\sigma^N_{\calV}\can$, $c\circ i_1=j\circ c^\prime$, 
see \Cref{sect:not:FSSubschemes},
where 
% $i_0,i_1\colon \calV\to \calV\times\A^1$ are the zero and unit sections, 
$\can\colon \calV\to X_{\underline{\eta}}$, $j\colon X_{\underline{\eta}}-Y\to X_{\underline{\eta}}$ are the canonical morphisms. 
% $\can\colon \calV\to X_{\underline{\eta}}$ is the canonical morphism, $j\colon X_{\underline{\eta}}-Y\to X_{\underline{\eta}}$ is the open immersion. 
Denote by $\widetilde W$ 
% \hookrightarrow\calV
the image closure 
% of 
% $Z\times_{X} Y$ %$g$
% along 
% the composite morphism
% $Z\rightarrow V\times \A^1 \rightarrow \calV$.
along the composite morphism
$Y\times_{X}Z\hookrightarrow %g^{-1}(Y)%$g$
Z\rightarrow \calV$. %\calV\times \A^1 \rightarrow \calV$.
% $\calV\leftarrow V\times \A^1 \leftarrow Z\xrightarrow{g} X$.
% the preimage of $Y$ along the morphisms
% $\calV\leftarrow V\times \A^1 \leftarrow Z\xrightarrow{g} X$.
Then $\operatorname{codim}{\widetilde W}\geq r-1$.
Since for each $l\in\mathbb Z$, 
% $E^l_{(-)}(-)$ defines an $\A^1$-invariant functor $\ZF_*(k)^\mathrm{pair}\to \mathrm{Ab}$ 
by \Cref{lm:A1invfibreparis},
there is the equality of homomorphisms 
\[(j\circ c^\prime)^*=\can^*\colon E^l_Y(X_{\underline{\eta}})\to E^l_{\widetilde W}(\calV).\] 
% are equal.
So since $(j\circ c^\prime)^*=0$, it follows that $\can^*=0$,
and morphism 
\eqref{eq:ElYcalVtovarinjcodVWpElWpcalV} 
% \[j^*\colon E^l_Y(X)\to \varinjlim_{\operatorname{codim}_{\calV}W^\prime=r-1}E^l_{W^\prime}(\calV)\] 
is trivial
% .
% Since the latter claim holds 
for all $X$ and $Y$ as above.
Hence
% the morphism 
\eqref{eq:ElYcalVtovarinjcodVWpElWpcalV} is trivial
for 
$X=U$, $Y=W$.
% $Y=W$, $X=\calV$.
\end{proof}
% \begin{remark}
% Note that the case of empty scheme $Z(f)$ in the latter lemma follows by the same argument with the use of the injectivity on local essentially smooth schemes \cite[Theorem 3.11]{hty-inv}.
% \end{remark}
\begin{corollary}\label{cor:shortexactseuence}
Under the assumptions of \Cref{lm:extensionofsupperts},
% for any $l\in\mathbb Z$ and $W\in\calV^{[r]}$,
% closed subscheme $W$ of $\calV$ such that $\operatorname{codim}{W}=r$,
the %short 
sequence
\[0\to \varinjlim_{\begin{smallmatrix}\codimcalVWprimermo,\\W^\prime\supset W\end{smallmatrix}} E^l_{W^\prime}(\calV) \to \varinjlim_{\codimcalVWprimermo} E^l_{W^\prime}(\calV-W)\to E^{l+1}_W(\calV)\to 0\]
is exact,
for any $l\in\mathbb Z$ and $W\in\calV^{[r]}$. 
% for any $l\in\mathbb Z$ and
% closed subscheme $W$ of $\calV$ such that $\operatorname{codim}{W}=r$.
% . Recall that $W$ here stands for any closed subscheme of $V$ such that $\operatorname{codim}{W}=r$.
\end{corollary}
\begin{proof}
The claim follows
from \Cref{lm:extensionofsupperts}
because of
the long exact sequence 
\[\dots\to \varinjlim_{\codimcalVWprimermo} E^l_{W^\prime}(\calV) \to \varinjlim_{\codimcalVWprimermo} E^l_{W^\prime}(\calV-W)\to E^{l+1}_W(\calV)\to \dots.
\tag*{\qedhere}\]
\end{proof}
\begin{theorem}\label{th:CousinsubsetsLocEssSm}
% Let
% $F\in \SH^{S^1}(\mathrm{Corr}^\fr(B))$, 
% and denote
% $E^l_{Z}(X)=\pi_{-l}\operatorname{fib}(F(X)\to F(X-Z))$.
The Cousin complex
\eqref{eq:Cous(calV,E)}
% \[0\to E^l(\calV)\to E^l(\calV^{(0)})\to \bigoplus_{x\in V^{(1)}}E_x^{l+1}(\calV)\to \dots \bigoplus_{x\in V^{(r)}}E_x^{l+r}(\calV)\to \dots \bigoplus_{x\in V^{(d)}}E_x^{l+d}(\calV)\to 0\]
is acyclic for %the scheme 
$\calV=U\times_B \eta$ for any essentially smooth local scheme $U$ %of dimension $d$ 
over $B$. 
% where $\calV^{(c)}$ denotes the set of points of $\calV$ of codimension $c$. 
% in particular, the claim holds for $\SH^{S^1}(\mathrm{Corr}^\fr(B))$-representable or $\SH(B)$- cohomology theories.
\end{theorem}
\begin{proof}
% Note that $E_x^l(\calV)=E_x^l(\calV_x)$ by the definition.
The claim follows because \Cref{cor:shortexactseuence} provides short exact sequences
\[0\to 
\varinjlim_{\codimcalVWprimermo}E_W^{l+r-1}(\calV)\to \bigoplus_{y\in \calV^{(r-1)}}E_y^{l+r-1}(\calV)\to \bigoplus_{\codimcalVWprimer}E_{W^\prime}^{l+r}(\calV)\to 0.
\tag*{\qedhere}\]
% It follows from \Cref{cor:shortexactseuence} that the short exact sequences
% \[0\to 
% \varinjlim_{\operatorname{codim}{W}=r-1}E_W^{l+r-1}(\calV)\to \bigoplus_{x\in V^{(r-1)}}E_x^{l+r-1}(\calV_x)\to \bigoplus_{\operatorname{codim}{W}=r}E_{W^\prime}^{l+r}(\calV)\to 0.
% \tag*{\qedhere}\]
\end{proof}

% \begin{corollary}\label{cor:CousinsubsetsLocEssSm(k)}
% % Let
% % $F\in \SH^{S^1}(\mathrm{Corr})$, and $E^l_{Z}(X)=\pi_{-l}\operatorname{fib}(F(X)\to F(X-Z))$.
% For a base field $k$,
% the Cousin complex \eqref{eq:Cous(calV,E)}
% % \[0\to E^l(V)\to E^l(\eta)\to \bigoplus_{x\in V^{(1)}}E_x^{l+1}(V)\to \dots \bigoplus_{x\in V^{(r)}}E_x^{l+r}(V)\to \dots \bigoplus_{x\in V^{(d)}}E_x^{l+d}(V)\to 0\]
% is exact for $\calV=U-Z(f)$ for an essentially smooth local scheme $U$ over $k$
% and $f\in \mathcal O(U)$ such that $Z(f)$ is smooth.
% \end{corollary}
% \begin{proof}
%     The claim follows by \Cref{th:CousinsubsetsLocEssSm}
%     applied to $B=(\A^1_k)_{\{0\}}$ and 
%     the essentially smooth $B$-scheme given by $f\colon U\to(\A^1_k)_0$.
% \end{proof}
% \begin{corollary}\label{cor:CousinsubsetsLocEssSm(k)}
% Let
% $F\in \SH^{S^1}(\mathrm{Corr})$, and $E^l_{Z}(X)=\pi_{-l}\operatorname{fib}(F(X)\to F(X-Z))$.
\begin{example}
For a base field $k$,
the Cousin complex \eqref{eq:Cous(calV,E)}
% \[0\to E^l(V)\to E^l(\eta)\to \bigoplus_{x\in V^{(1)}}E_x^{l+1}(V)\to \dots \bigoplus_{x\in V^{(r)}}E_x^{l+r}(V)\to \dots \bigoplus_{x\in V^{(d)}}E_x^{l+d}(V)\to 0\]
is exact for $\calV=U-Z(f)$ for an essentially smooth local scheme $U$ over $k$
and $f\in \mathcal O(U)$ such that $Z(f)$ is smooth.
\end{example}

\appendix

\section{Relative dimension over one-dimensional base.}\label{sect-app:RelDim}

% In this appendix, we prove some results on the relative dimension over one-dimensional local base schemes.

\begin{proposition}\label{prop:PosCodinXdimgeq1forallx}
% Let $X$ be an affine scheme such that $\dim X^\prime\geq 1$ 
% for each irreducible component $X^\prime=\overline{\eta}$, $\eta\in X^{(0)}$.
Let $X$ be an affine scheme such that $\dim X^\prime\geq 1$ 
for each irreducible component $X^\prime$ of $X$. 
Then there is a function $f\in \calO_X(X)$ such that $\codim_{X^\prime}(Z(f)\cap X^\prime)=1$ for each irreducible component $X^\prime$.
\end{proposition}
\begin{proof}
Since for each irreducible component $\calO_{X^\prime}(X^\prime)$ is not a field, it follows that there is a closed subscheme $F$ of $X$ that intersects non-emptily and does not contain each irreducible component $X^\prime$.
% Let $X_{2}=\bigcup_{\eta_0,\eta_1\in X^{(0)},\eta_0\neq\eta_1}(\overline{\eta_0}\cap\overline{\eta_1})$.
Let $X_{2}$ be the closed subscheme in $X$ that is the union of pairwise intersections of the irreducible components of $X$.
Then $X_{2}$ and $X_{2}\cup F$ do not contain any irreducible component $X^\prime$.
Hence for any $X^\prime$, there is a function $f_{X^\prime}\in \calO_{X^\prime}(X^\prime)$ 
such that $f_{X^\prime}\neq 0$, $f_{X^\prime}\big|_{(X_{2}\cup F)\cap X^\prime}=0$,
% Then
and combining the functions $f_{X^\prime}$ for all $X^\prime$ together,
% for all irreducible components $X^\prime$, 
we get the function $X$
such that $Z(f)\cap X^\prime\neq \emptyset$, and $f\big|_{X^\prime}\neq 0$
for any irreducible component $X^\prime$ of $X$.
Hence $Z(f)$ is of positive codimension in $X$ and intersect non-emptily each $X^\prime$.
Since 
%for any function $f$, 
$\codim_{X^\prime}(Z(f)\cap X^\prime)\leq 1$
% $\codim_X Z(f)\leq 1$
for any function $f$,
% the codimension of $Z(f)$ is less or equal to one,
$\codim_{X^\prime}(Z(f)\cap X^\prime)= 1$.
% the claim is proven.
\end{proof}

\begin{proposition}\label{prop:OnedimBaseEquidimSch}
Let $B$ be a one-dimensional irreducible local scheme, 
$z\in B^{(1)}$, $\eta\in B^{(0)}$.
% $z\in B$ be the closed point, $\eta$ be the generic point.
Given an irreducible scheme $X$ over $B$ such that the schemes $X_{\underline{z}}=X\times_B z$ and $X_{\underline{\eta}}=X\times_B \eta$ are non-empty,
then $X$ is equidimensional over $B$, i.e. $\dim X_{\underline{z}}=\dim X_{\underline{\eta}}$.
\end{proposition}
\begin{proof}
Since any scheme has a Zariski covering by affine schemes, we can assume that $X$ is affine, 
and consequently, $X_{\underline{\eta}}$ is of such type.
If $\dim X_{\underline{z}}=0$, then $X$ is quasi-finite over $B$, and $\dim X_{\underline{\eta}}=0$ as well.
Suppose $\dim X_{\underline{z}}=d\geq 1$, and the claim is proven for all $X$ such that $\dim X_{\underline{z}}<d$.

By \Cref{prop:PosCodinXdimgeq1forallx}, there is a regular function $f_{\underline{z}}\in \calO_{X_{\underline{z}}}(X_{\underline{z}})$ such that $\codim_{X_{\underline{z}}} Z(f_{\underline{z}}) = 1$.
Since $X$ is affine, there is a lifting $f\in \calO_X(X)$ of $f_{\underline{z}}$.
Then since $X$ is irreducible, it follows that $\codim_X Z(f)= 1$.
% , and so $\codim_X Z(f)=1$.
Moreover, since $X$ is irreducible, it follows that $X_{\underline{\eta}}$ is dense in $X$, and $X_{\underline{\eta}}$ is irreducible too.
Hence $f_{\underline{\eta}}:=f\big|_{X_{\underline{\eta}}}\neq 0$, 
% So $\codim_{X_{\underline{\eta}}} Z(f_{\underline{\eta}}) \geq 1$, 
and it follows that 
$\codim_{X_{\underline{\eta}}} Z(f_{\underline{\eta}})=1$.

Let $X_1$ denote the irreducible component of $Z(f)$. 
By the above $\codim_{X_{\underline{z}}} (X_1\times_B z)= 1$ and $\codim_{X_{\underline{\eta}}} (X_1\times_B \eta)= 1$.
Then $X_1$ is an irreducible affine scheme over $B$ such that 
$X_1\times_B z$ and $X_1\times_B \eta$ are non-empty. By the inductive assumption $\dim (X_1\times_B z) = \dim (X_1\times_B \eta)$.
Hence
$\dim X_{\underline{z}}=\dim X_{\underline{\eta}}$.
\end{proof}
\begin{corollary}\label{cor:OnedimBaseEquidimSch}
Let $B$ be a one-dimensional irreducible local scheme, 
$z\in B^{(1)}$, $\eta\in B^{(0)}$.
% $z\in B$ be the closed point, $\eta$ be the generic point.
Given a projective scheme $X$ over $B$ such that 
$\mathrm{Cl}_{X}({X_{\underline{\eta}}})=X$, where $X_{\underline{\eta}}=X\times_B \eta$,
% the scheme $X_{\underline{\eta}}=X\times_B \eta$ is dense,
the scheme $X$ is equidimensional over $B$, i.e. $\dim X_{\underline{z}}=\dim X_{\underline{\eta}}$.
\end{corollary}
\begin{proof}
Without loss of generality, we can assume that $X$ is irreducible and non-empty.
Suppose that $X\times_B z=\emptyset$, then the morphism $p\colon X\to B$ passes throw $\eta$, and since $p$ is proper, and $\eta$ equals $p(p^{-1}(\eta))$, it follows that $\eta$ is closed subscheme in $B$, 
% and consequently, $\dim B=0$.
% % $\eta$ 
% % the canonical open immersion $\eta\to B$ is proper 
% that contradicts to the fact that $B$ is one-dimensional.
that contradicts to that $\dim B=1$.
So $X\times_B z\neq\emptyset$,
and
the claim follows by \Cref{prop:OnedimBaseEquidimSch}.
\end{proof}

\section{Sections of line bundles}\label{sect-app:SerreTh}

% In this section, we summarise some results on the lifting and existence of section of line bundles used in the text.

\begin{lemma}\label{lm:semilocfinlocLinearbundle}
% For any semi-local affine scheme $Z$, any line bundle $\calL$ on $Z$ is trivial.
Any line bundle $\calL$ on a semi-local affine scheme $Z$ is trivial.
In particular, the claim holds for any finite scheme $Z$ over a local scheme $X$.
\end{lemma}
\begin{proof}
Since 
% $\calL$ is a coherent sheaf on the affine scheme 
$Z$
is
affine, 
the restriction homomorphism $\Gamma(Z,\calL)\to \Gamma(F,\calL)$ is surjective, where $F$ is the set of closed points of $Z$. Since $F$ is a union of spectra of fields there is an invertible section $f\in \Gamma(F,\calL)$. Then any section $s\in \Gamma(Z,\calL)$ such that $s\big|_F=f$ is invertible.
\end{proof}

\begin{lemma}\label{lm:b0b1binfty}
Let $X$ be a local scheme, and $Z$ be a finite scheme over $X$.
Assume the residue field at the closed point of $X$ has at least three elements.
Let $\calL_1$, $\calL_0$ be linear bundles on $Z$.
Then for any section $\delta\in \Gamma(Z,\calL_0)$, there are invertible sections $s_1,s_\infty\in \Gamma(Z,\calL_1)$, $s^+_0\in \Gamma(Z,\calL_1\otimes\calL_0^{-1})$, $s_1 = s_\infty + \delta s_0^+$.
\end{lemma}
\begin{proof}
Since $Z$ is finite over a local scheme, %spectrum of a local ring, 
$Z$ is semi-local and affine.
Denote by $F$ the union of closed points of $Z$. 

Suppose $Z=F$.
If $\delta=0$, we choose any invertible sections $s_0^+$ and $s_\infty$, 
and put $s_1=s_\infty$.
If $\delta$ is invertible, since by assumption %it follows that 
$\#k\geq 3$, 
there are invertible sections $s_0, s_1\in \Gamma(F,\calL_1)$ 
%, %$s_0\neq s_1$;
such that
% then 
the section $s_\infty:=s_1-s_0$ is %non trivial, and is 
invertible. Put $s_0^+=s_0 \delta^{-1}$.
Since
for any $\delta$, if $Z=F$,
then $Z=Z(\delta)\amalg (Z-Z(\Delta))$,
the claim 
% for $Z=F$ 
is proven.
% under the assumption $
% follows.

% For any $Z$ as in the lemma, 
Thus
by the above, 
for
any $Z$ as in the lemma,
% we choose
there are 
invertible sections 
$f_1,f_\infty\in \Gamma(F,\calL_1)$, $f^+_0\in \Gamma(F,\calL_1\otimes\calL_0^{-1})$, $f_1 = f_\infty + \delta f_0^+$.
Since $\calL_0$ and $\calL_1\otimes\calL^{-1}$ are coherent sheaves on the affine scheme,
% $Z$, 
%display
% \begin{equation}\label{eq:L0L1L0-1_ZF}\Gamma(Z,\calL_0)\twoheadrightarrow \Gamma(F,\calL_0),\quad 
% \Gamma(Z,\calL\otimes\calL^{-1})\to \Gamma(F,\calL\otimes\calL^{-1}).
% \end{equation} 
% where $F$ is the union of closed points of $Z$. 
there are sections
$s_\infty\in \Gamma(Z,\calL_1)$, $s^+_0\in \Gamma(Z,\calL_1\otimes\calL_0^{-1})$ 
such that $s_\infty\big|_F=f_\infty$, $s^+_0\big|_F=f_0$.
Define $s_1 = s_\infty + \delta s_0^+$, then $s_1\big|_F=f_1$.
Since 
$f_0$, $f_1$, $f_\infty$ are invertible
$s_0$, $s_1$, $s_\infty$ are invertible.
% Since $\calL_0$, and $\calL_1\otimes\calL^{-1}$ are coherent sheaves on the affine scheme $Z$, 
% \begin{equation}\label{eq:L0L1L0-1_ZF}\Gamma(Z,\calL_0)\twoheadrightarrow \Gamma(F,\calL_0),\quad 
% \Gamma(Z,\calL\otimes\calL^{-1})\to \Gamma(F,\calL\otimes\calL^{-1}).
% \end{equation} 
% % where $F$ is the union of closed points of $Z$. 
% Then there are sections
% $s_\infty\in \Gamma(Z,\calL_1)$, $s^+_0\in \Gamma(Z,\calL_1\otimes\calL_0^{-1})$, $s_\infty\big|_F=f_\infty$, $s_0\big|_F=f_0$.
% Define $s_1 = s_\infty + \delta s_0^+$, then $s_1\big|_F=f_1$.
% Since 
% $f_0$, $f_1$, $f_\infty$ are invertible
% $s_0$, $s_1$, $s_\infty$ are invertible.
%
% a preimage of 
% an invertible section 
% along %homomorphisms
% \eqref{eq:L0L1L0-1_ZF} 
% is invertible, the claim follows.
\end{proof}

\section{Framed correspondences and homotopies}\label{sect-app:FrHomotopies}

% In this section, we summarise some framed correspondences homotopies used in the article.
% which were known before, % in particular, in \cite{hty-inv}, \cite{Nesh-FrKMW}, \cite{DrKyllfinFrpi00}, \cite{rigid-AnanievDruzh}, \cite{FrRigidSmAffpairs}, \cite{surj-etale-exc}.

\begin{lemma}\label{lm:nvarepsilonclass}
Given a scheme $U$ over a base scheme $S$ and monic polynomial $f$ in $\calO(U)[t]$ of degree $l$, there is the equality 
% of framed correspondences
$[(Z(f),f,r)]=l_\varepsilon\in \ovZF_1(U,U)$,
where $r\colon \A^1_U\to U$ is the canonical projection.
\end{lemma}
\begin{proof}
The equality follows by the sequence
\[\begin{array}{lllll} 
[(Z(f),f,r)] &=& [(Z((t-1)^{l-1}), t(t-1)^{l-1},r)]+[(Z(t),t(t-1)^{l-1},r)]&&\\ 
&=&\phantom{|}[(Z((t-1)^{l-1}), (t-1)^{l-1},r)]+[(Z(t),(-1)^{l-1}t,r)]&&\\
&=&\phantom{|}[(l-1)_{\varepsilon}] + [\langle (-1)^{l-1}\rangle]\in \ovZF_1(U,U),
\end{array}\]
where the first equality follows because of the homotopy 
\[(Z(\widetilde f), \widetilde f, r)\in \Fr_1(U\times\A^1,U), \widetilde f= f(1-\lambda)+ t(t-1)^{l-1}\lambda\in \calO(U)[\lambda][t],\]
the second because of homotopies
\[\begin{array}{lll}
(Z((t-1)^{l-1})\times\A^1, \widetilde g_0, r)&\in& \Fr_1(U\times\A^1,U),\\ \widetilde g_0= t(t-1)^{l-1}(1-\lambda)+ (t-1)^{l-1}\lambda&\in& \calO(U)[\lambda][t],\\
(Z(t)\times\A^1, \widetilde g_1, r)&\in& \Fr_1(U\times\A^1,U), \\\widetilde g_1= t(t-1)^{l-1}(1-\lambda)+ t(-1)^{l-1}\lambda&\in& \calO(U)[\lambda][t],
\end{array}\]
and the third equality holds, since $[(Z((t-1)^{l-1}), (t-1)^{l-1},r)]=[(l-1)_{\varepsilon}]$ by the inductive assumption, and because $(Z(t),(-1)^{l-1}t,r)=\langle (-1)^{l-1}\rangle$.
\end{proof}

\begin{lemma}\label{lm:classtrivsigma}
% \[%todo%display
$
(l_1)_{\varepsilon}\circ(l_2)_{\varepsilon}=\sigma(l_1l_2)_{\varepsilon}\in \ovZF_2(U,U)
$
% \]
for any scheme $U$ and integers $l_1,l_2$.
\end{lemma}
\begin{proof}
By \Cref{lm:diagonaldeterminantframing},
% for any $l_1$,
$[(l_1)_{\varepsilon}\circ\langle (-1)^{j} \rangle]=\sigma[(-1)^{j}\cdot(l_1)_{\varepsilon}]$.
Note that % for any $l_1$,
$[(-1)^{j}\cdot(l_1)_{\varepsilon}] \rangle=[(-1)^{j l_1}\cdot(l_1)_{\varepsilon}],$
because for even $l_1$, we have
$[(-1)^{j}\cdot (l_1)_{\varepsilon}]=[(l_1)_{\varepsilon}]=[(-1)^{j l_1}\cdot(l_1)_{\varepsilon}]$,
and for odd $l_1$,
$(-1)^{j}=(-1)^{j l_1}.$
Thus
$[(l_1)_{\varepsilon}\circ\sum_{j=0}^{l_2-1}\langle (-1)^{j} \rangle]=
[\sigma[\sum_{j=0}^{l_2-1}((-1)^{j}\cdot(l_1)_{\varepsilon})]=
[\sigma[\sum_{j=0}^{l_1l_2-1}\langle (-1)^{j}\rangle]=
(l_2l_1)_{\varepsilon}.$
\end{proof}

Denote by $E_n(V)$ the subgroup of elementary matrices in $\GL_n(V)$.
\begin{lemma}\label{lm:ElemenatryMatrixesZaroSectionSupportframedhomotopy}
For a scheme $V$,
% and a matrix $G\in \GL_n(V)$, consider the framed correspondence
% $c_G=(0\times V, G(t_1,\dots,t_n), \mathrm{pr})$, where $\mathrm{pr}\colon \A^n_V\to V$ is the canonical projection.
If the matrices $G_1,G_2\in \GL_n(V)$ have the same class in $E_n(V)\backslash\GL_n(V)/E_n(V)$,
% If the matrices $G_1,G_2\in \GL_n(V)$ belong to the same orbit in $E(V)\backslash\GL_n(V)/E(V)$ with respect to the two-side action of $E_n(V)$ on $\GL_n(V)$,
then
$[\sigma^{G_1}_V]=[\sigma^{G_2}_V]\in \ovZF_n(V)$.
% $[c_{G_1}]=[c_{G_2}]\in \ovZF_n(V)$.
\end{lemma}
\begin{proof}
The claim follows because $E_n(V)$ is generated by elementary transvections $e_{v,l}(t)$, $t\in \mathcal O(V)$, $v\in \Gamma(V,\mathcal O^{\oplus n})$, $l\in \Gamma(V,(\mathcal O^{\oplus n})^\vee)$, and for each matrix $G\in \GL_n(V)$, each transvection $e_{v,l}(t)$ defines $\A^1$-homotopies $\sigma^{G e_{v,l}(t\lambda)}_V, \sigma^{e_{v,l}(t\lambda) G}_V\in \Fr_n(V\times\A^1,V)$, where $\mathcal O(V\times\A^1)=\mathcal O(V)[\lambda]$.
\end{proof}

\begin{lemma}\label{lm:diagonaldeterminantframing} 
$[\langle \lambda_1\rangle \circ\langle \lambda_2\rangle]=[\sigma\langle \lambda_1\lambda_2\rangle]\in \ovZF_2(U),$
for any scheme $U$ and $\lambda_1,\lambda_2\in \mathcal O^\times(U)$.
\end{lemma}
\begin{proof}
Consider the diagonal matrices $T_1=\langle \lambda_1,\lambda_2 \rangle$ and $T_2=\langle 1,\lambda_1\lambda_2 \rangle$
in $\GL_2(U)$,
and 
the rotation
$r_{1,2}=\left( \begin{smallmatrix}0 & 1\\ -1 & 0\end{smallmatrix} \right)$,
and elementary transvections $e_{2,1}(\alpha)$ for $\alpha\in\calO_U(U)$
in $E_2(U)$.
Since 
\[r_{1,2}e_{2,1}(-\lambda_1^{-1}) T_1 e_{1,2}(\lambda_2)e_{2,1}(-\lambda_2^{-1}) = T_2,\]
the claim follows from \Cref{lm:ElemenatryMatrixesZaroSectionSupportframedhomotopy},
because $\langle \lambda_1\rangle \circ\langle \lambda_2\rangle=c_{T_1}$, and $\sigma[\langle \lambda_1\lambda_2\rangle]=c_{T_2}$ in sense of notation from \Cref{lm:ElemenatryMatrixesZaroSectionSupportframedhomotopy}.
%
% Consider the diagonal matrices $T_1=\langle \lambda_1,\lambda_2 \rangle$ and $T_2=\langle 1,\lambda_1\lambda_2 \rangle$ in $\GL_2(U)$.
% Then
% \[r_{1,2}e_{2,1}(-\lambda_1^{-1}) T_1 e_{1,2}(\lambda_2)e_{2,1}(-\lambda_2^{-1}) = T_2,\text{ where }r_{1,2}=\left( \begin{smallmatrix}0 & 1\\ -1 & 0\end{smallmatrix} \right)\in\GL_2(U). \]
% Since $r_{1,2}\in E_2(U)$, %and the classes of $T_1$ and $T_2$ are equal in $E_2(U)\backslash\GL_2(U)/E_2(U)$,
% the claim follows from \Cref{lm:ElemenatryMatrixesZaroSectionSupportframedhomotopy},
% % Since the rotation $r_{1,2}$ is an elementary matrix, $T_1$ and $T_2$ belong to the same orbit in $\GL_2(U)$ with respect to the two-side action of elementary matrices subgroup $E_2(U)$.
% % Hence the claim follows from \Cref{lm:ElemenatryMatrixesZaroSectionSupportframedhomotopy},
% because $\langle \lambda_1\rangle \circ\langle \lambda_2\rangle=c_{T_1}$, and $\sigma[\langle \lambda_1\lambda_2\rangle]=c_{T_2}$ in sense of notation from \Cref{lm:ElemenatryMatrixesZaroSectionSupportframedhomotopy}.
\end{proof}

\begin{lemma}\label{lm:framedscheme}
% Let $X\in \Sm_K$ be a smooth affine scheme 
% % of dimension $d$ over a scheme $U$ 
% such that 
% % the tangent bundle of $X$ over $U$ is trivial.
% $T_U(X)\cong \mathbf{1}^d_X$.
% Then there is a $d$-dimensional framed correspondence $(X)^\fr=(S,\varphi,g)$ from $U$ to $X$ such that $g\big|_S\colon S\to X$ is an isomorphism, see \Cref{def:reldimFr}. \end{lemma}
Let $X\in \Sm_K$ be a smooth affine scheme, $T_U(X)\cong \mathbb{1}^d_X$. %\mathbbm
There is a $d$-dimensional framed correspondence $(X)^\fr=(S,\varphi,g)$ from $U$ to $X$ such that $g\big|_S\colon S\to X$ is an isomorphism, see \Cref{def:reldimFr}.
\end{lemma}
\begin{proof}
Since the tangent bundle of $X$ is trivial then there is a closed immersion 
$X\to \A^N_U$ 
such that the normal bundle $N_{X/\A^N_U}$ is trivial.
Then there is a vector of regular functions 
$f_{d+1},\dots, f_N\in \calO(\A^N_U)$, $d=\dim_B X_U$, such that
$Z(f_{d+1},\dots, f_{N})=X\amalg \check X$
for some $\check X$. 
By \cite[Theorem I.8]{Gru}, see also \cite{Elkiksoleqhens}, \cite[Lemma 3.11]{FrRigidSmAffpairs}, 
the canonical isomorphism $\underline{g}\colon S\to X$ admits a lifting $g\colon (\A^N_U)^h_{X}\to X$.
Thus we get the $d$-dimensional framed correspondence
% \[
$(X)^\fr=(X,f_{d+1},\dots, f_{N}, g)$.
% \]
\end{proof}

\printbibliography

@article {SS,
author = {Schmidt, Johannes and Strunk, Florian},
Title = {Stable $\A^1$-connectivity over Dedekind schemes}, 
Journal = {Annals of K-Theory},
Volume = {3}, 
Number = {2}, 
Pages = {331--367}, 
year = {2018}, 
DOI = "10.2140/akt.2018.3.331"
}

@incollection {Mor0,
author = {Morel, Fabien},
Title = {An introduction to $\mathbb{A}^1$-homotopy theory}, 
BOOKTITLE = {Contemporary developments in algebraic {$K$}-theory},
SERIES = {ICTP Lect. Notes, XV}, %ICTP Trieste, 
PAGES = {357--441},
PUBLISHER = {Abdus Salam Int. Cent. Theoret. Phys., Trieste},
YEAR = {2004},
MRCLASS = {19E08 (14F35 55-02)},
MRNUMBER = {2175638},
MRREVIEWER = {Daniel C. Isaksen},
year = {2002}
}

@article {Mor1,
author = {Morel, Fabien}, 
Title = {The stable-connectivity theorems},
Journal = {K-theory}, 
Volume = {35}, 
Number = {1}, 
Pages = {1--68}, 
year = {2005}
}

@article {Ayo06,
author = {Ayoub, Joseph},
Title = {Un contre-exemple \'a la conjecture de $\A^1$-connexit\'e de F. Morel}, 
Journal = {C. R. Math. Acad. Sci. Paris}, 
Volume = {342}, 
Number = {12}, 
Pages = {943--948}, 
year = {2006}, 
doi = "10.1016/j.crma.2006.04.017"
}

@Article{zbMATH07654583,
 Author = {Saito, Shuji},
 Title = {Reciprocity sheaves and logarithmic motives},
 FJournal = {Compositio Mathematica},
 Journal = {Compos. Math.},
 ISSN = {0010-437X},
 Volume = {159},
 Number = {2},
 Pages = {355--379},
 Year = {2023},
 Language = {English},
 DOI = {10.1112/S0010437X22007862},
 zbMATH = {7654583},
 Zbl = {1516.14048}
}

@Article{zbMATH07341096,
 Author = {Kahn, Bruno and Miyazaki, Hiroyasu and Saito, Shuji and Yamazaki, Takao},
 Title = {Motives with modulus. {I}: {Modulus} sheaves with transfers for non-proper modulus pairs},
 FJournal = {{\'E}pijournal de G{\'e}om{\'e}trie Alg{\'e}brique. EPIGA},
 Journal = {{\'E}pijournal de G{\'e}om. Alg{\'e}br., EPIGA},
 ISSN = {2491-6765},
 Volume = {5},
 Pages = {46},
 Note = {Id/No 1},
 Year = {2021},
 Language = {English},
 DOI = {10.46298/epiga.2021.volume5.5979},
 zbMATH = {7341096},
 Zbl = {1506.19002}
}

@Article{zbMATH07183731,
 Author = {Saito, Shuji},
 Title = {Purity of reciprocity sheaves},
 FJournal = {Advances in Mathematics},
 Journal = {Adv. Math.},
 ISSN = {0001-8708},
 Volume = {366},
 Pages = {70},
 Note = {Id/No 107067},
 Year = {2020},
 Language = {English},
 DOI = {10.1016/j.aim.2020.107067},
 zbMATH = {7183731},
 Zbl = {1437.19001}
}

@misc{binda2024logarithmicmotivichomotopytheory,
      title={Logarithmic motivic homotopy theory}, 
      author={Federico Binda and Doosung Park and Paul Arne Østvær},
      year={2024},
      eprint={2303.02729},
      archivePrefix={arXiv},
      primaryClass={math.AG},
      url={https://arxiv.org/abs/2303.02729}, 
}

@online {logDMk,
%% 26 Apr 2020]
%title = {Triangulated categories of logarithmic motives over a field},
%author = {Federico Binda AND Doosung Park AND Paul Arne {\O}stv{\ae}r},
%archivePrefix = "arXiv",
%Eprint = {2004.12298},%}

@Misc{Q,
 Author = {Quillen, Daniel},
 Title = {Higher algebraic {{\(K\)}}-theory. {I}},
 Year = {1973},
 Language = {English},
 HowPublished = {Algebr. {{\(K\)}}-{Theory} {I}, {Proc}. {Conf}. {Battelle} {Inst}. 1972, {Lect}. {Notes} {Math}. 341, 85-147 (1973).},
 DOI = {10.1007/BFb0067053},
 zbMATH = {3457106},
 Zbl = {0292.18004}
}

@InCollection{zbMATH05778592,
 Author = {Quillen, Daniel},
 Title = {Higher algebraic {{\(K\)}}-theory. {I}},
 BookTitle = {Cohomology of groups and algebraic \(K\)-theory. Selected papers of the international summer school on cohomology of groups and algebraic \(K\)-theory, Hangzhou, China, July 1--3, 2007},
 ISBN = {978-1-57146-144-5},
 Pages = {413--478},
 Year = {2010},
 Publisher = {Somerville, MA: International Press; Beijing: Higher Education Press},
 Language = {English},
 zbMATH = {5778592},
 Zbl = {1198.19001}
}

@Article{BlochOgus,
 Author = {S. {Bloch} and A. {Ogus}},
 Title = {{Gersten's conjecture and the homology of schemes}},
 FJournal = {{Annales Scientifiques de l'\'Ecole Normale Sup\'erieure. Quatri\`eme S\'erie}},
 Journal = {{Ann. Sci. \'Ec. Norm. Sup\'er. (4)}},
 ISSN = {0012-9593; 1873-2151/e},
 Volume = {7},
 Pages = {181--201},
 Year = {1974},
 Publisher = {Soci\'et\'e Math\'ematique de France (SMF) c/o Institut Henri Poincar\'e, Paris},
 MSC2010 = {14C35 14C15 14F99 18F25},
 Zbl = {0307.14008}
}

@misc{Sus-nonperftalk,
AUTHOR = {Andrei Suslin},
TITLE = {Motivic complexes over nonperfect fields},
YEAR = {2008},
month = {12},
howpublished = {
Talk on the second annual conference-meeting MIAN–POMI "Algebra and Algebraic Geometry",
PDMI RAS, St. Petersburg},
% howpublished = {
% Talk on the second annual conference-meeting MIAN–POMI "Algebra and Algebraic Geometry",
% St. Petersburg Department of Steklov Mathematical Institute of Russian Academy of Sciences, St. Petersburg
% },
note = {\url{http://www.mathnet.ru/php/presentation.phtml?option_lang=eng&presentid=265}, Video time: 32:15},
% note = {Video time: 32:15},
% URL = {http://www.mathnet.ru/php/presentation.phtml?option_lang=eng&presentid=265},
% annote = {annote},
}

@article {Sus17nonperfect,
AUTHOR = {Andrei Suslin},
TITLE = {Motivic complexes over nonperfect fields},
JOURNAL = {Ann. K-Theory},
VOLUME = {2}, 
NUMBER = {2},
YEAR = {2017}, 
PAGES = {277-302},
}

@book {cycles-book,
    AUTHOR = {Voevodsky, Vladimir and Suslin, Andrei and Friedlander, Eric
              M.},
     TITLE = {Cycles, transfers, and motivic homology theories},
    SERIES = {Annals of Mathematics Studies},
    VOLUME = {143},
 PUBLISHER = {Princeton University Press, Princeton, NJ},
      YEAR = {2000},
     PAGES = {vi+254},
   MRCLASS = {14F42 (14C25 19E15)},
  MRNUMBER = {1764197},
MRREVIEWER = {Spencer J. Bloch},
}

@article {Nesh-FrKMW,
    AUTHOR = {Neshitov, Alexander},
     TITLE = {Framed correspondences and the {M}ilnor-{W}itt {$K$}-theory},
   JOURNAL = {J. Inst. Math. Jussieu},
  FJOURNAL = {Journal of the Institute of Mathematics of Jussieu. JIMJ.
              Journal de l'Institut de Math\'{e}matiques de Jussieu},
    VOLUME = {17},
      YEAR = {2018},
    NUMBER = {4},
     PAGES = {823--852},
      ISSN = {1474-7480},
   MRCLASS = {14F42 (19D45)},
  MRNUMBER = {3835524},
       DOI = {10.1017/S1474748016000190},
       URL = {https://doi.org/10.1017/S1474748016000190},
}

@online{ccorrs,
% author={Andrei Druzhinin AND Håkon Kolderup},
% title={Cohomological correspondence categories},
% year={2018},
% archivePrefix = "arXiv",
% Eprint = {1808.05803},
% }

@article{framed-cancel,
author = {A. Ananyevskiy and G. Garkusha and I. Panin}, 
title = {Cancellation theorem for framed motives of algebraic varieties}, 
JOURNAL = {Advances in Mathematics}, 
volume = {383}, 
% pages = {107681}, 
NOTE = {article no. 107681}, 
year = {2021}, 
issn = {0001-8708}, 
doi = {https://doi.org/10.1016/j.aim.2021.107681}, 
% url = {https://www.sciencedirect.com/science/article/pii/S0001870821001195}, 
% keywords = {Motivic homotopy theory, Framed motives, Cancellation theorem},
% abstract = {The machinery of framed (pre)sheaves was developed by Voevodsky [17]. Based on the theory, framed motives of algebraic varieties are introduced and studied in [5]. An analog of Voevodsky's Cancellation Theorem [18] is proved in this paper for framed motives stating that a natural map of framed S1-spectraMfr(X)(n)→Hom_(G,Mfr(X)(n+1)),n⩾0, is a schemewise stable equivalence, where Mfr(X)(n) is the nth twisted framed motive of X. This result is also necessary for the proof of the main theorem of [5] computing fibrant resolutions of suspension P1-spectra ΣP1∞X+ with X a smooth algebraic variety. The Cancellation Theorem for framed motives is reduced to the Cancellation Theorem for linear framed motives stating that the natural map of complexes of abelian groupsZF(Δ•×X,Y)→ZF((Δ•×X)∧(Gm,1),Y∧(Gm,1)),X,Y∈Sm/k, is a quasi-isomorphism, where ZF(X,Y) is the group of stable linear framed correspondences in the sense of [5].}
}

@incollection {morel-trieste,
    AUTHOR = {Morel, Fabien},
     TITLE = {An introduction to {$\mathbb A^1$}-homotopy theory},
 BOOKTITLE = {Contemporary developments in algebraic {$K$}-theory},
    SERIES = {ICTP Lect. Notes, XV},
     PAGES = {357--441},
 PUBLISHER = {Abdus Salam Int. Cent. Theoret. Phys., Trieste},
      YEAR = {2004},
   MRCLASS = {19E08 (14F35 55-02)},
  MRNUMBER = {2175638},
MRREVIEWER = {Daniel C. Isaksen},
}

@online{five-authors,
%author={Elden Elmanto AND Marc Hoyois AND Adeel A. Khan AND Vladimir Sosnilo AND Maria Yakerson},
%title={Motivic infinite loop spaces},
%year={2018},
%archivePrefix = "arXiv",
%Eprint = {1711.05248},
%note = {to appear in CJM}
%}

@article{FramedGamma,
 Author = {Garkusha, Grigory A. and Panin, Ivan A. and {\O}stv{\ae}r, Paul Arne},
 Title = {Framed motivic {{\(\Gamma \)}}-spaces},
 FJournal = {Izvestiya: Mathematics},
 Journal = {Izv. Math.},
 ISSN = {1064-5632},
 Volume = {87},
 Number = {1},
 Pages = {1--28},
 Year = {2023},
 Language = {English},
 DOI = {10.4213/im9246e},
 zbMATH = {7733647}
}

@online{Panin_MovLoemmasAdz,
title = {Moving lemmas in mixed characteristic and applications},
author = {Panin, I.},
year={2022},
archivePrefix = "arXiv",
Eprint ={2202.00896v1},
}

@article{Framed,
author = {Garkusha, Grigory and Panin, Ivan},
year = {2021},
% month = {09},
pages = {261-313},
title = {Framed motives of algebraic varieties (after V. Voevodsky)},
volume = {34},
journal = {Journal of the American Mathematical Society},
doi = {10.1090/jams/958}
% year={2014 2018},
% archivePrefix = "arXiv",
% Eprint = {1409.4372},
% note = {To appear in \emph{J. Am. Math. Soc.}},
}

@article {Voe-cancel,
    AUTHOR = {Voevodsky, Vladimir},
     TITLE = {Cancellation theorem},
   JOURNAL = {Doc. Math.},
  FJOURNAL = {Documenta Mathematica},
      YEAR = {2010},
    NUMBER = {Extra vol.: Andrei A. Suslin sixtieth birthday},
     PAGES = {671--685},
      ISSN = {1431-0635},
   MRCLASS = {14F42 (19E15)},
  MRNUMBER = {2804268},
MRREVIEWER = {Oliver R\"ondigs},
}

@online{GWStrHomInv,
author={Andrei Druzhinin},
title={Strict homotopy invariance of {N}isnevich sheaves with {GW}-transfers},
%year={2018},
archivePrefix = "arXiv",
Eprint = {1709.05805},
}

@incollection {Voe-hty-inv,
    AUTHOR = {Voevodsky, Vladimir},
     TITLE = {Cohomological theory of presheaves with transfers},
 BOOKTITLE = {Cycles, transfers, and motivic homology theories},
    SERIES = {Ann. of Math. Stud.},
    VOLUME = {143},
     PAGES = {87--137},
 PUBLISHER = {Princeton Univ. Press, Princeton, NJ},
      YEAR = {2000},
   MRCLASS = {14F05},
  MRNUMBER = {1764200},
}

@incollection {Voe-motives,
    AUTHOR = {Voevodsky, Vladimir},
     TITLE = {Triangulated categories of motives over a field},
 BOOKTITLE = {Cycles, transfers, and motivic homology theories},
    SERIES = {Ann. of Math. Stud.},
    VOLUME = {143},
     PAGES = {188--238},
 PUBLISHER = {Princeton Univ. Press, Princeton, NJ},
      YEAR = {2000},
   MRCLASS = {14F42 (14C25)},
  MRNUMBER = {1764202},
}

@article {mot-functors,
    AUTHOR = {Dundas, Bj\o rn Ian and R\"ondigs, Oliver and {\O}stv{\ae}r, Paul Arne},
     TITLE = {Motivic functors},
   JOURNAL = {Doc. Math.},
  FJOURNAL = {Documenta Mathematica},
    VOLUME = {8},
      YEAR = {2003},
     PAGES = {489--525},
      ISSN = {1431-0635},
   MRCLASS = {55P42 (14F42)},
  MRNUMBER = {2029171},
MRREVIEWER = {Mark Hovey},
}

@article{hiWt,
% author = {\selectlanguage{russian}Дружинин\selectlanguage{english}/Druzhinin}, %Андрей Andrei 
% year = {2014},
% month = {01},
% pages = {181-182},
% title = {\selectlanguage{russian}О гомотопически инвариантных предпучках с Witt-трансферами\selectlanguage{english}
% /On the homotopy invariant presheaves with Witt-transfers},
% volume = {69},
% journal = {\selectlanguage{russian}Успехи математических наук\selectlanguage{english}
% /Russian Mathematical Surveys},
% doi = {10.4213/rm9593},
% }

@article{hty-inv,
author={Grigory Garkusha AND Ivan Panin},
title={Homotopy invariant presheaves with framed transfers},
journal = {Cambridge Journal of Mathematics},
year={2020},
volume = {8},
number = {1},
pages = {1--94},
}

@article{ElmantoKhannonperfect,
author = {Elmanto, Elden and Khan, Adeel},
year = {2020},
month = {01},
pages = {28-38},
title = {Perfection in motivic homotopy theory},
volume = {120},
journal = {Proceedings of the London Mathematical Society},
doi = {10.1112/plms.12280}
}

@article {Jardine-spt,
    AUTHOR = {Jardine, J. F.},
     TITLE = {Motivic symmetric spectra},
   JOURNAL = {Doc. Math.},
  FJOURNAL = {Documenta Mathematica},
    VOLUME = {5},
      YEAR = {2000},
     PAGES = {445--553 (electronic)},
      ISSN = {1431-0635},
   MRCLASS = {55P42 (14F42 55U35)},
  MRNUMBER = {1787949},
MRREVIEWER = {Jianqiang Zhao},
}

@incollection {Quillen,
    AUTHOR = {Quillen, Daniel},
     TITLE = {Higher algebraic {$K$}-theory. {I}},
 BOOKTITLE = {Algebraic {$K$}-theory, {I}: {H}igher {$K$}-theories ({P}roc.
              {C}onf., {B}attelle {M}emorial {I}nst., {S}eattle, {W}ash.,
              1972)},
     PAGES = {85--147. Lecture Notes in Math., Vol. 341},
 PUBLISHER = {Springer, Berlin},
      YEAR = {1973},
   MRCLASS = {18F25},
  MRNUMBER = {0338129},
MRREVIEWER = {Stephen M. Gersten},
}

@inproceedings {Voe98,
    AUTHOR = {Voevodsky, Vladimir},
     TITLE = {{$\bold A^1$}-homotopy theory},
 BOOKTITLE = {Proceedings of the {I}nternational {C}ongress of
              {M}athematicians, {V}ol. {I} ({B}erlin, 1998)},
   JOURNAL = {Doc. Math.},
  FJOURNAL = {Documenta Mathematica},
      YEAR = {1998},
    NUMBER = {Extra Vol. I},
     PAGES = {579--604 (electronic)},
      ISSN = {1431-0635},
   MRCLASS = {14F35 (14A15 55U35)},
  MRNUMBER = {1648048 (99j:14018)},
MRREVIEWER = {Mark Hovey},
}

@article {Morel-Voevodsky,
    AUTHOR = {Morel, Fabien and Voevodsky, Vladimir},
     TITLE = {{${\mathbb A}^1$}-homotopy theory of schemes},
   JOURNAL = {Inst. Hautes \'Etudes Sci. Publ. Math.},
  FJOURNAL = {Institut des Hautes \'Etudes Scientifiques. Publications
              Math\'ematiques},
    NUMBER = {90},
      YEAR = {1999},
     PAGES = {45--143 (2001)},
      ISSN = {0073-8301},
     CODEN = {PMIHA6},
   MRCLASS = {14F35 (19E08)},
  MRNUMBER = {1813224 (2002f:14029)},
MRREVIEWER = {Marc Levine},
       URL = {http://www.numdam.org/item?id=PMIHES_1999__90__45_0},
}

@article{Gru,
author = {L.~Gruson}, 
title = {Une propri\'et\'e des couples hens\'eliens},
%  Une propriété des couples henséliens
% Gruson, Laurent
Journal = {Publications math\'ematiques et informatique de Rennes}, 
volume = {4}, 
year = {1972}, 
% number = 10, 
% 13 p. 
}

@article{Elkiksoleqhens,
%author = {R.~Elkik}, 
%title = {Solutions d’\'equations \`a coefficients dans un anneau h\'enselien},
%Journal = {Annales scientifiques de l’\'E.N.S.}, 
%volume = {6},
%year = {1973}, 
%pages = {553-603},
%}

@online{Voe-notes,
%author={Voevodsky, Vladimir},
%title={Notes on framed correspondences},
%year={2001},
%% url = {},
%}

@book {Hartshorne-AlG,
    AUTHOR = {Hartshorne, Robin},
     TITLE = {Algebraic geometry},
      NOTE = {Graduate Texts in Mathematics, No. 52},
 PUBLISHER = {Springer-Verlag},
   ADDRESS = {New York},
      YEAR = {1977},
     PAGES = {xvi+496},
      ISBN = {0-387-90244-9},
   MRCLASS = {14-01},
  MRNUMBER = {0463157 (57 \#3116)},
MRREVIEWER = {Robert Speiser},
}

@online{surj-etale-exc,
author={Andrei Druzhinin AND Ivan Panin},
title={Surjectivity of the étale excision map for homotopy invariant framed presheaves},
year={2018},
archivePrefix = "arXiv",
Eprint = {1808.07765},
}

@article {BigFrmotives,
author={Grigory Garkusha AND Ivan Panin}, 
Title={The triangulated categories of framed bispectra and framed motives},
Journal={ St. Petersbg. Math. J.},
number = {6},
volume ={34}, 
pages = {991--1017}, 
year = {2023},
}

@article {ConeTheGNP,
Author={Grigory Garkusha AND Alexander Neshitov AND Ivan Panin},
Title={Framed motives of relative motivic spheres},
JOURNAL = {Trans. Amer. Math. Soc.}, 
Volume = {374}, 
Year = {2021}, 
pages = {5131-5161},
DOI = {https://doi.org/10.1090/tran/8386}
% MSC (2020): Primary 14F42, 55P42
% archivePrefix = "arXiv",
% Eprint = {1604.02732},
}

@online{SmModelSpectrumTP,
author={Druzhinin, A.},
Title={Geometric models for fibrant resolutions of motivic suspension spectra},
archivePrefix = "arXiv",
Eprint = {1811.11086},
}

@Article{MotiveswithmodulusIII,
 Author = {Kahn, Bruno and Miyazaki, Hiroyasu and Saito, Shuji and Yamazaki, Takao},
 Title = {Motives with modulus. {III}: {The} categories of motives},
 FJournal = {Annals of \(K\)-Theory},
 Journal = {Ann. \(K\)-Theory},
 ISSN = {2379-1683},
 Volume = {7},
 Number = {1},
 Pages = {119--178},
 Year = {2022},
 Language = {English},
 DOI = {10.2140/akt.2022.7.119},
 zbMATH = {7547325},
 Zbl = {1505.19003}
}

@Article{zbMATH07341097,
 Author = {Kahn, Bruno and Miyazaki, Hiroyasu and Saito, Shuji and Yamazaki, Takao},
 Title = {Motives with modulus. {II}: {Modulus} sheaves with transfers for proper modulus pairs},
 FJournal = {{\'E}pijournal de G{\'e}om{\'e}trie Alg{\'e}brique. EPIGA},
 Journal = {{\'E}pijournal de G{\'e}om. Alg{\'e}br., EPIGA},
 ISSN = {2491-6765},
 Volume = {5},
 Pages = {31},
 Note = {Id/No 2},
 Year = {2021},
 Language = {English},
 DOI = {10.46298/epiga.2021.volume5.5980},
 zbMATH = {7341097},
 Zbl = {1468.19005}
}

@online{FrRigidSmAffpairs,
%author={Druzhinin, A},
%Title={Rigidity for smooth affine pairs over a field}, 
%archivePrefix = "arXiv",
%Eprint = {1809.04158},
%}

@online{DKO:SHISpecZ,
%author={Druzhinin, Andrei AND Kolderup, H{\aa}kon  AND {\O}stv{\ae}r, Paul Arne}, %H{\aa}kon Kolderup Paul Arne %{\O}stv{\ae}r
%Title={Strict $\mathbb{A}^{1}$-invariance over the integers}, %Strict homotopy invariance over integers
%archivePrefix = "arXiv",
%Year = {2020},
%month = {12},
%Eprint = {2012.07365v1},
%}

@online{DrKyllfinFrpi00,
author={Andrei Druzhinin AND Jonas Irgens Kylling},
title={Framed motives and zeroth stable motivic homotopy groups in odd characteristic},
year={2018},
archivePrefix = "arXiv",
Eprint = {1809.03238},
}

@Article{GNThomSpectra,
 Author = {Garkusha, Grigory and Neshitov, Alexander},
 Title = {Fibrant resolutions for motivic {Thom} spectra},
 FJournal = {Annals of \(K\)-Theory},
 Journal = {Ann. \(K\)-Theory},
 ISSN = {2379-1683},
 Volume = {8},
 Number = {3},
 Pages = {421--488},
 Year = {2023},
 Language = {English},
 DOI = {10.2140/akt.2023.8.421},
 zbMATH = {7734998},
 Zbl = {1527.14048}
}

@online {SHzar,
title = {Zariski-local framed $\A^1$-homotopy theory},
author = {Andrei Druzhinin AND Vladimir Sosnilo},
archivePrefix = "arxiv",
Eprint = {2108.08257}
}

@article{ConnDodekindDomains,
Author = {Schmidt, J. and Strunk, F.},
Title = {Stable $\mathbb{A}^{1}$-connectivity over Dedekind schemes},
Journal = {Ann. K-Theory},
Volume = {3}, Number ={2}, Year = {2018}, Pages = {331-367},
}

@Article{ConnBase,
 Author = {Druzhinin, Anderi Eduardovich},
 Title = {Stable {{\(\mathbb{A}^1\)}}-connectivity over a base},
 FJournal = {Journal f{\"u}r die Reine und Angewandte Mathematik},
 Journal = {J. Reine Angew. Math.},
 ISSN = {0075-4102},
 Volume = {792},
 Pages = {61--91},
 Year = {2022},
 Language = {English},
 DOI = {10.1515/crelle-2022-0048},
 zbMATH = {7612787},
 Zbl = {1510.14018}
}

@online{ColumnsCousinlocesssmXoverB,
%author={Druzhinin Andrei AND Kolderup H{\aa}kon  AND {\O}stv{\ae}r Paul Arne}, 
%Title={Cousin complexes in motivic homotopy theory},
%archivePrefix = "arXiv",Eprint = {2402.10541},}
\end{document}